\documentclass[12pt,oneside]{amsbook}
\usepackage{amsmath, amssymb}
\usepackage{yfonts}
\usepackage[T1]{fontenc}
\usepackage{setspace}
\usepackage{verbatim} %allows for block commenting via \begin{comment}...\end{comment} 
\newtheorem{thm}{Theorem}
\newtheorem{prop}[thm]{Proposition}

\newtheorem{defn}[thm]{Definition}
\newtheorem{lem}[thm]{Lemma}
\newtheorem{cor}[thm]{Corollary}

\newtheorem{remark}[thm]{Remark}

\numberwithin{section}{chapter} \numberwithin{equation}{chapter}
\numberwithin{subsection}{chapter} \numberwithin{thm}{chapter}
\newcommand{\parab}{u_t = -\sum_{i,j=1}^m X^*_i \big( a_{ij}(x,t)X_j u \big)}

\newcommand{\ellip}{\lambda^{-1} |\xi|^2 \leq a_{ij}(x,t) \xi_i \xi_j \leq \lambda |\xi|^2 }
\newcommand{\normbr}{\frac{1}{|B_{r}|}}
\newcommand{\normbrp}{\frac{1}{|B_{r '}|}}
\newcommand{\normbrt}{\frac{1}{|B_{r} \times (0, \tau)|}}

\newcommand{\hg}{\nabla _{0}}
\newcommand{\nlparab}{ u_t = -\sum ^{m}_{j=1} X_j^*  A_{j}( x, \hg u)   }

\newcommand{\partiali}{\partial_{x_i}}
\newcommand{\partialj}{\partial_{x_j}}

 \allowdisplaybreaks[1]

\begin{document}
\pagenumbering{roman}
\thispagestyle{empty}
\hspace{5pt}\\ \vspace{100pt}
\begin{center}
A HARNACK INEQUALITY \\AND H\"OLDER CONTINUITY FOR \\WEAK SOLUTIONS TO PARABOLIC OPERATORS \\INVOLVING H\"ORMANDER VECTOR FIELDS 
\end{center}
\newpage
\thispagestyle{empty}
\begin{center}
A HARNACK INEQUALITY \\AND H\"OLDER CONTINUITY FOR \\WEAK SOLUTIONS TO \\ PARABOLIC OPERATORS \\INVOLVING H\"ORMANDER VECTOR FIELDS 
\\ \vspace{80pt}
\normalsize
A dissertation submitted in partial fulfillment\\
of the requirements for the degree of \\
Doctor of Philosophy in Mathematics\\ \vspace{75pt}
By \\ \vspace{75pt}
Garrett James Rea \\
University of Arkansas, 1997 \\
Bachelor of Arts in Computer Science \\
University of Arkansas, 2007 \\ 
Master of Science in Mathematics \\
\vspace{50pt}
May 2010\\
University of Arkansas
\end{center}
\newpage
\thispagestyle{empty}
\begin{center}
\textbf{Abstract}\end{center}
This paper deals with two separate but related results. First we consider weak solutions to a parabolic operator with H\"ormander vector fields. Adapting the iteration scheme of J\"urgen Moser for elliptic and parabolic equations in $\mathbb{R}^n$ we show a parabolic Harnack inequality. Then, after proving the Harnack inequality for weak solutions to equations of the form $u_t = \sum X_i (a_{ij} X_j u)$ we use this to show H\"older continuity. We assume the coefficients are bounded and elliptic. The iteration scheme is a tool that may be adapted to many settings and we extend this to nonlinear parabolic equations of the form $u_t = -X_i^* A_j(X_j u)$. With this we show both a Harnack inequality and H\"older continuity of weak solutions. 
\newpage
\thispagestyle{empty}
\singlespace 
\noindent This dissertation is approved for\\
Recommendation to the \\
Graduate Council\\ \vspace{25pt}

\noindent Dissertation Director:\\
\vspace{30pt}

\noindent \underline{\hspace{300pt}}\\
\textbf{Luca Capogna, Ph.D.} \\ \vspace{25pt}

\noindent Dissertation Committee:\\ \vspace{30pt}

\noindent \underline{\hspace{300pt}}\\
\textbf{Andrew Raich, Ph.D.} \\ \vspace{30pt}

\noindent \underline{\hspace{300pt}}\\
\textbf{Loredana Lanzani, Ph.D.} \\ \vspace{30pt}

\newpage
\thispagestyle{empty}
\begin{center}
\textbf{Dissertation Duplication Release}\end{center}
\vspace{30pt}
I hereby authorize the University of Arkansas Libraries to duplicate this dissertation when needed for research and/or scholarship.\\
\vspace{30pt}

\noindent Agreed \hspace{10pt} \underline{\hspace{300pt}} 

\noindent \hspace{55pt} \textbf{Garrett J. Rea}

\newpage
\setcounter{page}{5}
\doublespace
\begin{center}
\textbf{Acknowledgements}\end{center}

This dissertation would not have been possible without my advisor, Dr. Luca Capogna. From the time I took his Partial Differential Equations, my first graduate mathematics class, to stepping through calculations in my dissertation, he has been patient, encouraging and has pushed me to do more than I thought I could. Over the last few years, he has given me fantastic guidance and made himself available to meet or discuss far above and beyond what I could expect. I sincerely thank him and appreciate everything he has done. 

I would like to thank my committee, Dr. Loredana Lanzani for her help and Dr. Andy Raich for valuable input on my dissertation. 

My wife Brooke has been the greatest support to me from start to finish. I could not have done this without her. She says that she would do this again.

\newpage
\begin{center} 
\textbf{Contents}\end{center}

%\indent 0.1 \hspace{5pt} Notation \hfill 1 \\ 
\indent 0.1  \hspace{5pt} Introduction \hfill 1 \\ 
\indent 0.2. \hspace{5pt} Main results \hfill 5 \\ 
\indent 0.3. \hspace{5pt} Examples of a nonlinear parabolic operator \hfill 8 \\ 
\indent 0.4. \hspace{5pt} Chronology of the Elliptic and Parabolic Harnack inequality \hfill  10\\ 
\indent 0.5. \hspace{5pt} H\"ormander Vector Fields \hfill 19 \\ 
\indent 0.6. \hspace{5pt} Carnot Carath\'eodory Distance and Balls\hfill 21 \\ 
\indent 0.7. \hspace{5pt} Weak solutions \hfill 24 \\ 
Chapter 1. \hspace{5pt} Parabolic Harnack Inequality \hfill 27 \\
\indent 1.1 \hspace{5pt} Flow of the proof \hfill 27 \\
\indent 1.2 \hspace{5pt} Cacciopoli Inequality \hfill 27 \\
\indent 1.3 \hspace{5pt} Cacciopoli for other values of $p$ \hfill 33 \\
\indent 1.4 \hspace{5pt} Sobolev Theorem and Iteration \hfill 36 \\
\indent 1.5 \hspace{5pt} Bridge step \hfill 49 \\
\indent 1.6 \hspace{5pt} Lag mapping \hfill 52 \\
\indent 1.7 \hspace{5pt} Estimates on log(u) \hfill 55 \\
\indent 1.8 \hspace{5pt} H\"older Continuity \hfill 64 \\
Chapter 2. \hspace{3pt} Nonlinear Generalization \hfill 68 \\
\indent 2.1. \hspace{5pt} Introduction \hfill 68 \\
\indent 2.2. \hspace{5pt} Nonlinear Cacciopoli Inequality \hfill 70 \\
\indent 2.3. \hspace{5pt} Nonlinear Cacciopoli for $u^p$ \hfill 75 \\
\indent 2.4. \hspace{5pt} Nonlinear Iteration and bridge step \hfill 79 \\
\indent 2.5. \hspace{5pt} Nonlinear Estimates on log(u) and H\"older continuity \hfill 91 \\
Bibliography \hfill 91
\newpage

%---------------------------------------
%\begin{center} 
%\textbf{0.2 Notation}
%\end{center}
\begin{tabular}{l   p{11.3cm}}
Notation       & Use \\
\hline
$\mathbb{R}^n$               &  the $n$-dimensional Euclidean space \\
$\mathbb{R}^{n+1}$ & the product of the $n$-dimensional Euclidean space and a time interval $(0,T)$\\ 
$x$ & a point in $n$-dimensional Euclidean space. $x = \{ x_1, x_2,...x_n \}^{T}$\\
$d_{cc}(x,y)$ & the Carnot-Carath\'eodory (CC) distance from $x$ to $y$. (see pg. 20)\\
$B_{cc}(x,r)$ & a CC ball centered at $x$ with radius $r$. $B_{cc}(x,r) = \{ y \in \mathbb{R}^n \ | \ d_{cc}(x,y) < r \}$. Also $B(x,r), \ B_r$. (see pg. 20) \\
$|B_{cc}(x,r)| $ & the volume of the ball. $|B_{cc}(x,r)|= \int_B dx$. Here $dx$ is Lebesgue measure\\
$R$ & an open bounded rectangle in $\mathbb{R}^{n+1}$. It is the product of a Carnot-Carath\'eodory ball and a time interval: $B_{cc}(x,r) \times (0,T)$\\
$u(x,t) $ & a solution to a partial differential equation. $u(x,t):\mathbb{R}^{n+1} \rightarrow \mathbb{R}$\\
$X_i$ & a vector field.\\
$\{ X_i \}_{i=1}^{m}$ & a collection of vector fields.\\
$X_i u$ & the derivative of $u$ in the direction of $X_i$.\\
$X^*_i$ & the formal adjoint of $X_i$. (see pg. 18)\\
$\nabla u$ & the spatial gradient of $u$. \\
$|\nabla u|$ & the norm of the spatial gradient of $u$. $|\nabla u| = [\sum_{i=1}^n (\partial  u / \partial x_i)^2 ]^{1/2}$ \\
$\hg u$ & the horizontal gradient of $u$.  $\hg u = \sum_{i=1}^m X_i u X_i$. (see pg. 19)\\
$|\hg u|$ & the norm of the horizontal gradient of $u$. $|\hg u| = [\sum_{i=1}^m (X_i u)^2 ]^{1/2}$ \\
$\{ a_{ij}(x,t) \}$ & an $m \times m$ matrix of whose coefficients are functions of $x$ and $t$. \\ 
$L^2(R)$ & the Lebesgue space of functions on $R$. It consists of functions $u$ for which $\int_R u^2 dxdt < \infty$. (see pg. 23)\\
$L^{\infty}(R)$ & consists of all functions $u$ for which $ess \, \sup \, u < \infty$\\
\end{tabular} \newpage
\begin{tabular}{l   p{11.2cm}}
Notation       & Use \\
\hline
$W^{1,2}(R)$ & the Sobolev space of functions on $R$. It consists of functions $u \in L^2(R) $ with the weak derivatives of $u$ in $L^2(R)$.  \\
$S^{1,2}(R)$ & the Sub-elliptic Sobolev space of functions on $R$. It consists of functions $u \in L^2(R) $ with the weak derivatives of $u$ in $L^2(R)$. (see pg. 23)\\
$u_t$ & the derivative of $u$ in time. $u_t = \partial u / \partial t$\\

%$\pi^{\pi}$         & 3646   \\
%$(\pi^{\pi})^{\pi}$ & 806627 \\
\end{tabular} \\

\mainmatter %Used to the amsbook knows this is the main portion and will use arabic numerals
%*** You need the next 5 lines so that the page numbers are always in the upper right corner.  The blank page is so that
% the page whose number would be at the bottom center is skipped.  The \setcounter is so that the page numbering starts at 1
% on the second page, which is really the first page :).  You just don't include the blank page with your final turned in version.
\newpage
\thispagestyle{empty}
\hspace{5pt}
\pagebreak
\setcounter{page}{0}
%%%%%%%%%%%%%%%%%%%%%%%%%%%%%%%%%%%%%%%%%%%%%%%%%%%%%%%%%%%%%%%%%%%%%%%%%%%%%%%%%%
% You cannot use the \section{Introduction} command because it will put a header at the top of every page that the grad school
% doesn't allow.  

\begin{center} \textbf{0.1 Introduction}\end{center} \vspace{3pt}
This paper deals with two versions of the heat equation. The classical heat equation tells us the temperature at a given location in a material at a given time. It is simply modeled by
$$u_t = \Delta u$$
in some region. The derivatives here are the usual ones taken along the coordinates axes. This equation has been studied by many people and is very well known. This equation may also be written as:
$$u_t = \hbox{div} \, (\nabla u).$$

A typical property of interest is the bounds on the behavior of a solution to such a PDE. One classical result is that solutions to the heat equation have a Harnack-type inequality which says that the maximum temperature in a domain $R$ can be known in terms of the minimum. That is
$$\underset{R}{\max} \; u(\cdot, t_1) \leq C \underset{R}{\min} \; u(\cdot, t_2).$$
In this inequality, the constant $C$ can be calculated explicitly based on the first time observed, $t_1$, the second time observed, $t_2$, and the dimension of the space in which the heat diffuses. A second property that is often of interest is how smooth or well-behaved the solutions are. One goal is to show that solutions are H\"older continuous with respect to some metric. If $x$ and $y$ are any two points in the domain considered, then
$$|u(x) - u(y)| \leq C |x-y|^{\alpha}$$
where $\alpha$ is between zero and one.  

It is natural to extend the restrictions that one places on the equation to gain further understanding of the nature of the diffusion of heat. We may add the complication of specifying the nature of the medium through which the heat diffuses. We do this with a matrix that models the homogeneity of the material. The new equation would be
$$u_t = \hbox{div} \, (a \nabla u)$$
where $a = \{ a_{ij} \}$ is an $n \times n$ matrix of functions that model the material. This equation is well known when we assume the coefficients are smooth. This would be the easiest case. Even in the worst case, where the functions are merely measurable, much is known. 

A natural extension of the original problem would then be to take the derivatives not along the coordinates axes, but rather along vector fields. The vector fields could yield strange curves and even fractal-type behavior. In general, we consider less vector fields, $m$ than the dimension of the space, $n$. But in order to do much with them we ask that they possess a certain property: that they can span the entire space being considered when taken in any combination. If $X_i$ is any of the $m$ vector fields, then the equation takes the form
$$u_t = \sum_{ij}^{m} X_j (a_{ij} X_i u).$$ 
Because the solutions cannot be assumed to be too smooth, we add another complication. We assume that the derivatives do not necessarily exist as functions but only in a weak sense. This further muddies the original problem by making it more difficult to take the two derivatives we have in the above equations.      

A final wrinkle on the original equation would be to consider not a matrix $a_{ij}$, but a non-linear function $A$ in the equation. This gives us
$$u_t = \sum_{j}^{m} X^*_j (A_j(x, \hg u)).$$
Non-linear behavior is often difficult to deal with in PDE. This last change would give us a much wider class of functions to consider. 

This paper proves the two properties considered above regarding the behavior of solutions. We show a Harnack inequality and H\"older continuity of the solutions. We assume that the coefficients are only measurable, bounded functions and not smooth. We assume the derivatives are taken along vector fields that yield potentially strange geometry. 

This paper deals with two separate but related problems. First we will consider positive, weak solutions, $u(x,t) \in S^{1,2}(R)$, to $$\parab$$ in a domain $R = B_{cc}(x,r) \times (0,\tau) \subset \mathbb{R}^{n+1}$ where the $\{ X_i \} ^{m}_{i=1}$, $m \leq n$, are H\"ormander vector fields. The matrix of coefficients $\{ a_{ij}(x,t) \}$ are bounded, measurable functions on $R$ satisfying an ellipticity condition. We will call this the parabolic problem. For this problem we deal with two issues, establishing a Harnack-type inequality for positive, weak solutions and  proving H\"older continuity for those solutions. The problem is a classical one addressed by many authors. Jurgen Moser \cite{Mo1} dealt with weak solutions, $u(x,t) \in W^{1,2}(Q)$, to $$\sum_{i,j=1}^{n} {\partial \over \partial x_i} \big(a_{ij}(x,t) {\partial \over \partial x_j } v \big) =0$$ in a unit cube, $Q \subset \mathbb{R}^n$. In this paper, he introduced an iteration scheme that has formed the foundation of an approach to the issue of regularity of solutions adapted by many authors. His general approach was to define for some bounded domain $D \subset \mathbb{R}^n$,
$$M(p,D) := \big({1 \over |D|} \int_{D} u^p \big)^{1/p}.$$ Using the facts that 
$$\underset{D}{\max}\, u = M(\infty, D)$$
and
$$\underset{D}{\min}\, u = M(- \infty, D)$$
combined with Cacciopoli and Sobolev inequalities he was able to construct an iteration that accomplishes for some constant $C > 0$
$$\underset{D}{\max}\, u = M(\infty, D) = \underset{k \rightarrow \infty}{\lim} \, \big({1 \over |D|} \int_{D} u^{pk} \big)^{1/pk} \ \leq C \ \big({1 \over |D|} \int_{D} u^{pk} \big)^{1/pk} $$ along with
$$\underset{D}{\min}\, u = M(-\infty, D) = \underset{k \rightarrow \infty}{\lim} \, \big({1 \over |D|} \int_{D} u^{-pk} \big)^{1/-pk} \ \geq \ C \ \big({1 \over |D|} \int_{D} u^{-pk} \big)^{1/-pk} . $$  
These two inequalities hold for values of $pk$ close to $0$: $-pk < 0 < pk$. Then these are linked with 
$$ \big({1 \over |D|} \int_{D} u^{pk} \big)^{1/pk} \leq C \,  \big({1 \over |D|} \int_{D} u^{-pk} \big)^{1/-pk} $$
for $pk \rightarrow \epsilon$. 
Moser later \cite{Mo2} addressed the continuity of solutions to a parabolic operator in $\mathbb{R}^n$ using the same iteration. We will employ the Moser iteration here. 

The second problem we deal with is positive, weak solutions to $$\nlparab$$ in the same region $R$ as above. Here, $\{X_i \}_{i=1}^{m}$, $m \leq n$, are smooth vector fields satisfying H\"ormander's condition for hypoellipticity. $X_j^*$ is the formal adjoint of $X_j$, and $A(x, \hg u)$ is a measurable function satisfying certain requirements (see (\ref{nlheat}) and following). This is a nonlinear generalization of the parabolic problem. It is useful to extend the Moser technique to many situations as authors have since his work in 1961. Trudinger \cite{Tr1} used the Moser method in 1967 to establish Harnack type inequalities for quasilinear elliptic equations of the form
$$\hbox{div} \,   \vec{a} (x,u, \nabla u) + b(x,u,\nabla u) = 0$$ and
$$\hbox{div} \,   \vec{a} (x,t, u, \nabla u) + b(x,t,u,\nabla u) - u_t= 0.$$
Then in 1995 Saloff-Coste, \cite{SC} used a sketch version of Moser's method to show that the parabolic Harnack inequality is equivalent to both the doubling property for measures and the Poincar\'e inequality in the Riemannian manifold setting. This small list is far from exhaustive. 

We note that the following results are motivated by a direct proof, using the iteration of Moser, of the Harnack inequality and H\"older continuity of weak solutions to the two equations we examine. In the case of derivatives along H\"ormander vector fields, results are known in the case of certain elliptic operators, (see \cite{CDG3}) but not in the parabolic case. Also, results are not known for the case where the matrix of coefficients $a_{ij}$ are only measurable and $L^{\infty}$ bounded. Results have been known for parabolic-type operators since 1958 with constant $a_{ij}$ coefficients (see \cite{Had}). Then for measurable coefficients a parabolic and elliptic Harnack inequality was proven for non-divergence form operators (see \cite{KrS}). This current work brings together the following: derivatives along H\"ormander vector fields, merely bounded and measurable coefficients, parabolic-type operators and Moser's method. To accomplish this new tools are needed. A Sobolev inequality for H\"ormander vector fields from Capogna, Danielli, Garofalo, (see \cite{CDG3}), a Poincar\'e inequality for H\"ormander vectors fields from Jerison, (see \cite{Jer}) and the existence of smooth test functions of compact support on the Carnot-Carath\'eodory balls with bounded horizontal gradient from Citti, Garofalo and Lanconelli, (see \cite{CiL}). The iteration scheme of Moser is necessarily different. We take into account the geometry of the Carnot-Carath\'eodory balls that makes the integration tricky. These all result in explicit constants for both the Harnack inequality and H\"older continuity in this geometry.  
  
%-------------------------------------------------------------
% You cannot use the \chapter{} or \section{} commands because they put a header that the grad school doesn't allow
\begin{center} 
\textbf{0.2 Main results}
\end{center}
%-------------------------------------------------------------
   
The main results of this paper are the following four theorems. For the parabolic case we have the 

\begin{thm}\label{PHI} \emph{Harnack inequality for weak solutions }

\noindent Let $u \in S^{1,2}(R)$ be a positive, weak solution to 
\begin{equation}
\label{heat} u_t = -\sum_{i,j=1}^{m} X^*_i \big(a_{ij}(x,t) X_j u \big)
\end{equation}
 in $R = \{B_{cc}(x, r ) \times (0, \tau) \}$. Also define $R^{+} = \{B(x, r ' ) \times (\tau ^+, \tau ) \}$ and $R^{-} = \{B(x, r ' ) \times (\tau _1^-, \tau _2^- ) \}$ for $0 <  r ' < r$, $0 < \tau_1^- < \tau_2^-  < \tau^+  < \tau$. The $\{X_i \}_{i=1}^{m}$, $m \leq n$, are smooth H\"ormander vector fields and the matrix of coefficients $a_{ij}(x,t)$ are $L^{\infty}(R)$ bounded and elliptic satisfying
\begin{equation}
\label{ellipticity_condition}\lambda^{-1} |\xi|^2 \leq a_{ij}(x,t) \xi_i \xi_j \leq \lambda |\xi|^2 
\end{equation}
for every $(x,t) \in R$ and for every $\xi \in \mathbb{R}^n$. Then there exists a constant $C>0$ depending on $\lambda > 0$, $n$, and $R$ such that 
\begin{equation}
\label{harnack_inequality} \underset{R^-}{\max} \, u \,  \leq C \, \underset{R^+}{\min} \, u
\end{equation}
where the maximum and minimum are the essential supremum and infimum over their respective domains.
\end{thm} 
\noindent We use the notation for some function $v$, 
$$X^*_j v= -\sum_{i=1}^{n}(\frac{\partial b_i}{\partial x_i} v + \frac{\partial v}{\partial x_i} b_i). \ \hbox{(See section 0.5)}$$ 
The second is the
\begin{thm}\label{PHC} \emph{H\"older continuity of weak solutions with respect to the CC metric}

Let $u(x,t) \in S^{1,2}(R)$ be a positive, weak solution to $$\parab$$ in $R$ as defined above. Then for every $(x,t)$ and $(y,s) \in R  = B(x,r) \times (0, \tau)$, $0 < r < 1$ there exists a constant $C$ depending on $n$, $\lambda,$ and $R$ such that 
\begin{equation}
\label{holder_cty}|u(x,t) - u(y,s)| \leq \, C \, d\big((x,t),(y,s) \big)^{\alpha}
\end{equation}
for some $0 < \alpha <1$ and $d\big((x,t),(0,0)\big) = \max \, \{ d_{cc}(x,0), \sqrt{|t|} \}$.
\end{thm} 
\noindent Then for the nonlinear generalization we have the following
\begin{thm}\label{NLHI} \emph{Harnack inequality for weak solutions to a nonlinear operator}

Let $u \in S^{1,2}(R)$ be a positive, weak solution to 
\begin{equation}
\label{nlheat}u_t = -\sum_{j=1}^{m} X_j^* A_{j}(  x, \hg u)
\end{equation}
in $R = \{B(x, r ) \times (0, \tau) \}$. Also define $R^+ = \{B(x, r ' ) \times (\tau ^+, \tau ) \}$ for $0 < r' <r$, $0<\tau _1^-<\tau _2^-<\tau ^+<\tau$. $\{X_i \}_{i=1}^{m}$, $m \leq n$, are smooth vector fields satisfying H\"ormander's condition for hypoellipticity. $X_j^*$ is the formal adjoint of $X_j$. The measurable function $A( x, \hg u)$ fulfills the following structural requirements:
$$
\noindent (S) \quad   \quad \quad  \quad \left\{
 \begin{array}{rl}
  A_j( x, \hg u) \leq C|\hg u|\\
  \lambda |\hg u|^2 \leq A_j(x, \hg u) X_j \leq \Lambda |\hg u|^2\\
 \end{array} \right. 
 $$
%\end{equation*}
for some constants $C$, $\Lambda$, $\lambda > 0$. Then there exists a constant $C$ depending on $\lambda$, $n$, and $R$ such that 
\begin{equation}
 \underset{R^-}{\max} \, u \,  \leq C \, \underset{R^+}{\min} \, u
\end{equation}
where the max and min are the essential supremum and infimum over these domains.
\end{thm} 
\noindent Last we have 
\begin{thm}\label{NLHC} \emph{H\"older continuity of weak solutions to a nonlinear operator with respect to the CC metric}

\noindent Let $u(x,t)\in S^{1,2}(R)$ be a positive, weak solution to $$u_t =- \sum_{j=1}^{m} X_j^* A_{j}(  x, \hg u)$$ in $R$ as defined above. Then for every $(x,t)$ and $(y,s) \in R  = B(x,r) \times (0,\tau)$, $0 < r < 1$ there exists a constant $C>0$ depending on $n$, $\lambda,$ and $R$ such that 
\begin{equation}
 |u(x,t) - u(y,s)| \leq \, C \, d \big((x,t),(y,s)\big)^{\alpha}
\end{equation}
for some $0 < \alpha < 1$ and $d \big((x,t),(0,0) \big) = \max \, \{ d_{cc}(x,0), \sqrt{|t|} \}$.
\end{thm}  

%-------------------------------------------------------------
% You cannot use the \chapter{} or \section{} commands because they put a header that the grad school doesn't allow
\begin{center} 
\textbf{0.3 Examples of the nonlinear parabolic operator}
\end{center}
 
Let $u \in S^{1,2}(R)$ be a positive, weak solution to 
\begin{equation}
u_t = -X_j^* A_{j}(  x, \hg u)
\end{equation}
in $R = \{B(x, r ) \times (0, \tau) \}$, $R \subset \mathbb{R}^n$. Also define $R^+ = \{B(x, r ' ) \times (\tau ^+, \tau ) \}$ for $0 < r' <r$, $0<\tau _1^-<\tau _2^-<\tau ^+<\tau$. $\{X_i \}_{i=1}^{m}$, $m \leq n$, are smooth vector fields satisfying H\"ormander's condition for hypoellipticity. That is:
$$\hbox{Rank \big(Lie[} X_1, ... , X_m\hbox{] \big)}(x) = n$$
for any $x \in \mathbb{R}^n$. $X_j^*$ is the formal adjoint of $X_j$. The measurable function $A( x, \hg u): \mathbb{R}^m \rightarrow \mathbb{R}^m$ fulfills the following structural requirements:
$$
\noindent (S) \quad   \quad \quad  \quad \left\{
 \begin{array}{rl}
  A_j( x, \hg u) \leq C|\hg u|\\
  \lambda |\hg u|^2 \leq \sum_{j=1}^m A_j(x, \hg u) X_j u \leq \Lambda |\hg u|^2\\
 \end{array} \right. 
 $$
%\end{equation*}
for some constants $C$, $\Lambda$, $\lambda > 0$.
\begin{enumerate}
\item First we point out that the parabolic version we prove in Theorem \ref{PHI} is contained in this nonlinear example.
Let $A_j( x, \hg u) = a_{ij}(x,t) X_i u$. Then we have 
$$u_t = -\sum_j X_j^* A_{j}(  x, \hg u) = -\sum_{i,j}X^*_j \big(a_{ij}(x,t) X_i u \big).$$
Then we have the structural requirements $(S)$ met in the inequality (\ref{ellipticity_condition}) from the ellipticity and boundedness of $a_{ij}(x,t)$:
$$\ellip$$
and noticing that 
$$a_{ij}(x,t) X_i u  \leq \lambda |X_i u|.$$ 
We can apply the results of Theorems \ref{NLHI} and \ref{NLHC} to achieve the same Harnack inequality and H\"older continuity proved directly in the parabolic case.

\item Let $$A_j( x, \hg u) = \frac{X_j u}{\sqrt{1 + | \hg u|^2}}$$
with $|\hg u| < M$ in $R$ for some constant $M > 0$. This is similar to the mean curvature for intrinsic minimal graphs studied by several authors: Cheng, Hwang, Malchiodi and Yang \cite{ChH}, Manfredini \cite{Ma} and Pauls \cite{Pa}. We consider weak solutions of the equation
$$u_t = \sum_j^{m} X_j \big(\frac{X_j u}{\sqrt{1 + | \hg u|^2}}\big).$$
The structural requirements are met in the following:
$$A_j( x, \hg u) = \ \frac{X_j u}{\sqrt{1 + | \hg u|^2}} \leq \  \frac{|\hg u|}{\sqrt{1}},$$
$$\vec{A}( x, \hg u)\cdot \hg u  = \ \frac{X_j u}{\sqrt{1 + | \hg u|^2}} \cdot \hg u \leq  \frac{|\hg u|^2}{\sqrt{1}}$$
and 
$$\vec{A}( x, \hg u)\cdot \hg u  = \ \frac{\hg u}{\sqrt{1 + | \hg u|^2}} \cdot \hg u \geq  \frac{|\hg u|^2}{\sqrt{1 + M^2}}.$$
%$$\frac{|\hg u|^2}{|\hg u|^2} \leq A_j(\hg u)\cdot \hg u.$$
Applying the results from Theorem \ref{NLHI}, we have a Harnack inequality for weak  solutions to (\ref{nlheat}) in $R$. We also are guaranteed H\"older continuity of these solutions with respect to the CC metric. 

\item Let $$A_j( x, \hg u) = \frac{X_j u}{1 + | \hg u|^2}$$ 
with $|\hg u| < M$  in $R$ for some constant $M > 0$. That gives us 
$$u_t = -\sum_j X_j^* A_{j}(  x, X_j u) = -\sum_{j} X_j^* \frac{X_j u}{1 + | \hg u|^2} .$$
We can see that the structural conditions are satisfied in the following:
$$A_j( x, \hg u) = \ \frac{X_j u}{1 + | \hg u|^2} \leq \  \frac{| \hg u|}{1 + |  \hg u|^2} \ \leq |\hg u|,$$

$$\vec{A}( x, \hg u)\cdot \hg u = \sum_{j=1}^m A_j( x, \hg u)\cdot X_j u  = \frac{\hg u}{1 + | \hg u|^2} \cdot \hg u  \leq  \ \frac{|\hg u|^2}{1} ,$$
and 
$$\vec{A}( x, \hg u)\cdot \hg u = \sum_{j=1}^m A_j( x, \hg u)\cdot X_j u  = \frac{\hg u}{1 + | \hg u|^2} \cdot \hg u  \geq  \ \frac{|\hg u|^2}{1 + | M|^2}.$$
Applying the results from Theorem \ref{NLHI}, we have a Harnack inequality for weak  solutions to (\ref{nlheat}) in $R$. We also are guaranteed H\"older continuity of these solutions with respect to the CC metric. 

\end{enumerate}

%-------------------------------------------------------------
% You cannot use the \chapter{} or \section{} commands because they put a header that the grad school doesn't allow
\begin{center} 
\textbf{0.4 Chronology of the elliptic and parabolic Harnack inequality}
\end{center}
%-------------------------------------------------------------

This is a table showing some developments relevant to the current work in the study of elliptic and parabolic operators. We will look at different settings and with assumptions on the matrix of coefficients, $A(x,t)$, having components $ a_{ij}(x,t) $. 

{ \bf 1954 - J. Hadamard, \cite{Had}} shows the first parabolic Harnack inequality for operators with constant coefficients $a_{ij}$ in $\mathbb{R}^n$. Pini obtained the same result in 1954 separately.

{ \bf 1958 - J. Nash, \cite{Nas}} proves a Harnack inequality for an elliptic operator in $\mathbb{R}^n$ using heat kernel estimates. This was for measurable and uniformly elliptic $\{ a_{ij} \}$. Moser showed the same in 1964 using the iteration method that is employed in the current paper.

{ \bf 1961 - J. Moser, \cite{Mo2}} proves an elliptic Harnack inequality for solutions to 
$$\sum_{i,j=1}^n{\partial \over \partial x_i} (a_{ij} {\partial \over \partial x_j}\ u)= 0$$
and that the Harnack inequality implies H\"older continuity.

{ \bf 1964 - J. Moser, \cite{Mo1}} shows a parabolic Harnack inequality using his iteration scheme and the John-Nirenberg theorem. This paper provides the basis for many authors ( see \cite{FrL}, \cite{Tr1}, \cite{SC} ) as well as the current results.

{ \bf 1967 - Aronson and Serrin, \cite{AS}} treat the second order quasilinear parabolic equation
$$u_t = \hbox{div} \, A(x,t,u,\nabla u) = B(x,t,u,\nabla u)$$ 
and were the first to show two-sided Gaussian bounds for divergence form uniformly elliptic operators in $\mathbb{R}^n$. 

{ \bf 1967 - Trudinger, \cite{Tr1}} deals with two kinds of equations
$$\hbox{div} \, {\bf a} (x,u, \nabla u) + b(x,u,\nabla u) = 0$$ and
$$\hbox{div} \, {\bf a} (x,t, u, \nabla u) + b(x,t,u,\nabla u) - u_t= 0$$
and following the iteration used by Moser establishes a Harnack inequality for weak solutions to both. 

{ \bf 1968 - Trudinger, \cite{Tr2}} shows point wise estimates for quasilinear parabolic equations. 

{ \bf 1969 - Bony, \cite{Bo}} considers degenerate elliptic partial differential equations of the type
$$\sum^r_{i=1} X^2_i u + Yu + cu = 0$$
where $X_1,...,X_r, Y$ are first order homogeneous differential operators. Various properties of solutions are studied in connection with the Lie algebra generated by these operators. He proves a local form of the Harnack inequality.

{ \bf 1971 - Moser, \cite{Mo3}} shows a parabolic Harnack inequality using alternate methods besides the John-Nirenberg theorem. The John-Nirenberg theorem is used in the current results in a form modified by Aimar \cite{Ai}. Moser uses an approach by Bombieri based on minimal surfaces. 

{ \bf 1973 - Trudinger, \cite{Tr3}}  deals with operators of the type
$$Lu = - {\partial \over \partial x_i} \big(a_{ij}(x)u_{x_i} + a_i(x)u \big) + b_i(x)u_{x_i} + a(x)u = -{{\partial f^i \over \partial x_i}} + f$$
where the $a_{ij}$ are measurable functions on a domain $\Omega \subset \mathbb{R}^n$ and $f$ is an integrable function on $\Omega$. He proves the Harnack inequality for solutions and then shows H\"older continuity.
 
{ \bf 1980 - Krylov and Safonov, \cite{KrS}} prove a parabolic and elliptic Harnack inequality for non-divergence form operators $$L u = \ - \sum^{n}_{i,j=1} a_{ij} \partial _i \partial _j u$$ with $a_{ij}$ measurable, symmetric and uniformly elliptic.

{ \bf 1983 - Franchi and Lanconelli, \cite{FrL}} were the first to apply a Moser technique to obtain H\"older continuity of solutions for non-uniformly elliptic operators with measurable coefficients. 

{ \bf 1986 - Jerison and Sanchez-Calle, \cite{JeS}} treat degenerate elliptic operators with smooth coefficients. We give their main result.
\begin{thm}
Let $\Omega$ be an open subset of $\mathbb{R}^n$. 
Let $$L = \sum_{i,j=1}^{n} \big(1/k(x)\big) {\partial \over \partial x_i} \big(k(x) a_{ij} {\partial \over \partialj}\big)$$ where the coefficients $a_{ij}$ and $k$ are $C^{\infty} (\overline{\Omega})$, $k$ is positive and the matrix $a_{ij}$ is symmetric positive semi-definite for every $x \in \overline{\Omega} $. Also suppose that $L$ satisfies a sub-elliptic estimate: There is a constant $C$ and a number $\epsilon$ such that all $u \in C^{\infty}_0 (\Omega)$ satisfy
$$|| u ||_{2 \epsilon} \leq C(||Lu|| + ||u||).$$
Here $$||u||_{s} = \big(\int |\hat{u}(\xi)|^2 (1+ |\xi|^2)^s \, d \xi\big)^{1/2},$$ 
$$||u|| = \big(\int |u(x)|^2 \, d \mu (x)\big)^{1/2} $$
and  
$$d \mu (x) = h(x) dx.$$ 
Then for $({\partial \over \partial _t} - L)u(t,x) = 0$ we have for the heat kernel $h(t,x,y)$ the following:
$$h(t,x,y) \leq A \mu \big(B(x, t^{1/2})\big)^{-1} \, e^{-d(x,y)^2 / \gamma t}$$
$$h(t,x,y) \geq A' \mu \big(B(x,t^{1/2})\big)^{-1} \, e^{-d(x,y)^2 / \gamma ' t}$$
$$|\partial ^k _t X_{j_1} \cdots X_{j_p} h(t,x,y)| \leq C_{p,k} t^{-k-p/2} \mu \big(B(x, t^{1/2})\big)^{-1} \, e^{-d(x,y)^2 / \gamma t}$$
For some positive constants $A$, $A'$, $\gamma$, $\gamma '$, $C_{p,k}$ and every $x,y \in \Omega$ $0< t < 1$.
\end{thm} 
 
{ \bf 1988 - Kusuoka and Stroock, \cite{KS}} use probability methods for heat kernel estimates on H\"ormander vector fields. We give their main result.
\begin{thm}
Let $a_{ij} \in C^2_b$ be a symmetric, non-negative definite matrix-valued function. Define $$Lu(x) = \sum ^{N}_{i,j=1} [\partiali (a_{ij} \partialj u)](x).$$
Let $P(t,x, \cdot)$ be the unique transition probability function on $\mathbb{R}^N$ such that the associated Markov semigroup $\{ P_t : t > 0 \}$ satisfies 
$$P_t \phi (x) - \phi (x) = \int_0^1 [P_s L \phi ](x)ds$$ for all $ \phi \in C^{\infty}_0 (\mathbb{R}^N)$. Let $ \{ E_{ \lambda } : \lambda \in [0, \infty)  \}$ denote the resolution of the identity determined by $\{ \overline{P_t} : t > 0 \}$ and set $A = \int_{[0, \infty)} \lambda dE_{ \lambda }$. Let the Dirichlet form be given by 
$$S(f,f) = \int_{ [0, \infty ) } \lambda d (E) $$
Assume there exist $A \in (0, \infty)$, $\nu \in (0, \infty)$ and $\delta \in (0, \infty)$ such that 
$$||f||^{2+4/ \nu}_{2} \leq A (S(f,f) + \delta ||f||^2_2)||f||^{4/ \nu}_{1}$$
for $f \in L^2 (\mathbb{R}^N);$ or equivalently that there is a $B \in (0, \infty)$ such that 
$$||P_t||_{1 \rightarrow \infty } \leq Be^{\delta t } / t^{\nu / 2}$$
for $t > 0.$ Then $P(t, x, dy) = p(t, x,y )dy$ and there is a $C \in (0, \infty)$ depending only on $\nu$ such that for each $\rho \in (0,1]$ and all $(t,x) \in (0, \infty) \times \mathbb{R}^N:$
$$p(t,x, \cdot) \leq C(A/ \rho t )^{\nu/2} e^{\rho \delta t} exp[-D(x, \cdot)^2 / (1 + \rho )t] \ \ a.e. $$
\end{thm}

{ \bf 1992 - Varopolous, \cite{VSC}} shows various parabolic Harnack inequalities using heat kernel estimates. He does not use the approach taken by Moser.
\begin{thm}
Let $\{ X_i \}_{i=1}^k$ be a H\"ormander system of vector fields on the manifold $V$, $Y$ a $C^{\infty}$ vector field on $V$, and a, a $C^{\infty}$ function on $V$. Let $D$ be the differential operator
$$\sum_{i=1}^k X_i ^2 + Y - a - {\partial \over \partial t},$$
defined on $\mathbb{R} \times V.$  Let $I =  (\alpha,\beta) \subset \mathbb{R}$, $\Omega$ and open relatively compact subset of $V$, $K$ a compact subset of $\Omega$, $t_1$, $t_2$, such that $\alpha < t_1 < t_2 < \beta$, $J \in {\it I}(N)$, $m \in \mathbb{N}$.
Then there exists a constant $C$ such that every positive solution $u$ of $Du = 0$ in $I \times \Omega$ satisfies
$$\underset{x \in K}{sup} \ |({\partial \over \partial t})^m ({\partial \over \partial x})^J \ u(t_1, x)| \leq C \, \underset{x \in K}{inf} \ u (t_2, x)$$
with spatial derivatives taken with respect to some fixed system of local coordinates.
\end{thm}
 
{ \bf 1992 - Grigory'an, \cite{Gr}} uses heat kernel estimates to prove the parabolic Harnack inequality and does not rely on Moser's method. He shows the Harnack inequality is implied by the doubling property and the Poincar\'e inequality. He also shows a Harnack inequality without using the doubling property and Poincar\'e inequality. The behavior of the Green function $G(x, y, t)$ of the Cauchy problem for the heat equation on a connected, noncompact, complete Riemannian manifold is investigated. For manifolds with boundary it is assumed that the Green function satisfies a Neumann condition on the boundary. His results are similar to Saloff-Coste \cite{SC} but arrives at them with different methods.
\begin{thm}
Let $M$ be a geodesically complete, noncompact, smooth, connected Riemannian manifold of dimension $n$. Let $\Delta$ be the Laplace operator on $M$ associated with the Riemannian metric $g$. Suppose $M$ satisfies the following two hypotheses:

1. For any concentric geodesic balls $B(x,R)$, $B(x,2R)$ of radii $R$ and $2R$
$$\mu B(x, 2R) \leq A \mu (x,R)$$
Where the constant $A$ is the same for all balls.

2. For some constant $N >1$ and $a > 0$, for any geodesic ball $B(x,R) \subset M$ and for any smooth function $f$ defined in the ball $B(x, NR)$ the following inequality holds:
$$\int_{B(x, NR)} |\nabla f|^2 \, \geq \, {a \over R^2} \underset{\xi \in \mathbb{R}}{inf} \ \int_{B(x,R)} (f - \xi)^2$$
Then the uniform parabolic Harnack inequality is valid with the constant $P$ depending on $A$, $a$, $N$.

Fix a point $x \in M$. Let u(x,t) be a positive solution of $(\partial _t - \Delta) u = 0$ in $Q_{8R}=B(x,R) \times (0, 8R^2)$ which is smooth in $\overline{Q_{8R}}$ and satisfies the Neumann condition for $x \in \partial M$ (if $\partial M$ is non-empty). Set $\widetilde{Q} = B_R \times (3R^2, 4R^2)$, and suppose that $\underset{\widetilde{Q}}{sup} \, u = 1$. Then $u(x, 64 R^2) \geq \gamma$, $\gamma = \gamma (A, a, N)$. 
\end{thm}

{ \bf 1993 - Citti, Garofalo, and Lanconelli, \cite{CiL}  } establish a uniform Harnack inequality for weak solutions $u$ of 
$$
(-\sum_{i=1}^{m} X_i ^2 + V) \ u \ = \ 0
$$
where $V$ is a measurable function belonging to an analogue to the local Kato class for the $X_i$. They also introduce a very important cutoff function for CC balls which we make use of in the current paper.  We state this lemma for completeness.
\begin{lem}{Existence of cut-off functions}

For every ball $B(x,t) \subset R$, and every $s$,$t >0$ with $s < t$, there exists a function $\psi \in C^{\infty}_0 (B(x,t))$ such that $0 \leq \psi \leq 1$ in $B(x,t)$, $\psi \equiv 1$ in $B(x,s)$ and $|\hg \psi| \leq C/(t-s)$, where $C$ is a positive constant which is independent of $t$ and $s$.
\end{lem}

{ \bf 1993 - Capogna, Danielli and Garofalo, \cite{CDG3} } prove a version of the Sobolev embedding theorem for functions in $S^{1,p}_0$ using the fractional integral operators used by Gilbarg and Trudinger \cite{GiT}. With this result they prove a Harnack inequality for solutions to a subelliptic operator and H\"older continuity. Their main result:
\begin{thm}
Let $U \subset \mathbb{R}^n$ be a bounded open set with relative homogeneous dimension $Q > 0$. Let $\{ X_i \}_{i=1}^m$ be $C^{\infty}$ vector fields satisfying H\"ormander's condition for hypoellipticity. Assume $1 < p \leq Q$. Let $u \in S^{1,p}_{loc}$ be a non-negative solution to the equation
$$\sum^m_{j=1} X^*_j A_j (x, u, \hg u) = f(x,u, \hg u)$$
Then there exist $C>0$ and $R_0 > 0$ such that for any $B_R = B(x,R)$ with $B(x, 4R) \subset U$, and $R \leq R_0$,
$$\underset{B_R}{ess \ sup} \ u \ \leq \ C(\underset{B_R}{ess \ inf} \ u \ + K(R))$$
\end{thm}

{ \bf 1993 - Varopolous, \cite{VSC}} shows that the upper and lower Gaussian bounds on the fundamental solution are implied by the parabolic Harnack principle. This is done in the setting of H\"ormander vector fields.

{ \bf 1995 -  Saloff-Coste, \cite{SC}} shows the parabolic Harnack inequality is equivalent to the doubling property and the Poincar\'e inequality. This is done in Riemannian manifold setting. He assumes the $a_{ij}$ coefficients are bounded, measurable, symmetric and elliptic. His main result:
\begin{thm}
Let M be a $C^{\infty}$ connected manifold. Let L be a second order differential operator with real $C^{\infty}$ coefficients on M and no zero order term. Use as distance $\rho (x,y) = sup \{ |f(x) - f(y)|, f \in C^{\infty} (M), |\nabla f| \leq 1 \}$. Consider the two properties, the doubling property and Poincar\'e inequality:
\begin{equation}
 \label{doubling_property}\mu (B(x,2r)) \leq C_1 \mu (B(x,r)) 
\end{equation}
\begin{equation}
\label{poincare_property} \int_{B(x,r)} |f - f_B|^2 d \mu \leq C_2 r^2 \int_{B(x,2r)} |\nabla f|^2 d \mu 
\end{equation} 
for $0 < r < r_0$, $x \in M$, $f \in C^{\infty}(M)$ and $f_B = {1 \over |B|} \int_{B} f$
 
Then the following two statements are equivalent

A. The properties (\ref{doubling_property}) and (\ref{poincare_property}) above hold for some $ r_0$.

B. There exists $r_1 > 0$ and a constant $C$ depending only on the parameters 
$0 < \epsilon < \eta < \delta < 1$ such that for any $x\in M$, any real $s$, and any $0 < r < r_1$, any non-negative solution $u$ of $(\partial_t + L)u = 0$ in $Q = (s-r^2, s) \times B(x,r)$ satisfies
\begin{equation}
\underset{Q_{-}}{sup} \ u \ \leq \ C \underset{Q_{+}}{inf} \ u
\end{equation}
where $$Q_{-} = [s-\delta r^2, s - \eta r^2 ] \times B(x, \delta r)$$ and $$Q_{+} = [s-\epsilon r^2, s) \times B(x, \delta r)$$
\end{thm}  
It should be noted the similarities and differences between this paper's current results and those of Saloff-Coste. The current results are contained in those of Saloff-Coste. We examine a parabolic-type operators that satisfy properties (\ref{doubling_property}) and (\ref{poincare_property}) above. So by his results, the Harnack inequality holds in some region. He does not employ a direct proof which we do here. Also, using the geometry created by the H\"ormander vector fields, we calculate the constants for the Harnack inequality and H\"older continuity explicitly. A further advantage of the direct proof is to allow for extension into the generalized nonlinear case as we do here as well as other cases. 

{% \bf 1998 - Capogna, L., Garofalo, N. Regularity of quasi linear equations on the Heisenberg group. Comm. Pure Appl. Math, (1998)}

%Let $\mathbb{H}^n$ be the n-dimensional Heisenberg group. Let $X_1$, ... $X_{2n}$, $S$ be the usual basis of left invariant vector fields with $S$ in the center and for $\Omega \subset \mathbb{H}^n$ and open set, $k \in \mathbb{N}, 1 \leq a < \infty$ set $S^{k,n}_{loc}(\Omega)$ be the Sobolev space of $L^a_{loc}(\Omega)$ vector functions $u = (u_1, ... u_N) \in \mathbb{R}^N$, whose horizontal derivatives of oder less or equal to $k$ are in $L^a_{loc}(\Omega)$. Consider weak solutions to 
%$$u^{\alpha}_{t} - \sum_{\beta = 1}^{N} \sum^{i,j = 1}_{2n} -X_i(A^{\alpha \beta }_{ij} X_j u^{\beta} + f^{\alpha}_i (x)) = f^{\alpha} (x)$$ 
%for $\alpha = 1,...,N$, where $A^{\alpha \beta }_{ij}$ is a constant matrix satisfying 
%$$\sum_{\beta = 1}^{N} \sum^{i,j = 1}_{2n} A^{\alpha \beta }_{ij} \xi^{\alpha}_{i} \xi^{\beta}_{j} \geq \lambda |\xi|^2$$
%for every $\xi \in \mathbb{R}^{2n \times N}$, and $f^{\alpha}_{i}$, $f^{\alpha}$ are $L^2$ functions defined on $\Omega$. 

%The authors prove the analogue of he classical $L^2$ estimates for constant coefficients parabolic systems involving a family of H\"ormander type vector fields in a strified Nilpotent Lie group of step two. This result is used to prove partial regularity of solutions of nonlinear parabolic systems and estimates of the Hausdorff dimension of the singular set.

{ \bf 2002 - Uguzzoni, Bonfiglioli, and Lanconelli, \cite{UL1}} use heat kernel estimates with Gaussian bounds to prove the parabolic Harnack inequality on Carnot groups. 
 
{ \bf 2004 - Uguzzoni, Bonfiglioli, and Lanconelli, \cite{UL2}} construct the fundamental solutions for non-divergence form operators where the $X_i$ are H\"ormander vector fields generating a stratified group $G$ and $a_{ij}$ is a positive-definite matrix with H\"older continuous entries. They also provide Gaussian estimates of $\Gamma$ and its derivatives and some results for the relevant Cauchy problem. 

{ \bf 2006 - M. Bramanti et al., \cite{Br}} prove the existence of a fundamental solution for a class of H\"ormander heat-type operators. For this solution and its derivatives they obtain sharp Gaussian bounds that allow to prove an invariant Harnack inequality. The matrix $a_{ij} $ is real symmetric, uniformly positive definite and having H\"older continuous entries.

%-------------------------------------------------------------
% You cannot use the \chapter{} or \section{} commands because they put a header that the grad school doesn't allow
%\begin{center} 
%\textbf{0.4.1 H\"ormander Vector Field Setting }
%\end{center}
%-------------------------------------------------------------
 
A list of some developments relevant to the current work of those in the H\"ormander vector field setting. Chronologically:

{ \bf 1969 -} Bony showed local Harnack inequalities for degenerate elliptic operators. 

{ \bf 1988 -} Kusuoka and Stroock use probability methods for heat kernel estimates. 

{ \bf 1992 -} Grigory'an shows a Harnack inequality without using the doubling property of measures and the Poincar\'e inequality. 

{ \bf 1993 -} Capogna, Danielli and Garofalo prove a Harnack inequality  and H\"older continuity for solutions to a subelliptic operator.

{ \bf 1993 -} Varopolous shows that Gaussian bounds on the heat kernel are implied by the parabolic Harnack principle.

{ \bf 1995 -} Saloff-Coste shows the parabolic Harnack inequality is equivalent to the doubling property and the Poincar\'e inequality. 

{ \bf 2002 -} Uguzzoni uses heat kernel estimates with Gaussian bounds to prove a parabolic Harnack inequality. 
 
{ \bf 2004 -} Uguzzoni constructs the heat kernel for the non-divergence form operators. 

{ \bf 2006 -} Bramanti uses sharp Gaussian bounds to prove a Harnack inequality. 
%-------------------------------------------------------------
% You cannot use the \chapter{} or \section{} commands because they put a header that the grad school doesn't allow
\begin{center} 
\textbf{0.5 H\"ormander Vector Fields }
\end{center}
%-------------------------------------------------------------
 
We will use vector fields of the form $$X_i = b_i {\partial \over \partial x_i}$$
where the $b_i \in C^{\infty}(\mathbb{R}^n)$. We will consider solutions to $$\parab$$ and $$u_t = -\sum ^{m}_{j=1} X_j^* A_{j}( x, \hg u)$$ along smooth vector fields $\{ X_i \}_{i=1}^{m}$, $m \leq n$, satisfying H\"ormander's \cite{Ho} condition:
$$rank \big(\ Lie[X_1,...,X_m] \big)= n$$
at every $x \in \mathbb{R}^n$.   

We may take the derivative along a vector field. If we consider a simple example:
$$X_1 = b_1 \frac{\partial}{ \partial x_1}$$
then we have for some function $f$:
$$X_1 f = b_1 \frac{\partial f}{ \partial x_1}.$$
We also will make use of the formal adjoint of any $X_i$, $X^*_i$. The formal adjoint of a vector field $X_i$ is the operator $X^*_i$ such that for functions $f,g \in C^{\infty}_0(R)$
$$\int_R X_i f \cdot g \, dxdt = \int_R f \cdot X^*_i g \, dxdt.$$
We compute:
$$\int_R X_i f \cdot g \, dxdt = \int_R \sum_{j=1}^n (b_j \frac{\partial}{\partial x_j} f) \cdot g \, dxdt = \int_R \sum_{j=1}^n ( \frac{\partial}{\partial x_j} f) \cdot (b_j g) \, dxdt$$
And integrating by parts gives us:
$$= -\int_R  f \cdot \sum_{j=1}^n ( \frac{\partial}{\partial x_j}(b_j g) \, dxdt=-\int_R  f \cdot \sum_{j=1}^n ( \frac{\partial b_j}{\partial x_j} g + \frac{\partial g}{\partial x_j} b_j) \, dxdt $$
$$=\int_R  f \cdot X^*_i g \, dxdt$$
So that $X^*_i g$ takes the form $(-1)( \frac{\partial b_j}{\partial x_j} g + \frac{\partial g}{\partial x_j} b_j)$

In this paper, we use equations of the form
$$u_t =- \sum ^{m}_{j=1} X_j^* A_{j}( x, \hg u)$$
So for $\phi \in C^{\infty}_0(R)$ we can say 
$$\phi u_t = -\sum ^{m}_{j=1} X_j^* A_{j}( x, \hg u) \phi$$ 
And integrating over $R$
$$\int_R \phi u_t \, dxdt = -\sum ^{m}_{j=1} \int_R X_j \phi A_{j}( x, \hg u) \, dxdt.$$

We use throughout this paper the horizontal gradient of a function $f$. 
\begin{defn}{Horizontal gradient}\label{def_hg}

\noindent The horizontal gradient of a function $f$ is defined to be:
$$\hg f = X_1 f X_1 + X_2 f X_2 + \cdots + X_m f X_m.$$
\end{defn}
\noindent And the size of the horizontal gradient is
$$|\hg f| = \big(\sum_{i=1}^m (X_i f )^2\big)^{1/2}.$$ We note that for the $\{ X_i \}_{i=1}^{m}$ the horizontal gradient will not in general equal the Euclidean spatial gradient in $\mathbb{R}^n$, $\nabla f$. The vector fields $\{ X_i \}_{i=1}^{m}$ give us a natural distance. 
%-------------------------------------------------------------
% You cannot use the \chapter{} or \section{} commands because they put a header that the grad school doesn't allow
\begin{center} 
\textbf{0.6 Carnot Carath\'eodory Distance and Balls}
\end{center}
%-------------------------------------------------------------

The vector fields $\{ X_i \}_{i=1}^{m}$ give us a control distance. Consider a piecewise $C^1$ curve $\gamma : [0,T] \rightarrow \mathbb{R}^n$. Let $x, \, y \in \mathbb{R}^n$.  For $\delta > 0$ we define the class $C(\delta)$ of absolutely continuous paths $\gamma: [0, 1] \rightarrow \mathbb{R}^3$ with endpoints $\gamma(0) = x$ and $\gamma(1) = y$, so that $$\gamma ' (t) = \sum^{m}_{i=1} a_i(t)X_i |_{\gamma(t)} $$ and $$\sum^{m}_{i=1} a_i(t)^2 \leq \delta ^2$$ for a.e. $t \in [0, 1]$. Such paths are called horizontal. 
\begin{defn}{Carnot-Carath\'eodory metric}

We define the Carnot-Carath\'eodory (CC) metric to be 
\begin{equation}
d_{cc}(x, y) = \hbox{inf} \, \{ \delta \hbox{ such that }  C(\delta) = \emptyset \}.
\end{equation}
\end{defn}
\noindent A dual formulation is 
$$d_{cc}(x, y) = \hbox{ inf} \ \{ T: \exists \gamma:[0,T] \rightarrow \mathbb{R}^n, \gamma(0) = x, \gamma(T) = y, \sum^{m}_{i=1} a_i(t)^2 \leq 1  \ \ a.e.  \}.$$ We will refer to the CC distance as $d$ and $d_{cc}$. With the CC distance we can define balls in $\mathbb{R}^n$. We set 
$$B_r = B(x,r) = B_{cc}(x,r) = \{ y : d_{cc}(x,y) < r \}$$the ball centered at $x$ of radius $r$. We use as the volume of the ball $|B(x,r)| = \int_{B(x,r)} dx$. It is worth noting that the geometry created by the vector fields, $\{X_i \}_{i=1}^m$ is different from the Euclidean space, $\mathbb{R}^n$. One major issue we face in obtaining estimates on the behavior of functions is the geometry of the Carnot-Carath\'{e}odory balls. In the Euclidean setting the ratio of the volume of two balls is proportional to their radii. In the Carnot-Carath\'{e}odory setting the situation is different. The volume of two balls is related in the following way from a proposition of Nagel, Stein and Wainger \cite{NaW}.
\begin{prop}
Given a bounded set $U \subset \mathbb{R}^n$ there exist $Q \geq n$, $R_0 > 0$, $C_B > 0$, such that for every $x \in U$, $R \leq R_0$ and $0 < t < 1$
\begin{equation}
|B(x, tR)| \geq C_B t^Q |B(x, R)|
\end{equation}
\end{prop}
\noindent which we will make use of in our results. 
Test functions play an integral role in these results. Since the Sobolev theorem we use is assumed to be for functions of compact support on a Carnot-Carath\'eordory ball and our solutions are not supported in these domains, we work around the problem with certain test functions. Moser used such functions in the Euclidean case \cite{Mo1}, \cite{Mo2} and they are standard in applications of this type. But the geometry of the Carnot-Carath\'eordory balls is much different. In the Euclidean context, balls in $n$-dimensions are very well behaved, having symmetry as well as level sets which are also balls. This is not the case in the present context. The Carnot-Carath\'eordory balls can be badly behaved, having cusps and the level sets could potentially behave just as badly. When we use the test functions, we set two radii and construct the function to be $\phi = 1$ on the smaller of the concentric balls. Then the function decreases to $0$ with bounded gradient until it is $\phi =0$ outside of the larger radius. This fact has been an obstacle until 1993 when Citti, Garofalo, and Lanconelli \cite{CiL} constructed such test functions. Also Nagel and Stein \cite{Nag} in 2001 revisited the test functions along with showing the existence of smooth metric equivalent to the control metric over H\"ormander vector fields.

There is another difficulty present in the use of these test functions in the Carnot-Carath\'eodory setting. In the construction, the two radii may not be arbitrarily close to each other. We may not merely choose two radii $r > r' > 0$ and construct a useful test function. The radii must be bounded away from each other. We recall Nagel and Stein [Lemma 3.1.1] \cite{Nag} that there exist two constants $0 < c_1 < 1 < c_2$ such that $\phi = 1$ on $x < c_1 r$ and $\phi = 0$ on  $x > c_2 r$. Here $c_1$ and $c_2$ come about in the construction of level sets of the CC balls.

The facts of this construction are worth mentioning even though the details of these restrictions are not obvious when we use the test functions. In constructing the test functions Citti, Garofalo and Lanconelli construct a quasi-distance $D(x,y)$ that satisfies 
$$c_1 d_{cc}(x,y) \leq D(x,y) \leq c_2 d_{cc}(x,y)$$
for every $x,y \in B(x,r)$. Now construct a new ball $\overline{B}(x,R) = \{ y \; | \; D(x,y) < R \}$. This gives us
$$B_{\frac{R}{c_2}} \subset \overline{B}(x,R) \subset B_{\frac{R}{c_1}}$$
which says that we can create a $\overline{B}$ ball whose boundary is contained in a CC ball of a larger radius and containing a CC ball of a smaller radius. In order for us to use this fact we must construct nested balls. In the proof we use radii $r$ and $r'$ so let those constants be given. This gives us two inequalities based on
$$c_1d_{cc}(x,y) \leq D(x,y) \leq c_2 d_{cc}(x,y)$$ 
which are
$$c_1r \leq r \leq c_2 r$$
And the second is
$$c_1r' \leq r' \leq c_2 r'.$$
So we need balls that satisfy
$$B(x,r') \subset \overline{B}(x,r'c_2) \subset \overline{B}(x,rc_1)\subset B(x,r).$$
This requires our constants to satisfy $r'c_2 \leq rc_1$. Recalling that $0 < c_1 < 1 < c_2$ forces $r$ and $r'$ to be bounded away from each other and not arbitrarily close.

%-------------------------------------------------------------
% You cannot use the \chapter{} or \section{} commands because they put a header that the grad school doesn't allow
\begin{center} 
\textbf{0.7 Weak solutions}
\end{center}
%-------------------------------------------------------------
We will consider positive, weak solutions to $$\parab$$  and $$\nlparab$$ in a region $R \subset \mathbb{R}^{n+1}$. In order to create the context for these solutions we will begin with definitions.
\begin{defn}{$L^2(R)$}

$L^2(R)$ is defined to be the space of all functions such that 
$$||f||_{L^2(R)} = \big(\int_{R} |f|^2 \, dxdt \big)^{1/2} < \infty$$
\end{defn}
\begin{defn}{$L^2_{loc}(R)$}\\

$L^2_{loc}(R)$ is defined to be the space of all functions such that 
$$||f||_{L^2(R)} = (\int_{K} |f|^2 \, dxdt)^{1/2} < \infty$$
for all compact subsets $K$ of $R$.
\end{defn}
\noindent Throughout the paper we use the notion of weak derivatives. We define
\begin{defn}{Weak derivative}

If $u, \, v \in L^2(R)$, then $v$ is the weak derivative of $u$ if
$$\int_R u X^*_i \phi \, dxdt = - \int_R v \phi \, dxdt $$
for all $\phi \in C^{\infty}_0(R)$
\end{defn}
\begin{defn}{$S^{1,2} (R)$} 

\noindent We define $S^{1,2} (R)$ to be the
set of functions $u: R \rightarrow \mathbb{R}$ such that $u$, $u_t$ are in $L^2(R)$ and all of the horizontal derivatives of $u$ are in $L^2(R)$.
\end{defn}

The $S^{1,2}(R)$ norm is then given by
$$|| u ||_{{S^{1,2}}}=|| u
||_{L^2} + ||u_t||  +\sum_{i=1}^m|| X_{i}u ||_{L^2}.$$
We note that $S^{1,2}(R)$ is a Banach space under this norm.

In this paper we consider weak solutions. 
\begin{defn}{Weak solution to (\ref{heat})}

A weak solution to $$\parab$$ is defined as a function for which the first derivatives $u_t, X_1 u, ... , X_m u$ are square integrable in $R$ and satisfy 
$$\iint_{R} \phi u_t + \sum_{i,j=1}^m \big(a_{ij}(x,t) X_j u\big) \ X_i \phi \ dxdt \ = \ 0 $$ for every $\phi (x,t) \in C^{\infty}$  which has compact support in the $x$ variable for every fixed $t$.  
\end{defn}
\noindent We use an analogous definition for weak solutions to (\ref{nlheat})
\begin{defn}{Weak solution to (\ref{nlheat})}

A weak solution to $$\nlparab$$ is defined as a function for which the first derivatives $u_t, X_1 u, ... , X_m u$ are square integrable in $R$ and satisfy 
$$\iint_{R} \phi u_t + \sum_{j=1}^m \big(A_j(x, \hg u)\big) \ X_j \phi \ dxdt \ = \ 0 $$ for every $\phi (x,t) \in C^{\infty}$ which has compact support in the $x$ variable for every fixed $t$.  
\end{defn}
%put intro in here end.
\vfill

\pagebreak
%%%%%%%%%%%%%%%%%%%%%%%%%%%%%%%%%%%%%%%%%%%%%%%%%%%%%%%%%%%
% You cannot use the \chapter{} command because it puts a header at the top of every page which the grad school doesn't allow
\begin{center}
CHAPTER 1\\ \vspace{7pt}
\textbf{Parabolic Harnack Inequality}\\
\end{center} \vspace{7pt}
\setcounter{chapter}{1} %This makes sure the theorems and equations are labeled as 1.* instead of 0.*
\setcounter{equation}{0} %This makes sure the equations/eqnarrays etc are labeled starting with 1.1 instead of continuing the numbering
%------------------------- of the previous chapter.
\setcounter{thm}{0} %This makes sure the theorems/definitions etc are labeled starting with 1.1 instead of continuing the numbering
%-------------------------of the previous chapter.
%%%%%%%%%%%%%%%%%%%%%%%%%%%%%%%%%%%%%%%%%%%%%%%%%%%%%%%%%%%%%%%%%%
%-------------------------------------------------------------
% You cannot use the \chapter{} or \section{} commands because they put a header that the grad school doesn't allow
\begin{center} 
\textbf{1.1 Flow of the proof}
\end{center}
%-------------------------------------------------------------
 The flow of the proof will be as follows. We will first prove a Cacciopoli inequality for weak solutions $v$ to (\ref{heat}). We will show that these hold for various powers of the function. This will give us two inequalities for bounds on both the square of the horizontal gradient of $v$ and on the maximum of $v^2$ over all $t$ from above in terms of $v^2.$ We will the use a Sobolev embedding theorem by Capogna, Danielli and Garofalo \cite{CDG2} to allow for the estimate of higher powers of $v$ in terms of $v^2$. We will combine this with the bounds from the Cacciopoli inequality to form an inequality that we may iterate to achieve both the minimum and maximum of our function over specific domains. The iteration will work for integrals of $v^p$ and $v^{-p}$ for $p$ close to 0. We will bridge the gap between these two values with the results of Aimar \cite{Ai} and his John-Nirenberg style theorem for spaces of homogeneous type. In order to show that his results work, we need to show that $log(v)$ satisfies a certain bounded mean oscillation (BMO) condition. With this result, we will finish the proof of the parabolic Harnack inequality. Then, patterning our proof after Moser \cite{Mo1} we will show that this Harnack inequality implies H\"older continuity of solutions in a unit rectangle with respect to the CC metric. 
%-------------------------------------------------------------
% You cannot use the \chapter{} or \section{} commands because they put a header that the grad school doesn't allow
\begin{center} 
\textbf{1.2 Cacciopoli Inequality}
\end{center}
%-------------------------------------------------------------
 
In order to begin this proof we will construct a Cacciopoli inequality relating the value of $\int_{R'} v^2$ over $R' \subset R \subset \mathbb{R}^{n+1}$ to $\int_{R} |\hg v|^2$ and $\int_{R} v^2$.
 
\begin{thm}\label{localbounds}
Let $v \in S^{1,2}(R)$ be a positive, weak subsolution to 
\begin{equation}
 v_t = -\sum_{i,j=1}^{m} X^*_i \big(a_{ij}(x,t) X_j v\big)
\end{equation}
 in $R = \{B(x, r ) \times (0, \tau) \}$. Also define $R' = \{B(x, r ' ) \times (0, \tau ') \}$ for $0< r' <r$, $0 <\tau '< \tau$. Let $\{ X_i \}_{i=1}^m$, $m \leq n$, be smooth H\"ormander vector fields satisfying H\"ormander's condition for hypoellipticity. Let the matrix of coefficients $a_{ij}(x,t)$ be $L^{\infty}$ bounded and uniformly elliptic. That is 
$$\lambda^{-1} |\xi|^2 \leq a_{ij}(x,t) \xi_i \xi_j \leq \lambda |\xi|^2  $$
for every $(x,t) \in R$ and for every $\xi \in \mathbb{R}^n$. 
Then there exist constants $C_1$, $C_2$ depending on $\lambda$, $n$, and $R$ such that 
\begin{equation}
\label{caccgradv} \iint_{R'} | \hg v |^2 dxdt \ \leq C_1 \iint_{R} v^2 dxdt
\end{equation}

\begin{equation}
\label{caccmaxt} \underset{t \in [0 ,   \sigma ]}{\max} \int_{B(x, r')} v^2 dxdt \ \leq C_2 \iint_{R} v^2 dxdt
\end{equation}

Where the max is the essential supremum over this interval.
\end{thm}

%\begin{proof}
Consider the function $\phi (x,t) = v(x,t) \psi ^2(x,t)$. We will make the support of $\psi(x,t)$ compact in the $x$ variable for any fixed $t$ and differentiable. Notice this also makes $\phi \geq 0$. In order to accomplish the construction of $\phi$ we will recall the existence of test functions in Citti, Garofalo, and Lanconelli \cite{CiL}
\begin{lem}{ Existence of Carnot-Carath\'eodory test functions.}\label{CC_test_functions} 

There exists $r _0 > 0$ such that given a metric ball $B(x, r) \subset \subset R$, with $r \leq r _0$ and $0 < r ' < r$, there exists a function $\psi \in C^{\infty}_0 (B(x, r))$ such that 
$0 \leq \psi \leq 1$, $\psi \equiv 1$ in $B(x, r ')$ and $|\hg \psi| \leq \frac{C_{HG}}{r - r '}$. Here, $C_{HG} > 0$ is a constant independent of $r$ and $r '$.   
\end{lem}
\noindent Define $\psi (x,t) := \psi _1 (t) \psi _2 (x)$ in the following way:

%$$
\begin{equation}
\label{definition_of_psi}
\psi _1 (t) = \left\{
 \begin{array}{rl}
  1 & \text{if }  \tau' < t < \tau \\
   0 & \text{if } t \leq 0 \\
 \end{array} \right.
\end{equation}
%$$

$$
%\begin{equation*}
\psi _2 (x) = \left\{
 \begin{array}{rl}
  1 & \text{if } d_{cc}(x, 0) \leq r '\\
   0 & \text{if } d_{cc}(x, 0) > r \\
 \end{array} \right.
%\end{equation*}
$$
And $\psi _1 (t)$ is linearly interpolated from $0$ to $\tau'$. Since $v$ is a subsolution to (\ref{heat}) in $R$ we have 
\begin{equation}\iint_{R} \phi v_t \, dxdt \,+ \sum_{i,j=1}^m \iint_{R} X_i \phi \ \big(a_{ij}(x,t) X_j v\big) \ dxdt \leq 0
\end{equation}
\begin{equation}
\iint_{R} \psi ^2 v v_t  \, dxdt \,+ \sum_{i,j=1}^m \iint_{R} X_i (v \psi ^2) \ a_{ij}(x,t) X_j v \ dxdt \leq 0.
\end{equation}

\noindent We calculate $X_i(v \psi ^2) = X_i v \psi ^2 \ + \ 2 \psi X_i \psi v$ and substituting

\begin{equation}
\iint_{R} \psi ^2 v \, v_t  \, dxdt \, + \ \sum_{i,j=1}^m \iint_{R} X_i v \psi ^2 \,  a_{ij}(x,t) X_j v \ dxdt
\end{equation}
$$
 + \ \sum_{i,j=1}^m\iint_R 2 \psi X_i \psi v \ a_{ij}(x,t) X_j v \, dxdt \,\leq 0.
$$
\noindent Notice that $\frac{1}{2} \psi ^2 (v^2)_t = \psi ^2 v v_t$ which gives
\begin{equation}\label{cacc_ref_in_higher_powers}
\iint_{R} \frac{1}{2} \psi ^2 (v^2)_t \, dxdt+ \, \sum_{i,j=1}^m \iint X_i v \psi ^2 \, a_{ij}(x,t) X_j v \,dxdt
\end{equation}
$$
 \leq \, \sum_{i,j=1}^m -2 \iint_{R} v \psi X_i \psi \ a_{ij}(x,t) X_j v \, dxdt.
$$
\noindent Now consider the RHS of the equation 
\begin{equation}
\sum_{i,j=1}^m v \psi \, X_i \psi \, a_{ij}(x,t) \, X_j v 
\end{equation}
and the Schwarz inequality:
\begin{equation}
| \sum_{i,j=1} ^m v \psi \, X_i \psi \, a_{ij}(x,t) \, X_j v | \ \leq \ \big(\psi ^2 \, \langle \hg v, a \hg v \rangle \langle \hg \psi , a \hg \psi \rangle  \, v^2\big)^{ \frac{1}{2}}.
\end{equation}

\noindent And using Young's inequality we get:
\begin{equation}
\leq \frac{1}{4} \langle \hg v, a \hg v \rangle  \psi ^2 \ + \ v^2  \langle \hg \psi, a \hg \psi \rangle .
\end{equation}

\noindent Now substitute back into the RHS of (\ref{cacc_ref_in_higher_powers}) to get:
\begin{equation}
\iint_{R} \frac{1}{2} \psi ^2 (v^2)_t \ dxdt\ + \ \iint_{R}  \langle \hg v, a \hg v\rangle  \psi ^2 dxdt
\end{equation}

\begin{equation}
\leq 2 \iint_{R} \frac{\psi ^2}{4}  \langle \hg v, a \hg v \rangle  dxdt\ + \ v^2  \langle \hg \psi , a \hg \psi \rangle .
\end{equation}

\noindent Now gathering the $\psi ^2  \langle \hg v, a \hg v \rangle $ terms onto the LHS of (\ref{cacc_ref_in_higher_powers}) we get:

\begin{equation}
\iint_{R} \frac{1}{2} \psi ^2 (v^2)_t \ dxdt + \ \iint_{R} \frac{1}{2}  \langle \hg v, a \hg v \rangle  \psi ^2 dxdt 
\end{equation}
$$
\leq \ \iint_{R} v^2  \langle \hg \psi , a \hg \psi \rangle dxdt .
$$
\noindent Use the ellipticity and boundedness of $a_{ij}(x,t)$ to get
\begin{equation}
\iint_{R} \frac{1}{2} \psi ^2 (v^2)_t \, dxdt \ + \ \iint_{R} \frac{1}{2 \lambda} | \hg v|^2 \psi ^2 dxdt \ \leq \ 2 \lambda \iint_{R} v^2 | \hg \psi |^2 dxdt
\end{equation}

\noindent Notice that $\frac{1}{2}(\psi ^2 v^2)_t = \psi \psi _t v^2 + v v_t \psi ^2$. And use this fact by adding $v^2 \psi \psi_t $ to both sides: 

\begin{equation}
\iint_{R} \frac{1}{2} \psi ^2 (v^2)_t \, dxdt \ + \ v^2 \psi \psi _t  \, dxdt \ + \ \iint_{R} \frac{1}{2 \lambda} | \hg v|^2 \psi ^2  dxdt \ 
\end{equation} 
$$\leq \ 2 \lambda \iint_{R} v^2 (| \hg \psi |^2 +|\psi \psi _t|) dxdt $$
This gives
\begin{equation}
\label{caccopoli} \iint_{R} \frac{1}{2}  ( \psi ^2 v^2)_t \, dxdt  \ + \ \iint_{R} \frac{1}{2 \lambda} | \hg v|^2 \psi ^2  dxdt \ 
\end{equation}
$$
\leq \ 2 \lambda \iint_{R} v^2 \big(| \hg \psi |^2 +|\psi \psi _t|\big) dxdt .
$$
\noindent Now integrate this over $R_{\sigma} = \{B(x, r) \times 0 \leq t \leq \sigma < \tau \ \}$ and drop the second integrand on the LHS. We note that 
$$\iint_{R_{\sigma}} \frac{1}{2 \lambda} | \hg v|^2 \psi ^2  dxdt \geq 0.$$
We choose $\sigma$ so that 
\begin{equation}
\int_{B(x, r')} v^2 |_{t =  \sigma} \ dx \ \geq \ \frac{1}{2} \ \underset{t \in (0 ,  \sigma)}{\max} \int_{B(x, r ')} v^2 dx.
\end{equation}
 \noindent This gives us:
\begin{equation}
\label{maxtbounds}\iint _{R_{\sigma}} \frac{1}{2} (v^2 \psi ^2)_t \, dx dt \leq 2 \lambda \iint_{R_{\sigma}} v^2 \big(|\hg \psi|^2 \ + \ |\psi \, \psi _t|\big) \, dxdt.
\end{equation} 

\noindent Now consider the 
\begin{equation}
\underset{t \in (0 ,  \sigma)}{\max} \int_{B(x, r ')} v^2 \, dx \ \leq \ 2 \int_{B(x, r ')} v^2 \  |_{t = \sigma} \ dx
\end{equation} 
which we constructed above. $\psi = 1$ on $R'$ and we fix $t$, so 

$$\underset{t \in (0 ,  \sigma)}{\max} \int_{B(x, r ')} v^2 \, dx \, \leq \, 2 \int_{B(x, r ')} v^2 \,  |_{t = \sigma} \ dx $$

$$\leq \, 2 \int_{B(x, r ')} v^2 \psi ^2 |_{t =  \sigma} \, dx$$

$$\leq \ 2 \int_{B(x, r)} v^2 \psi ^2 |_{t =  \sigma} \, dx$$
And since $\psi$ vanishes at $t = 0$
$$= \ 2 \int_{B(x, r)} v^2 \psi ^2 |_{t =  \sigma} - v^2 \psi ^2 |_{t = 0}\, dx \ = \ 2 \iint_{R_{\sigma}} (v^2 \psi ^2)_t dxdt $$

\noindent Now recall (\ref{maxtbounds}) integrated over $R_{\sigma}$. We will combine that with the $\frac{1}{2 \lambda}$ factor as well as the fact that the above is four times the (\ref{maxtbounds}) inequality. Putting these together gives:
\begin{equation}
\underset{t \in (0 ,   \sigma)}{\max} \int_{B(x, r ')} v^2 dx \ \leq \ 8 \lambda \iint_{R_{\sigma}} v^2 \big(|\hg \psi|^2 + |\psi \psi _t|\big) \, dx dt.
\end{equation}

\noindent Then use the non-negativity of the RHS and take the integral over the larger rectangle $R$ for
\begin{equation}
\underset{t \in (0 ,  \sigma)}{\max} \int_{B(x, r ')} v^2 dx \ \leq 8 \lambda \iint_{R} v^2 \big(|\hg \psi|^2 + |\psi \psi _t|\big) \, dx dt.
\end{equation}

\noindent In creating the test function $\phi$ according to Lemma \ref{CC_test_functions} above, we use the result on the bounds of the horizontal gradient of test functions in Citti, Garofalo, and Lanconelli \cite{CiL} in order to attain a constant independent of  $ \psi $. Specifically, from Lemma \ref{CC_test_functions} we use $|\hg \psi| \leq \frac{C_{HG}}{r - r '}$. This allows us to say
\begin{equation}
|\hg \psi |^2 \ + \ |\psi \psi _t | \ \leq \ \big(\frac{C_{HG}^2}{(r - r ')^2} + \frac{1}{\tau - \tau '}\big)
\end{equation}

\noindent Which gives us
\begin{equation}
\underset{t \in (0 ,   \sigma)}{\max} \int_{R '} v^2 \, dx \ \leq \ 8 \lambda \big(\frac{C_{HG}^2}{(r - r ')^2} + \frac{1}{\tau - \tau '}\big) \iint_{R} v^2 \, dxdt.
\end{equation}
This proves (\ref{caccmaxt}). Now look at the second integrand on LHS of (\ref{caccopoli}). 
\begin{equation}
\iint_{R'} \frac{1}{2 \lambda} |\hg v|^2 \, \psi ^2 dxdt \ \leq \ \iint_{R'} \frac{1}{2 \lambda} |\hg v|^2 \, dxdt 
\end{equation}
\begin{equation}
\leq 2 \lambda \iint_{R'} v^2 \, (|\hg \psi |^2 + |\psi \psi _t|) \, dxdt
\end{equation}

\noindent Which gives us
\begin{equation}
\iint_{R'} |\hg v|^2 \, \psi ^2 dxdt \ \leq 4 \lambda ^2 \iint_{R'} v^2 \, (|\hg \psi |^2 + |\psi \psi _t|) \, dxdt
\end{equation}
\begin{equation}
\leq 8 \lambda ^2 (\frac{C_{HG}^2}{(r - r ')^2} + \frac{1}{\tau - \tau '}) \iint_{R} v^2   \, dxdt
\end{equation}

\noindent And the integral on the RHS is taken over R since it is non-negative. This gives 
\begin{equation}
 \iint_{R'} |\hg v|^2 \,   dxdt \ \leq 8 \lambda ^2 (\frac{C_{HG}^2}{(r - r ')^2} + \frac{1}{\tau - \tau '}) \iint_{R} v^2   \, dxdt
\end{equation} 
with $\psi = 1$ on $R'$. 
%\end{proof}
%-------------------------------------------------------------
% You cannot use the \chapter{} or \section{} commands because they put a header that the grad school doesn't allow
\begin{center} 
\textbf{1.3 Cacciopoli for other values of p}
\end{center}
%-------------------------------------------------------------

Above we proved a Cacciopoli type inequality, (Theorem \ref{localbounds}) for positive, weak solutions to 
$$\parab$$  
\noindent in $R \subset \mathbb{R}^{n+1}$ for smooth H\"ormander vector fields. Now we need to show that this result holds for powers of $u$. We will consider a function $v=u^p$ for $p>{1  \over 2}$. We also need similar inequalities for $p<0$ and $0<p<1$. This, along with the Sobolev style embedding inequality will provide the basis of the Moser-type iteration. 
\begin{prop}{(Cacciopoli for powers of u)}

Under the assumptions of Theorem \ref{localbounds} above, the results are true for powers of $v := u^{p}$, for $p > \frac{1}{2}$,
\begin{equation}
  \iint_{R'} | \hg (u^{p}) |^2 dxdt \ \leq C_1 \iint_{R} (u^{p})^2 dxdt
\end{equation}

\begin{equation}
 \underset{t \in [0 ,   \sigma ]}{\max} \int_{B(x, r')} (u^{p})^2 dxdt \ \leq C_2 \iint_{R} (u^{p})^2 dxdt
\end{equation}
If $u$ is a weak supersolution to (\ref{heat}) and $p<0$ then we have a similar result with the inequality reversed:
\begin{equation}
  \iint_{R'} | \hg (u^{p}) |^2 dxdt \ \geq C_1 \iint_{R} (u^{p})^2 dxdt
\end{equation}

\begin{equation}
 \underset{t \in [0 ,   \sigma ]}{\max} \int_{B(x, r')} (u^{p})^2 dxdt \ \geq C_2 \iint_{R} (u^{p})^2 dxdt
\end{equation}
\end{prop}
\noindent \emph{Proof:}

Let $v:= u^p$ be a positive subsolution to $$\parab$$ in the rectangle $R$ as above and $\phi = \psi^2 u^p$. Then we have 
\begin{equation}
\iint_{R} \phi v_t \ dxdt+ \sum_{i,j=1}^m \iint_{R} \, X_i \phi (a_{ij}(x,t) X_j v) \, dxdt \, \leq 0
\end{equation} 
Calculate 
$$v_t = (u^p)_t = pu^{p-1}u_t$$
and
$$X_i(\phi) = X_i(\psi^2 u^p) = p u^{p-1} \psi ^2 X_i u + 2X_i \psi u^p \psi. $$
We substitute the definitions of $\phi$ and $v$:
$$
\iint_{R} \psi ^2 u^p p u^{p-1} u_t \ dxdt 
$$

$$
+  \sum_{i,j=1}^m \iint_{R} (p u^{p-1} \psi ^2 X_i u + 2X_i \psi u^p \psi) a_{ij}(x,t) (p u^{p-1} X_j u) \, dxdt \, \leq 0
$$ 

$$
\iint_{R} \psi ^2 u^{2p-1} u_t +  p^2 u^{2p-2} \psi ^2 \langle \hg u, a \hg u \rangle  dxdt 
$$

$$
+ \iint_R 2p \psi u^{2p-1} \langle \hg \psi, a \hg u \rangle  \, dxdt \, \leq \, 0.
$$
We recognize that $2p u^{2p-1} u_t = \frac{1}{2} \psi ^2 (u^{2p})_t \, $:

$$
\iint_{R} \frac{1}{2} \psi ^2 (u^{2p})_t \, dxdt +  p^2 u^{2p-2} \psi ^2 \langle \hg u, a \hg u  \rangle    
$$

$$
+ \iint_{R}  2p \psi u^{2p-1}  \langle \hg \psi, a \hg u \rangle  \, dxdt \, \leq \, 0
$$ 

\noindent Examining the RHS we have by the Schwarz and Young's inequalities:
\begin{equation}
-2p \psi u^{2p-1}  \langle \hg \psi, a \hg u \rangle  \, \leq \, |2p \psi u^{2p-1}  \langle \hg \psi, a \hg u \rangle |
\end{equation}
\begin{equation}
\leq 2p( \langle \hg \psi, a \hg \psi \rangle u^{2p} \ \psi ^2  \langle \hg u, a \hg u  \rangle  u^{2p-2})^{\frac{1}{2}}
\end{equation} 
\begin{equation}
\leq 2p ( \langle \hg \psi, a \hg \psi \rangle  u^{2p} + \frac{1}{4}  \langle \hg u, a \hg u  \rangle  \psi ^2 u^{2p-2}).
\end{equation}

\noindent Substituting back in and combining like terms we have
$$
\iint_{R} \psi ^2 (u^{2p})_t \, dxdt  + \iint_{R} (p^2 -\frac{p}{2}) u^{2p-2} \psi ^2  \langle \hg u, a \hg u \rangle  \, dxdt \, $$

$$\leq \, \iint_{R} 2p \langle \hg \psi, a \hg \psi \rangle  u^{2p} \, dxdt$$
 
\noindent At this point, observe $v = u^p, \ \hg v = p u^{p-1} \hg u$ and use both boundedness and ellipticity of $a_{ij}(x,t)$
\begin{equation}
\iint \psi ^2 (v^{2})_t  \, dxdt + \iint \frac{1}{\lambda} \frac{2p-1}{p} \psi ^2 |\hg v|^2 \, dxdt \, \leq \, \iint 2p \lambda |\hg \psi|^2 \, v^2 \, dxdt
\end{equation}
 
\noindent At which point, we have the same inequality as (\ref{caccopoli}) and the proof follows as above. 
For the case $p < p_0 < 0$ we let $v= u^{p/2}$ be a supersolution to (\ref{heat}) in $R$. We then have for $\phi = v \psi ^2$ the following:
\begin{equation}
v_t = (u^{p/2}) = \frac{p}{2} u^{p/2 -1} u_t,
\end{equation} 
\begin{equation}
X_i v = X_i (u^{p/2}) = \frac{p}{2} u^{p/2 -1} X_i u
\end{equation}
and 
\begin{equation}
X_i \phi = X_i (v \psi ^2) = X_i (u^{p/2} \psi ^2) = \frac{p}{2} u^{p/2 -1} X_i u \psi ^2 + 2 \psi X_i \psi u^{p/2}
\end{equation}
And we substitute these into 
\begin{equation}
\iint_{R} \phi v_t  \, dxdt + \sum_{i,j=1}^m \iint_{R} \, X_i \phi (a_{ij}(x,t) X_j v) \, dxdt \geq 0.
\end{equation} 
to get
$$
\iint_{R} u^{p/2} \psi^2 (\frac{p}{2} u^{p/2 -1} u_t) dxdt $$

$$+ \; \iint_R  [\frac{p}{2} u^{p/2 -1} X_i u \psi ^2 + 2 \psi X_i \psi u^{p/2}] \, a_{ij}(x,t) (\frac{p}{2} u^{p/2 -1} X_i u) \, dxdt \,\geq 0.
$$ 
Since 
\begin{equation}
\frac{1}{2}(v^2)_t =  v v_t = (u^{p/2})(\frac{p}{2} u^{p/2 -1} u_t)
\end{equation} 
we can arrive at
$$
\iint_{R} \frac{1}{2}(v^2)_t \psi ^2 + \psi ^{2} \langle \hg v, a \hg v \rangle  \, dxdt \,$$

$$
 \geq -2 \iint_{R} v \psi  \langle \hg \psi, a \hg v \rangle  \, dxdt
$$
At this point, we have inequality (\ref{cacc_ref_in_higher_powers}) and we may proceed in the proof exactly as above. The case for $0 < p < 1$ is modified in the following way: we must reverse $t$ and use $v:= u^p(x,-t)$. Then our test function must be $\phi = u^p \psi^2$. The proof proceeds in the same way except there is a factor of $\frac{1}{\lambda} \frac{2p-1}{p^2} $ in front of the $|\hg v|^2$ term. That is: 

\begin{equation}
\iint \psi ^2 (v^{2})_t  \, dxdt + \iint \frac{1}{\lambda} \frac{2p-1}{p^2} \psi ^2 |\hg v|^2 \, dxdt \, \leq \, \iint 2p \lambda |\hg \psi|^2 \, v^2 \, dxdt
\end{equation}

 %-------------------------------------------------------------
% You cannot use the \chapter{} or \section{} commands because they put a header that the grad school doesn't allow
\vspace{5pt}
\begin{center} 
\textbf{1.4 Sobolev Theorem and Iteration}
\end{center} Above we proved the following estimates, (Theorem \ref{localbounds}):
 
$$
\iint_{R'} | \hg v |^2 dxdt \ \leq C_1 \iint_{R} v^2 dxdt
$$
$$
\underset{t \in [0 , \sigma ]}{\max} \int_{B(x, r')} v^2 dx \ \leq C_2 \iint_{R} v^2 dxdt
$$
for positive, weak solutions to (\ref{heat}) in $R$. Then we showed that these held for powers of $u$, $v:= u^p$. Now we will combine these with the Sobolev embedding theorem to begin the Moser iteration. We introduce the
\begin{thm}{(Sobolev Embedding Theorem)} Let $U \subset \mathbb{R}^n$ be a bounded open set, then there exists a $C_S>0$, $r _0 > 0$ such that for any $x \in U$, $B_{r} = B_{cc}(x, r)$, with $r \leq r _0$, here $1 < k \leq 2$
\begin{equation}
(\frac{1}{|B_{r}|} \int_{B_{r}} |u|^{2k} dx)^{\frac{1}{2k}} \ \leq \ C_S(\frac{1}{|B_{r}|} \int_{B_{r}} |\hg u|^2 dx)^{\frac{1}{2}}
\end{equation}
\noindent for any $u \in S_0^{1,2} (B_r)$
\end{thm}
\noindent Squaring both sides, we will use this formula in the form
\begin{equation}
(\frac{1}{|B_{r}|} \int_{B_{r}} |u|^{2k}dx )^{\frac{1}{k}} \ \leq \ C_S(\frac{1}{|B_{r}|} \int_{B_{r}} |\hg u|^2 dx)
\end{equation}
In order to use the Sobolev embedding theorem in the setting where our function $v$ is not compactly supported but merely defined on the ball $B(x,r)$ we will recall the cut off function $\psi(x,t)$ defined above:

\begin{equation}
\psi _1 (t) = \left\{
 \begin{array}{rl}
  1 & \text{if } \tau' < t < \tau \\
   0 & \text{if } t \leq 0 \\
 \end{array} \right.
\end{equation}
%$$

$$
%\begin{equation*}
\psi _2 (x) = \left\{
 \begin{array}{rl}
  1 & \text{if } d_{cc}(x, 0) \leq r '\\
   0 & \text{if } d_{cc}(x, 0) > r \\
 \end{array} \right.
%\end{equation*}
$$
and as before we have a bound on the horizontal gradient. So we set $w = v \cdot \psi \in S^{1,2}_0(B(x,r))$ and notice that on $B(x, r')$ we have:

$$w(x,t) = v(x,t) \cdot \psi (x,t) =v(x,t) \cdot \psi_1(t) \cdot \psi_2(x) = v(x,t) \cdot 1 = v(x,t).
$$
Also we have 

$$\hg w = \hg(v \cdot \psi) = \hg v \psi + \hg \psi v = \hg v \psi + 0 \cdot v = \hg v \psi = \hg v.$$
Then on $B(x, r)$ we have:

$$w(x,t) = v(x,t) \cdot \psi (x,t).$$
Also we have 
$$\hg w = \hg(v \cdot \psi) = \hg v \psi + \hg \psi v $$
Now we calculate for the function $w$ above:

$$(\frac{1}{|B_{r'}|} \int_{B_{r'}} |v|^{2k} dx)^{\frac{1}{k}} \ \leq (\frac{1}{|B_{r'}|} \int_{B_{r'}} |v \psi|^{2k} dx)^{\frac{1}{k}} \underset{def}{=}  (\frac{1}{|B_{r'}|} \int_{B_{r'}} |w|^{2k} dx)^{\frac{1}{k}}$$
$$
= (\frac{|B_r|}{|B_r||B_{r'}|} \int_{B_{r'}} |w|^{2k} dx)^{\frac{1}{k}}$$
Then on the larger domain we have 
\begin{equation}\label{for_balls}\leq  (\frac{|B_r|}{|B_{r'}|}\frac{1}{|B_{r}|} \int_{B_{r}} |w|^{2k} dx)^{\frac{1}{k}} =  (\frac{|B_r|}{|B_{r'}|})^{\frac{1}{k}}(\frac{1}{|B_{r}|} \int_{B_{r}} |w|^{2k} dx)^{\frac{1}{k}}
\end{equation}
which is the domain over which we may apply the Sobolev embedding theorem since $w$ has compact support. Now we will deal with the size of the CC balls with the following proposition from Nagel, Stein and Wainger, \cite{NaW} and use 
\begin{prop}\label{citti}
Given a bounded set $U \subset \mathbb{R}^n$ there exist $Q \geq n$, $R_0 > 0$, $C_B > 0$, such that for every $x \in U$, $R \leq R_0$ and $0 < t < 1$
\begin{equation}
\label{bounds_on_balls}|B(x, tR)| \geq C_B t^Q |B(x, R)|
\end{equation}
\end{prop}
 
\noindent In our use of the proposition, we have $t = \frac{r '}{r}$ and use
\begin{equation}
\frac{|B(x, r)|}{|B(x, t r)|} = \frac{|B(x, r)|}{|B(x, \frac{r '}{r } r)|} \leq C_B(\frac{r}{r '})^Q 
\end{equation}
Now estimate the size of the balls with $C_B(\frac{r}{r '})^Q$. From (\ref{for_balls})
\begin{equation} 
\leq \ C_S [C_B(\frac{r}{r '})^{Q/k}](\frac{1}{|B_{r}|} \int_{B_{r}} |\hg w|^2 dx)  
\end{equation}
Now substitute the definitions of $w = v \psi$, and $\hg w = \hg v \psi + \hg \psi v$
\begin{equation}
\label{sob_pre_schwarz_youngs}= C_S[C_B(\frac{r}{r '})^{Q/k}](\frac{1}{|B_{r}|} \int_{B_{r}} |\hg v \psi + \hg \psi v|^2 dx ).  \ 
\end{equation}
We notice that 
$$|\hg v \psi + \hg \psi v|^2 = |\hg v|^2 \psi^2 + v^2 |\hg \psi|^2 + 2 \psi \langle (\hg v)v \hg \psi \rangle $$
Now apply the Schwarz inequality
$$\leq |\hg v|^2 \psi^2 + v^2 |\hg \psi|^2 + 2[(|\hg v|^2 \psi^2 v^2 |\hg \psi|^2)^{1/2} ]$$  
And Young's inequality
$$\leq |\hg v|^2 \psi^2 + v^2 |\hg \psi|^2 + [(|\hg v|^2 \psi^2) +  (|v|^2)|\hg \psi|^2]$$
$$= 2|\hg v|^2 \psi^2 + 2v^2 |\hg \psi|^2$$
Since we are on $B(x,r)$, $\psi \leq 1$ and $|\hg \psi| \leq C_{HG}/(r-r')^2$
 $$= 2|\hg v|^2   + 2C_{HG}v^2 /(r-r')^2.$$
We replace this back into (\ref{sob_pre_schwarz_youngs}):
\begin{equation}
\leq	C_S[C_B(\frac{r}{r '})^{Q/k}](\frac{1}{|B_{r}|} \int_{B_{r}} [2|\hg v|^2   + 2C_{HG}v^2 /(r-r')^2]  dx)
\end{equation}
$$
=  2 \ C_S[C_B(\frac{r}{r '})^{Q/k}](\frac{1}{|B_{r}|} \int_{B_{r}} |\hg v|^2 + (C_{HG}v^2/(r-r')^2) dx\, )
$$
This allows us to say
\begin{equation}
\label{sobolev_wo_compact}
(\frac{1}{|B_{r'}|} \int_{B_{r'}} |v|^{2k} dx)^{\frac{1}{k}} \ \leq 2 \ C_S[C_B(\frac{r}{r '})^{Q/k}](\frac{1}{|B_{r}|} \int_{B_{r}} |\hg v|^2 + (C_{HG}v^2/(r-r')^2) \,dx ).
\end{equation}
Now we will use this estimate from the Sobolev embedding theorem to prove 

\begin{thm} 
Let $u \in S^{1,2}(R)$ be a positive, weak subsolution to (\ref{heat}) in $R = \{B(x, r ) \times (0, \tau) \}$. Also define $$R' := \{B(x, r ' ) \times (0, \tau ') \}$$ 
for $0<r'< r $, $0<\tau'< \tau$. $\{ X_i \}^{m}_{i=1}$, $m \leq n$, are smooth H\"ormander vector fields and the matrix of coefficients $a_{ij}(x,t)$ are $L^{\infty}$ bounded, measurable and uniformly elliptic satisfying
$$\lambda^{-1} |\xi|^2 \leq a_{ij}(x,t) \xi_i \xi_j \leq \lambda |\xi|^2  $$
for every $(x,t) \in R$ and for every $\xi \in \mathbb{R}^n$. Then there exists a constant $\gamma$ depending on $\lambda$, $n$, and $R$ such that
\begin{equation}
\label{maxuleqpnorm}\underset{R'}{\max} \ u \ \leq \gamma (\frac{1}{|R|} \iint_{R} u^p \, dxdt )^{1/p} 
\end{equation}
for $p>1$. 
\end{thm}
\emph{Proof:}
In order to prove the theorem, we need a lemma modeled after Moser \cite{Mo1}. Now for $v \in S^{1,2}(R)$, define 
\begin{equation}
H_{r,  \tau}(v):=  \normbrt \iint_R v^2 \, dxdt
\end{equation}
\begin{equation}
M_{r, \tau} (v) :=  \normbr \, \underset{t \in [0, \tau ]}{\max} \int_{B_{r}} v^2 dx
\end{equation}
\begin{equation}
D_{r, \tau}(v) := \normbrt \, \iint_R |\hg v|^2 dxdt
\end{equation}

\begin{lem}
Let $v \in S^{1,2}(R)$. Let $k=1 + \frac{2}{n}$ for $n \in \mathbb{N}$, $\alpha = \frac{n}{n-2}$, $\beta = \frac{n}{2}$ Then, we have for some constants $C_1$, $C_2 >0$: 
\begin{equation}
\label{preliminaryiteration}H_{r',  \tau}(v^k) \leq C_1 \, M_{r,  \tau}(v)^{\frac{2}{n}} \bigl[ C_2 H_{r',  \tau}(v) + D_{r',  \tau}(v) \bigr]
\end{equation}

\end{lem}
\noindent \emph{Proof:}

Notice that 
$$\frac{1}{\alpha} + \frac{1}{\beta} = \frac{n-2}{n} + \frac{2}{n} = 1.$$ In particular, $\alpha$ and $\beta$ are H\"older conjugates. First consider

$$\int_{B_{r'}} v^{2k} \, dx = \int_{B_{r'}} v^{2(1+ \frac{2}{n})} \, dx = \int_{B_{r'}} v^2 v^{\frac{4}{n}} \, dx$$

\noindent And by H\"older's inequality we have

$$\leq (\int_{B_{r'}} v^{2 \alpha} \, dx)^{\frac{1}{\alpha}}(\int_{B_{r'}} v^{\frac{4 \beta}{n}} \, dx)^{\frac{1}{\beta}} = (\int_{B_{r'}} v^{2 \alpha} \, dx)^{\frac{1}{\alpha}}(\int_{B_{r'}} v^2 \, dx)^{\frac{2}{n}}$$

$$ \leq (\int_{B_{r'}} v^{2 \alpha} \, dx)^{\frac{1}{\alpha}} M(v)^{\frac{2}{n}}|B_r|^{\frac{2}{n}}$$

\noindent Since the max over all $t$ is greater than or equal to the value for any $t$. We will use Sobolev's embedding theorem on the RHS. We notice now the step mentioned earlier. The integral
\begin{equation}
\label{integral_on_which_to_use_sobolev}(\int_{B_{r'}} v^{2 \alpha} \, dx)^{\frac{1}{\alpha}}
\end{equation}
is taken over the ball, $B_{r'}$ and we will use the results from above. Use the Sobolev embedding theorem on the (\ref{integral_on_which_to_use_sobolev}) term in the form of the inequality (\ref{sobolev_wo_compact}) to achieve:
$$
\int_{B_{r'}} v^{2k} \, dx \, \leq  2 \ C_S[C_B(\frac{r}{r '})^{Q/k}]( \int_{B_{r}} |\hg v|^2 + (C_{HG}v^2/(r-r')^2) \,dx ) \, M(v)^{\frac{2}{n}} |B_r|^{\frac{2}{n}}
$$
 
\noindent Integrate in t over  $[0'$,$   \sigma]$ to get:
$$
|B_{r'} \times (0,\tau')|H_{r', \tau'}(v^{k}) 
$$
$$ \leq 2 \ C_S[C_B(\frac{r}{r '})^{Q/k}]|B_{r} \times (0,\tau')|$$
$$\times \bigl[  D_{r, \tau'}(v) +  (1/r^2) H_{r, \tau'}(v) \ \bigl] M_{r, \tau}(v)^{\frac{2}{n}}|B_{r}|^{\frac{2}{n}}.
$$
Now we combine constants and divide the size of the $|B_{r'} \times (0,\tau')|$:
$$ H_{r', \tau'}(v^{k}) \leq 2 \ C_S [C_B(\frac{r}{r '})^{Q/k}] \frac{|B_{r} \times (0,\tau')|}{|B_{r'} \times (0,\tau')|} $$
$$\times \bigl[  D_{r, \tau'}(v) \, + \, (1/r^2) H_{r, \tau'}(v) \, \bigl] M_{r, \tau}(v)^{\frac{2}{n}}|B_{r}|^{\frac{2}{n}}.
$$
Before we finish with the creation of the iteration inequality we will deal with the size of the CC balls with the proposition from Nagel, Stein and Wainger, \cite{NaW} and use  

\begin{equation}
\frac{|B(x, r)|}{|B(x, t r)|} = \frac{|B(x, r)|}{|B(x, \frac{r '}{r } r)|} \leq C_B(\frac{r}{r '})^Q 
\end{equation}
Now estimate the size of the balls with $C_B(\frac{r}{r '})^Q$.
$$ H_{r', \tau'}(v^{k}) \leq 2 \ C_S [C^2_B(\frac{r}{r '})^{Q + Q/k} ] \bigl[  D_{r, \tau'}(v) \, + \, (1/r^2) H_{r, \tau'}(v) \, \bigl] M_{r, \tau}(v)^{\frac{2}{n}}|B_{r}|^{\frac{2}{n}}.
$$
\begin{equation}
\label{before_iteration}H_{r', \tau'}(v^{k}) \, \leq 2 \ C_S [C^2_B(\frac{r}{r '})^{Q + Q/k} ]|B_{r}|^{\frac{2}{n}} \bigl[  D_{r, \tau'}(v) \, + \, ({1 \over r^2}) H_{r, \tau'}(v) \, \bigl] M_{r, \tau}(v)^{\frac{2}{n}}.
\end{equation}
This proves the lemma. In terms of $H_{r, \tau}(v)$, $D_{r, \tau}(v)$ and $M_{r, \tau}(v)$ we calculate our local bounds (\ref{caccgradv}) and (\ref{caccmaxt}) as:
\begin{equation}
M_{r , \tau } (v) = \normbrp \underset{t \in (0, \tau )}{\max} \int_{B_{r }} v^2 dx 
\end{equation}
\begin{equation}
 \leq  8 \lambda (\frac{C_{HG}^2}{(r - r ')^2} + \frac{1}{\tau - \tau '}) \normbr \iint_{R} v^2 dx dt
\end{equation}
\begin{equation}
\leq  8 \lambda (\frac{C_{HG}^2}{(r - r ')^2} + \frac{1}{\tau - \tau '})  H_{r, \tau}(v)
\end{equation}

\noindent Then for some $C_3$ depending on $\lambda$, $R$ and $n$ we have 
\begin{equation}
M_{r , \tau } (v) \leq C_3 \, H_{r ,\tau}(v)
\end{equation}

\noindent And
$$D_{r , \tau '} (v) = {1 \over |B_{r} \times (0, \tau ')|}  \, \iint_{R} |\hg v|^2 \,dxdt  $$
where we need to replace the $\tau '$ with a $\tau$ 
$$D_{r , \tau '} (v) \leq \, {1 \over |B_{r} \times (0, \tau ')|} \, \frac{| (0, \tau)|}{| (0, \tau)|} \, 8 \lambda ^2 (\frac{C_{HG}^2}{(r - r ')^2} + \frac{1}{\tau - \tau '}) \iint_{R} v^2 dxdt$$
$$\leq \, {1 \over |B_{r} \times (0, \tau )|} \, \frac{| (0, \tau)|}{| (0, \tau')|} \, 8 \lambda ^2 (\frac{C_{HG}^2}{(r - r ')^2} + \frac{1}{\tau - \tau '}) \iint_{R} v^2 dxdt$$

$$= \frac{ \tau}{\tau '} \, 8 \lambda ^2 (\frac{C_{HG}^2}{(r - r ')^2} + \frac{1}{\tau - \tau '}) \, H_{r ,\tau} (v)$$

gives for some $C_4$ dependent on $\lambda$, $n$ and $R$
\begin{equation}
\label{localboundsonD}D_{r ', \tau '} (v) \leq \, C_4 \, H_{r ,\tau} (v)
\end{equation}

\noindent So substituting these into the results of the lemma (\ref{preliminaryiteration}) above, we get

$$H_{r ', \tau '} (v^{k}) \leq  C_1 C_3 H_{r,\tau }(v)^{\frac{2}{n}} [C_4 H_{r, \tau }(v)+ C_2H_{r, \tau }(v)]  $$

$$=  C_1 C_3[C_4 + C_2 ] H_{r,\tau }(v)^{\frac{2}{n}}  H_{r, \tau }(v)$$

$$= \gamma \, H_{r , \tau }(v)^{\frac{2}{n} + 1}$$
\begin{equation}
= \gamma \, H_{r , \tau }(v)^k
\end{equation}
where 
$$
\label{gamma}\gamma = 2 C_S [C^2_B(\frac{r}{r '})^{Q + Q/k} ] |B_r|^{2/n}$$
$$ \times \ [  8 \lambda ^2 \frac{ \tau}{\tau '} (\frac{C_{HG}^2}{(r - r ')^2} + \frac{1}{\tau - \tau '}) +1/r^2 ] \bigl(  8 \lambda   (\frac{C_{HG}^2}{(r - r ')^2} + \frac{1}{\tau - \tau '}) \,     \bigr)^{\frac{2}{n}}.
$$

%$H_{r ' \tau '} (w^{2k}) \leq  \Bigl( 8 \lambda ^2 (\frac{C_{HG}^2}{(r - r ')^2} + \frac{1}{\tau - \tau '}) +1 \Bigr) H_{r  \tau } (w) \Bigl( 8 \lambda (\frac{C_{HG}^2}{(r - r ')^2} + \frac{1}{\tau - \tau '}) H_{r  \tau }(w) \Bigr) ^{\frac{2}{n}}$

%$H_{r ' \tau '} (w^{2k}) \leq  \Bigl( 8 \lambda ^2 (\frac{C_{HG}^2}{(r - r ')^2} + \frac{1}{\tau - \tau '}) +1 \Bigr)[8 \lambda (\frac{C_{HG}^2}{(r - r ')^2} + \frac{1}{\tau - \tau '})]^{\frac{2}{n}} H_{r  \tau } (w)   H_{r  \tau }(w) ^{\frac{2}{n}}$

%$H_{r ' \tau '} (w^{2k}) \leq  C(\lambda, R, n) \,  H_{r  \tau }(w) ^{\frac{2}{n} + 1} = C(\lambda, R, n) \,  H_{r  \tau }(w) ^k $

\noindent This shows that for a positive, weak subsolution $v$ in $R$ that is square integrable, then $v^k$ is square integrable in $R'$. We can then repeat this exponentiation by $k$ to achieve the max of $v$ in a smaller domain. For example, take
\begin{equation}
\label{vnu} v_{\nu} := v^{p_{\nu} /2}
\end{equation}
a positive, weak subsolution in $R$ where $p_{\nu} = p_0 k^{\nu}$, $p_{\nu} \geq p_0 \geq p' >1$ and $k = 1 + 2/n$, $n \in \mathbb{N}$ is as above. Notice that $k > 1$ for all $n$ and since $p_0$ is greater than 1, $p_{\nu} > 1$ for all $\nu$. The construction allows $p_{\nu} \rightarrow \infty$ as $\nu \rightarrow \infty$. Likewise, when $p_0 < 0$ as will be the case for $\underset{R}{\min} \, u(x,t)$, then $p_{\nu} \rightarrow - \infty$. 

Now let us set up the domains for each of the iterations. Each iteration will be defined in its own rectangle. 

\noindent Let 
\begin{equation}
\tau _{\nu} := \tau \frac{1 + \tau ' \nu}{1+ \tau \nu},
\end{equation} 
which converges to $\tau '$ as $\nu \rightarrow \infty$. 
Also let 
\begin{equation}
r _{\nu} := r \frac{1 + r ' \nu}{1 + r \nu}
\end{equation}
which converges to $r '$ and $\nu \rightarrow \infty$. 
\noindent Let 
\begin{equation}
R_{\nu} = B_{r _{\nu}} \times ( 0, \tau _{\nu}  ).
\end{equation}
\noindent Notice that since $\nu$ starts at 0, we have
\begin{equation}
\tau _0 = \tau \frac{1 + \tau ' \cdot 0}{1 + \tau \cdot 0} = \tau
\end{equation} 
and 
\begin{equation}
\tau _1 = \tau \frac{1 + \tau '}{1 + \tau} < \tau
\end{equation}
which gives us smaller $\tau _{\nu}$ as $\nu \rightarrow \infty$. Also notice that 
\begin{equation}
r _0 = r \frac{1 + r \ \cdot 0}{1+ r \cdot 0} = r > r _1.
\end{equation}
Likewise, 
\begin{equation}
R_0 = B_0 \times ( 0, \tau _0 ) = B(x, r) \times (0, \tau)
\end{equation}
and $R_{\nu} \rightarrow R'$ as $\nu \rightarrow \infty$. 

\noindent Set 
\begin{equation}
H_{\nu} := H_{r _{\nu} \tau _{\nu}} (v_{\nu}) = \frac{1}{|R_{\nu}|} \iint_{R_{\nu}} v^{p_{\nu}} dx dt.
\end{equation}
\noindent Now we will show the iteration. Take for example, for any positive value of $\nu$
\begin{equation}
H_{\nu + 1} = H_{r _{\nu +1}, \tau _{\nu +1}} (v_{\nu +1}) = \frac{1}{|R_{\nu +1}|} \iint_{R_{\nu + 1}} v^{p_{\nu +1}} dx dt \leq \gamma H_{\nu} ^k
\end{equation}
To arrive at a constant for $\gamma$ we need to go back to the local bounds calculated earlier for higher powers of p. Using the same conventions for $H_{\nu} = H_{r _{\nu}, \tau _{\nu}}$ we have
$$M_{\nu +1} = \frac{1}{|B_{\nu +1}|} \underset{t}{\max} \int_{B_{\nu +1}} v_{\nu +1} dx = \frac{1}{|B_{\nu +1}|} \frac{|R_{\nu}|}{|R_{\nu}|} \underset{t}{\max} \int_{B_{\nu +1}} v^{p_0 k^{\nu +1} /2} dx $$

$$\leq 4p_{\nu +1} \lambda (\frac{C_{HG}^2}{(r _{\nu} - r _{\nu +1})^2} + \frac{1}{\tau _{\nu} - \tau _{\nu +1}}) \frac{|R_{\nu}|}{|B_{\nu +1}|} H_{\nu}$$

\begin{equation}
\label{mnubounds}\leq 4p_{\nu +1} \lambda (\frac{C_{HG}^2}{(r _{\nu} - r _{\nu +1})^2} + \frac{1}{\tau _{\nu} - \tau _{\nu +1}}) \bigl( C_B (\frac{r _{\nu}}{r _{\nu +1}})^Q \bigr) \tau _{\nu}  H_{\nu}
\end{equation}
Where $C_{HG}$ was the constant from the bounds on the horizontal gradient. $C_B$ was the constant from the bounds of the CC balls. 

$$D_{\nu +1} = \frac{1}{|R_{\nu +1}|} \iint_{R_{\nu +1}} |\hg v_{\nu +1}|^2 dxdt = \frac{1}{|R_{\nu +1}|} \frac{|R_{\nu }|}{|R_{\nu }|} \iint_{R_{\nu +1}} |\hg v_{\nu +1}|^2 dxdt  $$

\begin{equation}
\label{Dnubounds}\leq 4 \lambda ^2 (p_{\nu +1})^2 / (2p_{\nu} -1) (\frac{C_{HG}^2}{(r _{\nu} - r _{\nu +1})^2} + \frac{1}{\tau _{\nu} - \tau _{\nu +1}}) \bigl( C_B (\frac{r _{\nu}}{r _{\nu +1}})^Q \bigr) H_{\nu}
\end{equation}

Now substitute (\ref{mnubounds}) and (\ref{Dnubounds}) into the lemma above. This gives

$$ \leq [ 4 \lambda ^2 (p_{\nu +1})^2 / (2p_{\nu} -1) (\frac{C_{HG}^2}{(r _{\nu} - r _{\nu +1})^2} + \frac{1}{\tau _{\nu} - \tau _{\nu +1}}) \bigl( C_B (\frac{r _{\nu}}{r _{\nu +1}})^Q \bigr) H_{\nu} + H_{\nu} ] \times $$

$$[ 4p_{\nu +1} \lambda (\frac{C_{HG}^2}{(r _{\nu} - r _{\nu +1})^2} + \frac{1}{\tau _{\nu} - \tau _{\nu +1}}) \bigl( C_B (\frac{r _{\nu}}{r _{\nu +1}})^Q \bigr) \tau _{\nu}]^{\frac{2}{n}}  H_{\nu}^{\frac{2}{n}}$$

$$\leq [ 4 \lambda ^2 (p_{\nu +1})^2 / (2p_{\nu} -1) (\frac{C_{HG}^2}{(r _{\nu} - r _{\nu +1})^2} + \frac{1}{\tau _{\nu} - \tau _{\nu +1}}) \bigl( C_B (\frac{r _{\nu}}{r _{\nu +1}})^Q \bigr) + 1 ] \times $$

$$[ 4p_{\nu +1} \lambda (\frac{C_{HG}^2}{(r _{\nu} - r _{\nu +1})^2} + \frac{1}{\tau _{\nu} - \tau _{\nu +1}}) \bigl( C_B (\frac{r _{\nu}}{r _{\nu +1}})^Q \bigr) \tau _{\nu} ]^{\frac{2}{n}}  H_{\nu} ^k$$

\noindent Which is 
\begin{equation}
\label{hnuplusoneleqghnu}H_{\nu + 1} \leq \gamma H_{\nu} ^k
\end{equation}
with $\gamma$ depending on $R$, $n$, $\lambda$ only. Now we will use this fact to achieve a maximum. First we will iterate down to $H_0$. Take any $H_{\nu }$ and notice that 
$$H_{\nu} \leq \gamma _1 H^k_{\nu -1} \leq \gamma _2 ^{\nu} H^k _{\nu -1}$$ 
for some $\gamma_2 \geq \gamma^{\frac{1}{\nu}}_1.$
So $$H_{\nu } \leq \gamma _2 ^{\nu} H^k _{\nu -1} $$
Now use $(H_{\nu-1})^k$ to iterate again.
$$= \gamma _2 ^{\nu} (H_{\nu-1})^k \leq \gamma _2 ^{\nu } (\gamma_2 ^{\nu - 1} H^k_{\nu - 2})^k$$

$$= \gamma _2 ^{(\nu ) + (\nu -1)k }  H^{k^2}_{\nu - 2}$$

$$\leq \gamma _2 ^{(\nu ) + (\nu -1)k }  (\gamma _2 ^{(\nu - 2)}H^k_{\nu - 3})^{k^2}$$

$$= \gamma _2 ^{(\nu ) + (\nu -1)k + (\nu -2)k^2} H_{\nu - 3}^{k^3}$$

$$\leq \gamma _2 ^{(\nu) + (\nu -1)k + (\nu -2)k^2 + ... + (1)k^{\nu - 1}} H_{0}^{k^{\nu}}$$

$$\leq \gamma _3 ^{k^{\nu}} H_{0}^{k^{\nu}}$$ 
for some 
$$\gamma_3 \geq \gamma_2^{\frac{\nu + (\nu -1)k + ... + (1)k^{\nu -1}}{k^{\nu}}}$$
So we have 
$$H_{\nu} \leq \gamma_3^{k^{\nu}} H_{0}^{k^{\nu}}$$
which in terms of the integrals says that
$$ \frac{1}{|R_{\nu}|} \iint_{R_{\nu}} v^{p_0 k^{\nu} / 2 } dxdt \leq \gamma_3^{k^{\nu}} \bigl( \frac{1}{|R|} \iint_{R} v^{p_0 / 2 } dxdt \bigr)^{k^{\nu}} $$
Then taking $k^{\nu}$ roots: 
\begin{equation}\label{before_sup_norm}\bigl( \frac{1}{|R_{\nu}|} \iint_{R_{\nu}} |v^{p_0 / 2 }|^{k^{\nu}} dxdt \bigr)^{1/k^{\nu}} \leq \gamma_3  \frac{1}{|R|} \iint_{R} v^{p_0 / 2 } dxdt. 
\end{equation} 
%\noindent The LHS tends to the maximum of $v^{p_0/2}$ in $R'$, 
We now use the following fact from $L^p$ theory:
$$\underset{p \rightarrow \infty}{\lim} \ \ (\int_{R} v^p \, dxdt)^{1/p} = ||v||_{\infty} = \underset{R}{ess \ sup} \ u.$$
So on the LHS of (\ref{before_sup_norm}) we have
$$ \underset{R'}{\max} \, v \, = \, \underset{\nu \rightarrow \infty}{\lim} ( \frac{1}{|R_{\nu}|} \iint_{R_{\nu}} |v^{p_0 / 2 }|^{k^{\nu}} dxdt)^{1/k^{\nu}} $$
So we conclude:
\begin{equation} 
\label{max}\underset{R'}{\max} \, u \, \leq \gamma \bigl( \frac{1}{|R|} \iint_{R} u^{p } dxdt \bigr)^{1/p}.
\end{equation}

We will consider the case for $p < 0$ in order to achieve a similar result for the minimum of $v^{p_0/2}$ in $R'$, where $v$ is a weak supersolution. So assume that $p_{\nu} < p_0 < 0$ where $-p_0$ is the same value as above. Let $p_{\nu} = p_0 k^{\nu}/2$. Then $\underset{{\nu} \rightarrow \infty}{\lim} p_{\nu} = - \infty$. Notice also that at $\nu = 0$, we have $p_{\nu = 0} = -p_0/2$ as before. Each of the calculations are the same as the $0 < p_0 < p_{\nu} $ case except that we have 
$$H_{\nu +1} \geq \gamma H_{\nu}$$
since $v$ is a supersolution we have:
$$\iint \phi v_t \, dxdt+ \sum_{i,j = 1}^{m} \iint X_i \phi (a_{ij}(x,t) X_j v) \, dxdt\geq 0.$$
\noindent Then by a similar argument we have 
\begin{equation}
\label{hnuGEQghzero}H_{\nu} \geq \gamma ^{k^{\nu}} H_0^{k^{\nu}}
\end{equation}

\noindent  and taking $k^{\nu}$ roots the RHS of (\ref{hnuGEQghzero}) tends to the minimum of $v^{p_0/2}$ in $R'$,
$$\underset{\nu \rightarrow \infty}{\lim} ( \frac{1}{|R_{\nu}|} \iint_{R_{\nu}} |v^{p_0 / 2 }|^{k^{\nu}} dxdt)^{1/k^{\nu}} \, = \, \underset{R'}{\min} \, v$$
where $p_0<0$.

\noindent So we conclude:
$$\gamma \bigl( \frac{1}{|R|} \iint_{R} u^{-p } dxdt \bigr)^{1/-p} \, \leq \, \underset{R'}{\min} \, u $$

\noindent Now it is only left to establish an inequality that will bridge the gap between positive and negative $p_0$,
$$ \frac{1}{|R_{-}|} \iint_{R_{-}} v^{p_0 k^{\nu} / 2 } dxdt \leq \gamma \frac{1}{|R^{+}|} \iint_{R^{+}} v^{-p_0 k^{\nu} / 2 } dxdt $$

$$  $$
 %-------------------------------------------------------------
% You cannot use the \chapter{} or \section{} commands because they put a header that the grad school doesn't allow
\begin{center} 
\textbf{1.5 Bridge step}
\end{center}
%-------------------------------------------------------------

Above we established appropriate estimates on the maximum and minimum of a positive, weak solution to (\ref{heat}) in a rectangle $R$. We turn to the gap for $-p_0 < 0 < p_0$. Our earlier results gave us (\ref{max})
$$ \underset {R'}{\max} \, u^p \, \leq \gamma \iint_R u^p dxdt$$
so we start with $v = u^{- \epsilon}$ which was proven to be a positive subsolution under the $p<0$ case. Applying the results of (\ref{maxuleqpnorm}) above we get

$$ \underset{R'}{\max} \, u^{-\epsilon}  \, \leq \gamma ( \iint_R u^{- \epsilon p} dxdt)^{1/p }$$

\noindent Now taking the $-1 / \epsilon $ root  and setting $p=2$ we find

$$ \underset{R'}{\min}  \, u \, \geq \, \gamma^{-1/ \epsilon} ( \iint_R u^{-2 \epsilon } dxdt)^{-1/2 \epsilon } $$

\noindent Let $\gamma^{-1 / \epsilon} = \gamma_1$ and use $\epsilon$ for $2 \epsilon$
\begin{equation}
\label{min_geq_M_ep_R}\underset{R'}{\min}  \, u \, \geq \, \gamma_1 ( \iint_R u^{- \epsilon } dxdt)^{-1/ \epsilon } 
\end{equation}
Introduce the notation for some domain $D \subset R$
\begin{equation}
\label{Mpd} M(p,D) := (\frac{1}{|D|} \int_{D} u^p \, dxdt )^{1/p}
\end{equation}
and use this for (\ref{min_geq_M_ep_R}) which gives us 
\begin{equation}
M(- \infty, R') = \underset{R'}{\min} \, u \, \geq \, \gamma_1 ( \iint_R u^{- \epsilon } dxdt)^{-1/ \epsilon } = \gamma_1 M(- \epsilon , R). 
\end{equation}
We also want to establish for some constant $\gamma ^*$
\begin{equation}
\label{max_LEQ_M_epsilon}\underset{R^{*}}{\max \, u} = M( \infty, R^{*}) \, \leq \, \gamma^{*} ( \iint_R u^{ \epsilon } dxdt)^{1/ \epsilon } = \gamma^{*} M( \epsilon , R). 
\end{equation}
So we let $R^{*} \subset R^{-}$ where the boundary of $R^{*}$ does not touch the top of $R^{-}$. Now apply the results of (\ref{maxuleqpnorm}) to $u$ in $R^*$ with $p' > 1$, then for some $\gamma'$: 
\begin{equation}
\label{maxRstarLEQ_M_pprime_Rplus}\underset{R^*}{\max} \, u \leq \ \gamma' M(p', R^-)
\end{equation}
Recall that we proved the local bounds:
$$\iint |\hg v|^2 dxdt\leq C_1 \iint v^2dxdt$$ and 
$$\underset{t \in (0, \tau)}{\max} \int  v^2 dx\leq C_2 \iint v^2dxdt$$ 
for $0 < p < 1,$ in the form $v = u^{p/2}(-t, x)$. We create an iteration argument like the one above. This gives us, for some $\gamma$ the inequality
$$({1 \over |R^{-}|} \iint_{R^-} u^{pk}dxdt )^{1/pk} \leq \gamma({1 \over |R|} \iint_{R} u^{p} dxdt)^{1/p}$$ like before. Now let $pk = p'$ and $p = {p' \over k}.$ This gives

$$({1 \over |R^{-}|} \iint_{R^-} u^{p'} dxdt)^{1/p'} \leq \gamma({1 \over |R |} \iint_{R } u^{p'/k} dxdt)^{k/p'}$$ 
The issue here is that the $p'$ will stay fixed and we will iterate the RHS to achieve ${p' \over k} \rightarrow \epsilon$ with larger $k$. Using the $k=1 + 2/n$ from the previous iteration, we set $p' = {k +1 \over 2}$ and $p_{\nu} = {1 \over k^{\nu}} {k+1 \over 2}$. This gives

$$M(p', R^-) \leq \gamma_{\nu} \ M(p_{\nu}, R)$$
for some $\gamma_{\nu}$.
Now let $\nu$ be the smallest integer such that  
\begin{equation}
\label{p_required_for_bridge} p_{\nu} = {1 \over k^{\nu}} {k+1 \over 2} < \epsilon
\end{equation}
We find
\begin{equation}
\label{M_p_prime_Rplus_leq_M_ep_R}M(p', R^-) \leq \gamma M(\epsilon, R)
\end{equation}
Combining (\ref{M_p_prime_Rplus_leq_M_ep_R}) and (\ref{maxRstarLEQ_M_pprime_Rplus}) gives 
$$\underset{R^{*}}{\max \, u} \leq \ \gamma'  M(p', R^-) \leq \ \gamma \gamma' M( \epsilon , R)$$ which proves (\ref{max_LEQ_M_epsilon}) with $\gamma ^* = \gamma \gamma'$. Take this along with
$$
M(- \epsilon , R) \leq \, \gamma \, \underset{R'}{\min} \, u
$$
proved earlier. To finish the proof of the Harnack inequality, we need to bridge the gap between $M(\epsilon, R^*)$ and $M(- \epsilon , R')$. In order to do this we will use the results of Hugo Aimar \cite{Ai}. 
 %-------------------------------------------------------------
% You cannot use the \chapter{} or \section{} commands because they put a header that the grad school doesn't allow
\begin{center} 
\textbf{1.6 Lag mapping}
\end{center}
%-------------------------------------------------------------

First, some definitions. 

\begin{defn}[Lag Mapping] Let $R(x_0, r_0)$ be a fixed ball in the space of homogeneous type $(X, d, \mu)$. We use $X=\mathbb{R}^{n+1}$. We shall say that a function
$$T: R(x_0, r_0) \times (0, r_0] \rightarrow X \times R^{+}; \, \, \, T(x,r) = ( \xi  , \rho),$$
is a lag mapping if there exist three constants $K_i$, $i=1,2,3$ such that the inequalities
\begin{equation}
\label{K_for_dist_x_xi_leq_Kr}d_{cc}(x, \xi) \leq K_1 r,
\end{equation}
\begin{equation}
\label{lag_bounds_on_rho}K_2 \rho \leq r \leq K_3 \rho
\end{equation}
hold for every $x \in R(x_0, r_0)$ and $r \in (0, r_0]$.
\end{defn}
\noindent In the following we will consider rectangles 
$$R((x,t),r) = R(x,r)= B(x,r) \times (t -r, t+r).$$
In the definition of BMO we shall use non-negative functions $h$ of a real variable $t$, which are continuous, non-negative, $h(t) =0$ if $t<0$, $h(t)$ is increasing for $t>0$, $h(t+s) \leq h(t) + h(s)$ and $e^{-\epsilon h(t)}$ is an integrable function on $(0, \infty)$ for every $\epsilon > 0$.

\begin{defn} [BMO condition] Let $f$ be a real-valued measurable function on the space of homogeneous type $(X, d, \mu)$. We shall say that $f$ satisfies a BMO condition with lag mapping $T$ on the rectangle $R(x_0, r_0)$, with respect to $h(s)= \sqrt{s} $ if there exists a real function $C(x,r)$ on $R(x_0 r_0) \times (0, R_0]$ such that

\begin{equation}
\label{averageh_F-C_leqconst}{1\over |R(x,r)|} (\int_{R(x,r)}\sqrt{(f-C(x,r))^+} \ dxdt) \leq N(f),
\end{equation}
\begin{equation}
\label{averageh_C-F_leqconst}{1\over |R(x,r)|} (\int_{R(x,r)}\sqrt{(C(x,r)-f)^+} \ dxdt) \leq N(f),
\end{equation}
for some constant $N(f)$ and every $(x,r) \in R(x_0, r_0) \times  (0, R_0]$. Let $BMO(x_0, r_0, T, h)$ be the class of all such functions. 

\end{defn}

\noindent The main result of the paper and the one we will make use of is the 

\begin{thm} [John Nirenberg for BMO with Time Lag] Let $T$ be a one to one lag mapping on $R(x_0, r_0)$ with the following property: there exists $0 < \gamma < 1$ such that for every $r \leq r_0$
\begin{equation}
\label{ball_subset_lag_of_ball}R(x_0, \gamma r) \times (0, \gamma r] \subset T(R(x_0) \times (0,r]).
\end{equation}
 
Then there exist two constants $\delta$ and $C$ depending on $K_i$, $A$, $\gamma$ and a lag mapping $S$ on $B(x_0, \delta r_0) \, \, (S(x,r) = (\zeta, \tau))$ such that if $f \in BMO(x_0, r_0, T, h)$ and $u=e^{-f}$, the inequality 
\begin{equation}
\label{bridge_inequality_in_aimar}(m_{R(x,r)} (u^{- \epsilon})) \cdot (m_{R(\zeta,\tau)} (u^{ \epsilon})) \leq C
\end{equation}
holds for some $\epsilon > 0$, depending only on $N(f)$, and every $(x,r) \in R(x_0, \delta r_0) \times (0, \delta r_0]$. Moreover, the mapping S can be given explicitly in terms of $T$: 
$$ S(x,r) = (\xi (x, \frac{r}{ \eta}); \eta \rho (x, \frac{r}{\eta})),$$
where $\xi$ and $\rho$ are the components of $T$ and $\eta$ is a constant depending only on $K_i$, $i=0,1,2,3$.
\end{thm}
We will use this in the following way. The bridge part of the theorem is from (\ref{bridge_inequality_in_aimar}) and takes the form
\begin{equation}\label{bridge_inequality}
{1 \over |R((x,t),r)|} \iint_{R((x,t),r)} u^{- \epsilon} dxdt \cdot {1 \over |S(R)|} \iint_{S(R)} u^{ \epsilon} dxdt \leq C
\end{equation}
Where $$S((x,t),r) = (x, t-2r^2/ \eta ^2),r)$$

\noindent Now we want to show that this theorem applies to our setting. First our space $(X, d, \mu)$ is a space of homogenous type with 
$$X= \mathbb{R}^{n+1}$$ 
$$d((x,t),(y,s)) = d_{cc}(x,y) + |t-s|^{1/2}$$
and $\mu = (n+1)$-dimensional Lebesgue measure. We will use two notations: $(x,t)\in X \subset \mathbb{R}^{n+1}$ or $x \in X \subset \mathbb{R}^{n+1}$ depending on whether we want to draw attention to the time dimension of a point.

The constants $K_2 \rho \leq r \leq K_3 \rho$ are both 1. Now we need to show that the doubling property holds. We want an $A$ such that 

$$0 < \mu (R(x,2r)) < A \mu (R(x,r)) $$ 
holds for $r>0$ and every $(x,t) \in X$. An earlier proposition [\ref{citti}] gives us 

$$|B(x, tr)| \geq C t^Q |B(x, r)|$$
for every $x \in U$ bounded, $r \leq R_0$ and $0 < t < 1$. So set $t=1/2$ and use

$$|B(x, {1 \over 2}r)| \geq C_B ({1 \over 2})^Q |B(x, r)|$$
$$|B(x, {1 \over 2}r)| \times 2r \geq C_B ({1 \over 2})^Q |B(x, r)| \times 2r$$
$$|B(x, {1 \over 2}r) \times ({-r \over 2}, {r \over 2})| \geq C_B ({1 \over 2})^{Q+1} |B(x, r) \times (-r,r)|$$
$$|R(x,{r \over 2})| \geq C_B ({1 \over 2})^{Q+1} |R(x,r)|$$
Which gives the desired inequality for $A > {1 \over C_B}(2)^{Q+1}$

We will use the lag map $T:\mathbb{R}^n \times \mathbb{R} \times [0, \infty) \rightarrow \mathbb{R}^n \times \mathbb{R} \times [0, \infty)$. This map acts on pairs from $\mathbb{R}^{n+1}$ and $\mathbb{R}^+$ in the form $T((x,t), r)$ where $(x,t)$ is a point in space-time and $r$ is a radius for the rectangle in $\mathbb{R}^{n+1}$. Define $T((x,t); r):= ((x, t-2r^2);r) = ((\xi, s ), \rho)$.  We want to check that $T$ is a lag map. Notice that 
$$d((x,t), (\xi, \rho)) = d_{cc}(x, x) + |t-\rho |^{1/2} = 0 + |t -(t-2r^2)|^{1/2} = r \sqrt{2}.$$
Last we check that there exists $0 < \gamma < 1$ such that for every $r \leq R_0$
\begin{equation}
R(x_0, \gamma r) \times (0, \gamma r] \subset T(R(x_0,r) \times (0,r]).
\end{equation}
$$R(x_0, \gamma r) \times (0, \gamma r] = B_{cc}(x, \gamma r) \times (t -\gamma r,t+ \gamma r ) \times (0, \gamma r].
$$
$$
T(R(x_0, r) \times (0, r]) = B_{cc}(x, r) \times (t -2r^2 -r,t+r -2r^2) \times (0, \gamma r]
$$
Clearly, $B_{cc}(x, \gamma r) \subset B_{cc}(x, r)$ and $(0, \gamma r] \subset (0, r]$ for every $0 < \gamma < 1$. And  $(t -\gamma r,t+ \gamma r ) \subset (t -2r^2 - r,t+ r -2r^2 )$ for $\gamma < 1$.

 %-------------------------------------------------------------
% You cannot use the \chapter{} or \section{} commands because they put a header that the grad school doesn't allow
\begin{center} 
\textbf{1.7 Estimates on log(u)}
\end{center}
%-------------------------------------------------------------

\noindent Now we calculate estimates on $v:=- \log u$. We will use the results of these estimates to show that $v \in BMO(x_0, r_0, T, h)$ in order to use the results (\ref{bridge_inequality_in_aimar}) of Aimar which will conclude the proof. We will consider positive supersolutions $u \geq \epsilon > 0$ in $B(x,r) \times (-r < t < r)$. Let $v(x,t):= - log(u)$. Recall from the definition of weak supersolution to $$\parab$$ in the region $R$ defined before:
\begin{equation}\label{start_of_logu}
\iint_{R} \phi u_t \ dxdt+  \sum^{m}_{i,j=1}\iint_{R}(a_{ij}(x,t) X_j u) X_i\phi \ dxdt \ \geq 0.
\end{equation}
We will use $\phi = \psi ^2 u^{-1}$ where $\phi(x,t) \geq 0$ is of compact support in the $x$ variable and $\psi$ is independent of $t$ and of compact support. 

\noindent Then we have with 
$$v_t = (-log(u))_t = -\frac{u_t}{u}$$
$$X_i v = X_i(-log(u)) = \frac{-X_i u}{u}$$
$$X_i \phi = 2 \psi X_i \psi u^{-1} - u^{-2}X_i u \psi^2$$

\begin{equation}
\iint_{R} \psi^2 \frac{u_t}{u} \   + <2 \psi \hg \psi u^{-1}, a \hg u>  \ dxdt 
\end{equation}
$$
- \iint_{R} \ <\frac{\hg u}{u^2} \psi^2, a \hg u> dxdt \, \geq 0.
$$
And substituting for $v_t$ and $\hg v$
$$\iint_{R} -(\psi ^2 v_t) - <\hg(\psi^2), a \hg v> - <\psi^2 \hg v, a \hg v> dxdt \,  \geq 0$$
$$\iint_{R} \psi ^2 v_t + < \hg (\psi^2), a \hg v> + <\psi^2 \hg v, a \hg v> dxdt \,  \leq 0$$
Now consider the integral in $t$ from $t_1$ to $t_2$ in the interval:
$$\int_{B(x,r)} \psi^2 v  \ dx|_{t_1}^{t_2}  + \int_{t_1}^{t_2} \int_{B(x,r)} 2\psi<\hg \psi, a \hg v> + \psi ^2<\hg v, a \hg v> \ dxdt \, \leq 0.$$
And taking to the RHS the second integrand
$$\int_{B(x,r)} \psi^2 v dx|_{t_1}^{t_2}  + \int_{t_1}^{t_2} \int_{B(x,r)} \psi ^2<\hg v, a \hg v> dxdt  \, $$
$$\leq \,  -2 \int_{t_1}^{t_2} \int_{B(x,r)} \psi<\hg \psi, a \hg v> dxdt .$$
\begin{equation}
\label{int_t1_to_t2_for_schwarz} \int_{B(x,r)} \psi^2v  \ dx |_{t_1}^{t_2} + \int_{t_1}^{t_2} \int_{B(x,r)} \psi ^2<\hg v, a \hg v> \ dxdt  
\end{equation} 
$$
\leq \,  |2 \int_{t_1}^{t_2} \int_{B(x,r)} \psi<\hg \psi, a \hg v>| dxdt .
$$
Now we look at the second term in the integral and use Schwarz' and Young's inequality:
$$2\int_{t_1}^{t_2} \int_{B(x,r)} \psi <\hg \psi, a \hg v> dxdt \ $$

$$\leq 2(\int_{t_1}^{t_2} \int_{B(x,r)} \psi^2 <\hg v, a \hg v> dxdt $$

$$ \times \int_{t_1}^{t_2} \int_{B(x,r)} <\hg v, a \hg v> \ dxdt)^{1/2} $$

\begin{equation}
\label{schwarzboundsonpsi}\leq 2 ( \frac{1}{4} \int_{t_1}^{t_2} \int_{B(x,r)} \psi^2 <\hg v, a \hg v> dxdt 
\end{equation}
$$
+  \int_{t_1}^{t_2} \int_{B(x,r)} <\hg v, a \hg v> dxdt ) 
$$
When we combine (\ref{int_t1_to_t2_for_schwarz}) with (\ref{schwarzboundsonpsi}) we get:
\begin{equation}
\label{preliminary_before_poincare}\int_{B(x,r)} \psi^2 v dx |_{t_1}^{t_2} + \frac{1}{2} \int_{t_1}^{t_2} \int_{B(x,r)} \psi ^2 <\hg v, a \hg v> dxdt \, 
\end{equation}
$$
\leq \, 2 \int_{t_1}^{t_2} \int_{B(x,r)} <\hg \psi, a \hg \psi> dx dt
$$
Now we make use of Jerison's Poincar\'e inequality for H\"ormander vector fields \cite{Jer}. We will point out that Jerison assumes $v$ is $C^{\infty}(B(x,r))$. In order to improve this result for functions $v \in S^{1,2}  (B(x,r))$, we will use the fact that the smooth functions are dense in $S^{1,2} (B(x,r))$. I.e. given a $v \in S^{1,2} (B(x,r))$, there exists a sequence of smooth functions on the ball, $\{ v_n \}_{n=1}^{\infty} \in C^{\infty}(B(x,r))$, such that $$\underset{n \rightarrow \infty}{\lim} v_n = v$$ and $$\underset{n \rightarrow \infty}{\lim} \hg v_n = \hg v$$ in the $S^{1,2}$ norm.
We will let $\Omega$ be a neighborhood of the closure of the unit ball $B(x,1)$, in $\mathbb{R}^n$.  

\begin{lem}[Poincar\'e  inequality]
Let $\{ X_i \}^{m}_{i=1}$, $m \leq n$, be smooth vector fields satisfying H\"ormander's condition. Then there is a positive constant $C_P$ and a radius $r_0  > 0$ such that for every $x \in B(x,1)$ and every $0 < r < r_0$ for which $B(x, 2r) \subset \Omega$ 

\begin{equation}
\int_{B(x,r)} |v -V_B|^2 \ dx \,  \leq C_Pr^2 \int_{B(x,r)} |\hg v|^2  \ dx
\end{equation}
holds for every function, $v \in S^{1,2}(\overline{B(x,r)} \times (0, \tau))$ with fixed $t$. Here $V_B(t)= {1 \over |B(x.r)|} \int v(x,t) dx $. 

\end{lem}

%\noindent We will let $$\psi = \psi (x) = \prod _{\nu = 1}^n \chi _{\nu} (x_{\nu})$$ and  
%$$
%\begin{equation*}
%\chi _{\nu} (x_{\nu}) = \left\{
% \begin{array}{rl}
%  1 & \text{if } |x_{\nu} - x_{\nu}^0| < r \\
%   0 & \text{if } |x_{\nu} - x_{\nu}^0| \geq 2r \\
% \end{array} \right.
%\end{equation*}
%$$

%\noindent and is linearly interpolated otherwise. Here we use $(x_{1}^0, x_{1}^0, %...x_{n}^0)$ to be the center of the ball in space. 
%--------------------------------------------------------
\noindent Now we construct a test function whose existence is guaranteed by Lemma (\ref{CC_test_functions}) \cite{CiL}. Let $\psi(x)$ be defined as follows: 
$$
%\begin{equation*}
\psi (x) = \left\{
 \begin{array}{rl}
  1 & \text{if } d_{cc}(x, 0) < r \\
   0 & \text{if } d_{cc}(x, 0) \geq 2r \\
 \end{array} \right.
%\end{equation*}
$$
We note that $|\hg \psi | \leq C_{HG}/r^2$. Now we recall the bounds on $a_{ij}(x,t)$, integrating in $t$ on the RHS since $\psi$ is independent of t and use (\ref{preliminary_before_poincare}). 
$$\int_{B(x,r)} \psi ^2 v dx \, |_{t_1}^{t_2} \, + \, \frac{1}{2 \lambda} \int_{t_1}^{t_2} \int_{B(x,r)} |\hg v|^2 \psi ^2 dxdt \, \leq 2 \lambda (t_2 - t_1) \int_{B(x,r)} |\hg \psi|^2 dx.$$
%%%%%%%%%%%%%%%%%%%%%%%%%%%%%%%%%%%%%%%%%%%%%%%%%%%%%%%%%%%%%%%%%%%%%%%%%%%%%
%%%%%%%%%%%%%%%%%%%%%%%%%%%%%%%%%%%%%%%%%%%%%%%%%%%%%%%%%%%%%%%%%%%%%%%%%%%%
%%% adding the alternative argument here. ending %%%%%%%%%%%%%%%%%%%%%%%%%
%%%%%%%%%%%%%%%%%%%%%%%%%%%%%%%%%%%%%%%%%%%%%%%%%%%%%%%%%%%%%%%
We observe that we defined $\psi(x) = 1$ on $B(x,r)$ so this leads to
\begin{equation}\label{before_assuming_diffblty_of_VB}
\int_{B(x,r)}  v dx \, |_{t_1}^{t_2} \, + \, \frac{1}{2 \lambda} \int_{t_1}^{t_2} \int_{B(x,r)} |\hg v|^2  dxdt \, \leq 2 \lambda (t_2 - t_1) \int_{B(x,r)} |\hg \psi|^2 dx.
\end{equation}
Substitute the value  $V_B = {1 \over |B(x,r)|}\int_{B(x,r)} v dx$ and divide by $(t_2 - t_1)$,
$$({V_B(t_2) - V_B(t_1) \over (t_2 - t_1)})|B(x,r)| + \, \frac{1}{2 \lambda} \int_{t_1}^{t_2} \int_{B(x,r)} |\hg v|^2  dxdt \, \leq 2 \lambda \int_{B(x,r)} |\hg \psi|^2 dx.$$
Now use the Poincar\'e lemma
$$({V_B(t_2) - V_B(t_1) \over (t_2 - t_1)})|B(x,r)| + {C_Pr^2 \over 2 \lambda (t_2 - t_1)} \int_{t_1}^{t_2} \int_{B(x,r)} (v - V_B)^2 dxdt \leq 2 \lambda \int_{B(x,r)} |\hg \psi|^2 dx.$$
On the RHS we have $|\hg \psi| \leq 1 /r$ so we get 
$$({V_B(t_2) - V_B(t_1)  \over (t_2 - t_1)})|B(x,r)| + {C_Pr^2 \over 2 \lambda (t_2 - t_1)} \int_{t_1}^{t_2} \int_{B(x,r)} (v - V_B)^2 dxdt \leq {2 \lambda \over r^2} \int_{B(x,r)}  dx$$
Now divide through by $|B(x,r)|$ 
$${V_B(t_2) - V_B(t_1) \over (t_2 - t_1)}  + {C_Pr^2 \over 2 |B(x,r)| \lambda (t_2 - t_1)} \int_{t_1}^{t_2} \int_{B(x,r)} (v - V_B)^2 dxdt \leq {2 \lambda \over |B(x,r)| r^2} |B(x,r)|  $$
Take $t_1 \rightarrow t_2$. We note that $V_B$ need not be differentiable although we will assume so for the present argument. We present a remark (\ref{steklov_proof}) at the end of the proof over an alternative technique using the Steklov average. This allows us not to assume the $t$ derivative of $V_B$ exists.  
\begin{equation}
\label{dV_plus_poincare_LEQ_const}\frac{dV_B}{dt} + {C_Pr^2 \over 2 |B(x,r)|\lambda} \int_{B(x, r)} (v(x,t) - V_B)^2 dx \leq {2 \lambda \over r^2}   =: m 
\end{equation}
Let $w(x,t) = v(x,t) - mt$ and $W = V_B(t) - mt$. Then
$${dW \over dt} = {dV_B \over dt } -m$$ and 
$$v(x,t) - V_B = w(x,t) + mt -(W + mt) = w(x,t) - W$$
Substituting this back into (\ref{dV_plus_poincare_LEQ_const}) we get
\begin{equation}
{dW \over dt} + {C_Pr^2 \over 2|B(x,r)| \lambda} \int_{B(x, r)} (w(x,t) - W)^2 dx \leq 0. 
\end{equation}
\begin{equation}
{dW \over dt} \leq - {C_Pr^2 \over 2 |B(x,r)| \lambda} \int_{B(x, r)} (w(x,t) - W)^2 dx 
\end{equation}
Let $$B_s(t) = \{ x \in B(x,r) : w(x,t) > s \}$$
Then $(s-W)^2 < (w-W)^2$, hence
\begin{equation}
{dW \over dt} \leq - {C_Pr^2 \over 2 |B(x,r)| \lambda} \int_{B_s(t)} (s - W)^2 dx  
\end{equation}
\begin{equation}
{dW \over dt} \leq - {C_Pr^2 \over 2|B(x,r)| \lambda} (s - W)^2 |B_s(t)|  
\end{equation}
\begin{equation}
{dW \over dt}(s - W)^{-2} \leq - {C_Pr^2 \over 2|B(x,r)| \lambda}  |B_s(t)|  
\end{equation}
We integrate from $t_1$ to $t_2$:
\begin{equation}
\int_{t_1}^{t_2}{dW \over dt}(s - W)^{-2} \leq - {C_Pr^2 \over 2|B(x,r)| \lambda}  |B_s(t) \times (t_1, t_2)|  
\end{equation}
%-------------------
$$-({1 \over (s - W(t_2))} - {1 \over (s - W(t_1))}) \leq - {C_Pr^2 \over 2|B(x,r)| \lambda}  |B_s(t) \times (t_1, t_2)| $$
$${1 \over (s - W(t_2))} - {1 \over (s - W(t_1))} \geq  {C_Pr^2 \over 2|B(x,r)| \lambda}  |B_s(t) \times (t_1, t_2)| $$
We defined $W(t) = V_B(t) - mt$ and from (\ref{dV_plus_poincare_LEQ_const}) we see that 
\begin{equation}
\label{bounds_on_dv}{dV_B \over dt}  \leq {2 \lambda \over r^2}   = m.
\end{equation}
So we can say $W(t) \leq 0$ on $(t_1, t_2)$. So that leads us to say:
$$ {1 \over s} \geq {1 \over (s - W(t_2))} \geq {1 \over (s - W(t_2))} - {1 \over (s - W(t_1))}.$$
Hence
\begin{equation}
\label{prelim_bounds_on_s-W}|B_s(t) \times (t_1, t_2)|  \leq {2|B(x,r)| \lambda \over C_Pr^2 s}  =: {c_3 \over r^2 s}
\end{equation}

%-----------------
%\begin{equation}
%\label{measure_of_ball_before_v-V}{-1 \over (s - W)} \leq - {C_Pr^2 \over 2|B(x,r)| \lambda}  |B_s(t) \times (t_1, t_2)| 
%\end{equation}
%\begin{equation}
%\label{prelim_bounds_on_s-W}{2|B(x,r)| \lambda \over C_Pr^2 (s - W)} \geq  |B_s(t) %\times (t_1, t_2)| 
%\end{equation}
%Now we would like to replace the $(s-W)$ term by something smaller and not dependent on $t$. 

%Since the second term on the LHS of (\ref{dV_plus_poincare_LEQ_const}) is greater than or equal to zero we can say that
%\begin{equation}
%\label{bounds_on_dv}{dV_B \over dt}  \leq {2 \lambda \over r^2}   = m
%\end{equation}
%Let $t_0 \in \mathbb{R}$, $x_0 \in \mathbb{R}^n$ and $0 < r < 1$ be fixed. Set $t_1 = t_0 - r^2$, $t_2 = t_0 + r^2$ and consider $t \in (t_1, t_2)$.
%We set $V - V_1 := V_B(t) - V_B(t_1)$ and notice that integrating (\ref{bounds_on_dv}) from $t_1$ to $t$ gives $V_B(t) - V_B(t_1) \leq m(t - t_1)$. Putting these together we get 
%$$V_B(t) - V_B(t_1) \leq m(t - t_1) \leq m(t_2 - t_1) \leq 2m r^2 .$$ And hence
%\begin{equation}
%\label{bounds_on_V_of_t}V_B(t)  \leq V_B(t_1) + 2m r^2
%\end{equation}
%We defined $W(t)=V_B(t) - mt$ so using this and (\ref{bounds_on_V_of_t}) together with (\ref{prelim_bounds_on_s-W}) to get
%$$|B_s(t) \times (t_1, t_2)| \leq {2|B(x,r)| \lambda \over C_Pr^2 (s - W)}$$
%$$= {2|B(x,r)| \lambda \over C_Pr^2 (s - (V_B(t)-mt))} \leq {2|B(x,r)| \lambda \over C_Pr^2 (s - ((2mr^2 + V_1)+mt))}$$
%\begin{equation}
%\label{bounds_on_rectangle_in_logu} \leq {2|B(x,r)| \lambda \over C_Pr^2 s }  =: {c_3 \over r^2 s}
%\end{equation}
\noindent Set  
$$\widetilde {R}((x,t),r) = B(x,r) \times (t - r^2, t+ r^2)$$
So we can say that
\begin{equation}
\label{bounds_on_measure_of_rectangle}\mu \{ (x,t) \in \widetilde {R}((x_0, t_0),r) : w > s \} \leq {c_3 \over r^2 s}
\end{equation}
for every $s > 3mr^2$. We let $\mu$ be the $n+1$ dimensional Lebesgue measure on the rectangle. We will prove this for $w$ since we defined $w = v -mt$. 
So
\begin{equation}
\int_{\widetilde {R}((x_0,t_0),r)} \sqrt{(v(x,t) - V_B)^{+}} d \mu (x,t) \leq \int_{\widetilde {R}((x_0,t_0),r)} \sqrt{w^{+}} d \mu (x,t) 
\end{equation}
Use the Lebesgue-Stieltjes integral: 
$$\int_{\widetilde {R}((x_0,t_0),r)} \sqrt{w^{+}} d \mu (x,t) = - {1 \over 2} \int_{3mr^2}^{\infty} s^{-1/2} d \mu (s)$$
And integration by parts
$$= -{1 \over 2} \bigl[ \mu (s) s^{-1/2} |_{3mr^2}^{\infty} - {1 \over 2} \int_{3mr^2}^{\infty} s^{-3/2} \mu (s) ds \bigr]$$
with the boundary terms equal to 0.
$$= {1 \over 4} \int_{3mr^2}^{\infty} s^{-3/2} \mu (s) ds$$
Now we use the bounds on $\mu$:
$$\leq {1 \over 4} \int_{3mr^2}^{\infty} s^{-3/2} {c_3 \over r^2 s} ds$$

\begin{equation}
\label{convergingintegral}= {1 \over 4} {c_3 \over r^2} \int_{3mr^2}^{\infty} {1 \over s^{5/2}} ds
\end{equation}
Which converges since 
$$
C \int_{3mr^2}^{\infty} {1 \over s^{5/2}} ds \ \leq \ C \int_{\mathbb{R}} {1 \over s^{5/2}} ds \  < \infty.
$$
This leads us to say
$$\int_{\widetilde {R}((x,t),r)} \sqrt{(v(x,t) - V_B)^{+}} d \mu (x,t) \leq \, c_4$$
Now notice that for any slice of time $t_1$ to $t_2$ we have $|t_1 - t_2| = 2r^2$. If we replace $t$ with $2t_1-t$ in the above we find 
$$t_1 \leq 2t_1 - t \leq t_2$$
$$-t_1 \leq -t \leq t_2 -2t_1$$
$$t_1 \geq t \geq 2t_1-t_2  = 2t_1 -(t_0 + r^2)$$
and that gives $$|(t_1) - [2t_1 -(t_0 + r^2)]| = |t_0 - t_1 -r^2| = 2r^2$$
This is the correct size of the time slice with $t_1$ as the upper and $t_2$ as the lower bounds. Now we consider $2t_1 -t \in (t_2, t_1)$. We take the argument over at (\ref{bounds_on_dv}).
$${dV_B \over dt}  \leq {2 \lambda \over r^2}   = m$$ 
and we integrate from $t_1$ to $t$ to get
\begin{equation}
\label{negative_t_bounds_on_VTone_minus_VB}V_B(t_1) - V_B \leq m[(2t_1 - t) -t_1 ] \leq 2mr^2.
\end{equation}
Now for $s > 3mr^2 $ set $$B_s(t) = \{ x \in B(x,r) : V_B - v(x,t) > s  \}$$
which implies $V_B(t_1) - s > v(x,t)$
Putting these requirements together with (\ref{negative_t_bounds_on_VTone_minus_VB}) gives us:
$$V_B - v(x,t) > V_B -(V_B(t_1) -s) > s -(V_B(t_1) - V_B) > s - 2mr^2 > 0$$
Set  
$$\widehat {R}((x,t),r) = B(x,r) \times ((2t_1 -t) - r^2, (2t_1 - t)+ r^2)$$
We use the above inequality with (\ref{dV_plus_poincare_LEQ_const}) to achieve following the integration above:
$$\mu( \{ x \in \widehat {R}: V_B- v(x,t) \} ) \leq {c \over s}.$$
This converges by the argument above. This gives us 
$$\int_{B(x_0, t_0 - 2r^2) \times (-r, r)} \sqrt{(V_B - v(x,t)  )^{+}} d \mu (x,t) \leq c_5$$
where the center of the new ball is $t_0 - 2r^2$ as expected.

This shows $v = -log(u)$ belongs to BMO with $h(s) = \sqrt{s^+}$ and lag mapping $T$. Aimar's main theorem then gives the desired result

$$(\int_{R((x,t);r)} u^{- \epsilon} d \mu) \cdot (\int_{R((x,\tau);r)} u^{ \epsilon} d \mu) \leq C \mu (R)^2$$
where $\tau = t - 2(r/ \eta)^2$. In the above 
$$\eta = [K_0(1+2[K^5_0(2(K_2^{-1} + 2 K_1)K_3 K_2^{-1} + 5 K_0^2 K_2^{-1} + K_1 ) + 1] K_0)]^{-1}$$
This bridges the gap if we choose from (\ref{p_required_for_bridge}) $p_{\nu} = \epsilon$.

\begin{remark}\label{steklov_proof}
We begin our discussion of an alternate proof with recalling equation (\ref{start_of_logu}) for a weak solution $u \in S^{1,2}(R),$ $R  \subset \mathbb{R}^{n+1}$.

\begin{equation}
\iint_{R} \phi u_t \ dxdt +  \sum^{m}_{i,j=1}\iint_{R}(a_{ij}(x,t) X_j u) X_i\phi \ dxdt \ \geq 0.
\end{equation}
From (\ref{start_of_logu}) we proved the inequality 
\begin{equation}\label{int_v-V_B_leq_Const}
\int_{\widetilde {R}((x,t),r)} \sqrt{(v(x,t) - V_B)^{+}} \; d \mu (x,t) \leq \, c_4
\end{equation} 
Now we will alter the argument from the beginning with the Steklov average. 
\begin{defn}{Steklov average}

\noindent For a positive, weak solution $u\in S^{1,2}(R)$ to (\ref{heat}) we define 
$$u_h \ = \ \frac{1}{h} \int_{t}^{t+h} u(x, \tau) d \tau$$
\end{defn}
\noindent Observe that if $u$ is a weak solution to (\ref{heat}) then so is $u_h$. Also, we have that $u_h$ is differentiable in $t$ so we can use the following form of the weak solution: For each $t \in [t_1 , t_2] \subset [0, T]$ we have 
\begin{equation}\label{logu_with_steklov}
\iint_R \phi \ \partial_t u_h(x,t) \ dx dt + \sum^{m}_{i,j=1} \iint_R (a_{ij}(x,t) X_j u_h(x,t)) X_i \phi \ dxdt = 0 
\end{equation}
for every $\phi(x,t) \in C^{\infty}(R)$ for each fixed $t$.

We also have the following results from $L^p$ theory concerning the Steklov average \cite{L}:
\begin{enumerate} 
\item If $u\in L^p(\Omega \times (0,T))$ then $u_h \to u$ in $L^p(\Omega \times (0,T-\delta))$ as $h \to 0^+$ for $0 < \delta < T$.
\item If $\nabla u \in L^p(\Omega \times (0,T))$ then $\nabla u_h \to \nabla u$ in $L^p(\Omega \times (0,T-\delta))$ as $h \to 0^+$ for $0 < \delta < T$.
\end{enumerate}
Here $\nabla u$ denotes the Euclidean spatial gradient of $u$.  Notice that the analogous result for the horizontal gradient of $u$ follows from (ii).

We may continue from (\ref{logu_with_steklov}) to achieve the same bounds on $log(u)$. Then, at the end of the proof of the bounds on $log(u)$ before the application of the modified John-Nirenberg theorem of Aimar we will have from (\ref{bridge_inequality})
\begin{equation}
{1 \over |R((x,t),r)|} \iint_{R((x,t),r)} u_h^{- \epsilon} dxdt \cdot {1 \over |S(R)|} \iint_{S(R)} u_h^{ \epsilon} dxdt \leq C
\end{equation}
Now we may apply the results from the $L^p$ theory of Steklov averages:
$$
\underset{h \rightarrow 0^+}{\lim} {1 \over |R((x,t),r)|} \iint_{R((x,t),r)} u_h^{- \epsilon} dxdt \cdot {1 \over |S(R)|} \iint_{S(R)} u_h^{ \epsilon} dxdt 
$$
$$= {1 \over |R((x,t),r)|} \iint_{R((x,t),r)} u^{- \epsilon} dxdt \cdot {1 \over |S(R)|} \iint_{S(R)} u^{ \epsilon} dxdt \ \leq C.
$$

%\noindent Notice that $R$ is bounded in $\mathbb{R}^{n+1}$ and we have 
%$$$$
\end{remark}
 %-------------------------------------------------------------
% You cannot use the \chapter{} or \section{} commands because they put a header that the grad school doesn't allow
\begin{center} 
\textbf{1.8  H\"older Continuity}
\end{center}
%-------------------------------------------------------------

We now show that weak solutions in a bounded domain are H\"older continuous. This section follows closely the approach by Moser \cite{Mo1}. Let us construct appropriate domains to examine the continuity. First let $u$ be a positive, weak solution to $$\parab$$ in
$$R = B(x,r) \times (0, \tau)$$ and subdomains
$$R^+ = B(x, r/4) \times (7 \tau /8, \tau )$$
$$R^- = B(x, r/4) \times ( \tau /8, \tau / 6)$$
We will use the distance $d((x,t),(0,0)) = \max \, \{ d_{cc}(x,0), \sqrt{|t|} \}$ for $x \in \mathbb{R}^n$ and $t \in \mathbb{R}$.  
Also, we will form an $n+1$ dimensional vector $x' = (x_1, x_2, ... t) \in R \subset \mathbb{R}^{n+1}$. With this extra notation we will now restate the distance $d((x,t),(y,s))$ as 
$$d(x', y') := d((x,t),(y,s)).$$
Then we have for a positive, weak solution in $R$, $u(x,t):=u(x')$ a corollary to the main theorem: 
 
\begin{cor}(H\"older Continuity of weak solutions) Let $u \in S_{Loc}^{1,2}(R)$  be a weak positive solution to $$\parab$$ in a domain $R$. Then using the notation from above, there exist constants $C >0$, and $0 < \alpha < 1$ such that for every $x'$, $y' \in R' \subset R$, $x' \neq y'$, we have 
\begin{equation}
|u(x') - u(y')| \leq C \,d(x',y')^{\alpha}
\end{equation}
\end{cor}
\noindent Set the average of $u$ in $R^-$: 
$$\mu^- := {1 \over |R^-|} \iint_{R^-} u(x,t)\, dxdt \, \leq \underset{R^-}{\max} \, u $$ 
We proved earlier in the main theorem that 
\begin{equation}
\label{min_GEQ_mu}\underset{R^+}{\min} \, u \, \geq \gamma ^{-1} \mu ^{-1} = \gamma ^{-1} {1 \over |R^-|} \iint_{R^-} u(x,t) \, dxdt
\end{equation}
For $p=1.$
For any positive, weak solution $u$ in a bounded domain $R$, let
$$M = \underset{R}{ess \max} \, u$$
$$m = \underset{R}{ess \min} \, u$$
$$M^+ = \underset{R^+}{ess \max} \, u$$
$$m^+ = \underset{R^+}{ess \min} \, u$$
Then $M-u$ is a non-negative weak solution to $$\parab$$ since it is $const - u \geq 0$. Likewise with $u -m$. The mean values of $M - u$, $u - m$ in $R^-$ are
$${1 \over |R^-|} \iint_{R^-} (M-u) \, dxdt = M - \mu^-$$
$${1 \over |R^-|} \iint_{R^-} (u-m)\,  dxdt = \mu^- - u$$
Apply (\ref{min_GEQ_mu}) to both 
\begin{equation}
\label{min_of_M_minus_u}\underset{R^+}{\min}\, (M-u) \, = (M-M^+) \geq \gamma^{-1} (M - \mu^-).
\end{equation} Since $(M-u)$ is minimized when $u$ is maximized in $R^+$.
\begin{equation}
\label{min_of_u_minus_m}\underset{R^+}{\min}\, (u-m) \, = (m^+ - m) \geq \gamma^{-1} (\mu^- -m).
\end{equation} Since $(u-m)$ is minimized when $u$ is minimized.
Let the oscillations $\omega = M - m$ and $\omega^+ = M^+ - m^+$ represent the difference between the max and min on a given domain. Adding (\ref{min_of_M_minus_u}) and (\ref{min_of_u_minus_m}) we find that
\begin{equation}
(M - M^+) + (m^+ -m) \geq \gamma ^{-1} (\mu ^- -m + M - \mu^-) 
\end{equation}
implies
$$(M - m) - (M^+ -m^+) \geq \gamma ^{-1} (M -m) $$
$$\omega - \omega^+ \geq \gamma^{-1} \omega $$
$$\omega - \gamma^{-1} \omega \geq \omega^+ $$
$$\omega (1 - \gamma^{-1}) \geq \omega^+ $$
So we can say
\begin{equation}\label{omega_plus_leq_theta_omega}\omega ^+ \leq ({\gamma - 1 \over \gamma}) \omega = \theta \omega
\end{equation} for some $0 < \theta < 1.$ Let $x' \in  R' \subset R$ with a $k$-neighborhood also contained in $R'$.
So for any $x'$, $y' \in R'$ such that $d(x',y')< k$ then the oscillation of $u$ in $R'$ is $\omega.$ Also we have if $d(x',y')< {k \over 4}$ then the oscillation of $u$ in $R'$ is $\theta \omega.$ Finally this allows us to say if $d(x',y')< {k \over 4^{\nu}}$ then the oscillation of $u$ in $R'$ is $\theta^{\nu} \omega.$

\noindent Let $R' \subset R$ be any region. Let $x' \in R'$ with a $k$-neighborhood contained wholly in $R$. Then for $x'$, $y' \in R'$ satisfying $d(x' , y') < k$, choose $\nu$ such that ${k \over 4^{\nu +1}} \leq d(x' , y') < {k \over 4^{\nu}}$. Then we have ${1 \over 4^{\nu}} \leq {4 \over k}d(x', y')$. Putting these together we have

\begin{equation}\label{holder_cty}
|u(x') - u(y')| \leq \theta^{\nu} \omega \leq ({1 \over 4})^{\alpha \nu} \omega = [({1 \over 4})^{\nu}]^{\alpha} \omega \leq [{4 \over k} d(x' , y')]^{\alpha} \omega
\end{equation}
$$= [({4 \over k})^{\alpha} \omega] d(x', y')^{\alpha} = Cd(x',y')^{\alpha}$$
with $$\alpha = {-log(\theta) \over log(4)} = {-log({\gamma - 1 \over \gamma}) \over log(4)}$$ and $0<\alpha < 1$. This proves the H\"older continuity with respect to the CC metric and so our second main result. We notice that with the above result and a proposition of Nagel, Stein and Wainger \cite{NaW} we have that the solutions to (\ref{heat}) are H\"older continuous with respect to the Euclidean metric. They proved that there exists a constant $c > 0$ and a constant $\epsilon$ depending on the number of commutators it takes to span the space such that
$$c |x-y| \leq d_{cc}(x,y) \leq c^{-1} |x-y|^{\epsilon}.$$
So with (\ref{holder_cty}) we have
$$|u(x,t) - u(y,s)| \leq Cd(x',y')^{\alpha} \leq c^{-1} C|x-y|^{\epsilon \alpha}. $$

\pagebreak
%%%%%%%%%%%%%%%%%%%%%%%%%%%%%%%%%%%%%%%%%%%%%%%%%%%%%%%%%%%%%%%%%%%%%%%%%%%%%%%%%%%%%%%%%%%%%%%%%%%%%%%%%%%
% You cannot use the \chapter{} or \section{} commands because they put a header that the grad school doesn't allow
\begin{center}
CHAPTER 2\\ \vspace{7pt}
\textbf{Nonlinear Generalization}\\ \vspace{7pt}
\textbf{2.1. Introduction}
\end{center}
\vspace{7pt}
\setcounter{chapter}{2} %This makes sure all the equations and theorems are labeled 2.* instead of 1.*
\setcounter{equation}{0} %This makes sure the equations/eqnarrays start with 2.1 instead of continuing the previous numbers
\setcounter{thm}{0} %This makes sure the theorems/definitions etc. start with 2.1 instead of continuing the previous numbers
%%%%%%%%%%%%%%%%%%%%%%%%%%%%%%%%%%%%%%%%%%%%%%%%%%%%%%%%%%%%%%%%%%%%%%%%%%%%%%%%%%%%%%%%%%%%%%%%%%%%%%%%%%%%%
 Now we will consider a more general non-linear case of which the above results are a special case. We will examine the operator
\begin{equation}
\label{nonlinear_operator}Lu = u_t + \sum^{m}_{j=1} X^{*}_j A(x, \hg u)
\end{equation}
And corresponding solutions to 
\begin{equation}
\label{non_linear_equation} Lu = 0
\end{equation}
Where $A:\mathbb{R}^{m} \rightarrow \mathbb{R}^m$ is a measurable function, and $X^{*}_j$ is the formal adjoint of $X_j$.
Let $u \in S^{1,2}_{loc}(\mathbb{R}^{n+1})$. We want to consider weak solutions $u$ to (\ref{non_linear_equation}) in a rectangle $R((x,t),r) = B_{cc}(x,r) \times (0,r)$. Here, $r \in \mathbb{R}$, $x \in \mathbb{R}^{n}$. By $B_{cc}(x,r)$ we mean a Carnot-Carath\'eodory ball $B(x,r) := \{ y \in \mathbb{R}^n : d_{cc}(x,y) < r  \}$ where $d_{cc}$ is the natural control metric associated with the vector fields $\{ X_i \}^{m}_{i=1}$. In the following we will refer to $d_{cc}(x,y)$ as $d(x,y)$. By a weak solution to $Lu = 0$ in $R((x,t),r)$ we mean
\begin{equation}
\label{definition_of_weak_solution} \int_{R((x,t),r)} \phi u_t \ dxdt +  \sum ^{m}_{j=1} \int_{R((x,t),r)} A_j (x, \hg u) X_j \phi \,  dx dt = 0.
\end{equation} 
for all $\phi \in C^{\infty}_{0}(B(x,r))$ with $\phi (x) \geq 0$.

We will consider the following structural conditions on the function $A$ which we will denote (S).
%\begin{equation*}

$$\noindent (S) \quad   \quad \quad  \quad \left\{
 \begin{array}{rl}
  A_j(x, \hg u) \leq C|\hg u|\\
  \lambda |\hg u|^2 \leq A_j(x, \hg u) X_j u \leq \Lambda |\hg u|^2\\
 \end{array} \right.$$
%\end{equation*}
for some constants $C$, $\Lambda$, $\lambda > 0$.

Notice that in the parabolic case we considered for the Moser iteration, 
$$A_j(x, X_j u) = a_{ij}(x,t) X_j u.$$

We will now develop minimum and maximum inequalities for weak solutions to (\ref{non_linear_equation}) in a rectangle. 
\begin{thm}
Let $R$ be a rectangle, $R((x,t),r) = B_{cc}(x,r) \times (0,r)$ in $\mathbb{R}^{n+1}$. Let
$$ Lu = u_t + \sum^{m}_{j=1} X^{*}_j A(x, \hg u)
$$
and consider positive, weak solutions $u \in S^{1,2}(R)$ to $Lu = 0$ in $R$. Here $\{ X_i \}^{m}_{i=1}$, $m \leq n$, are vector fields satisfying H\"ormander's condition for hypoellipticity. Under the assumption of the structural conditions $(S)$ we have the following
\begin{equation}
\underset{R'}{\max} \ u \ \leq \gamma_1({1 \over |R|} \int_R u^p \  dxdt)^{1/p} 
\end{equation} 
\begin{equation}
 ({1 \ \over |R|} \int_R u^p \ dxdt)^{1/p} \  \leq \  \gamma_2 \  \underset{R'}{\min} \ u
\end{equation}
for some constants $\gamma_1$,  $\gamma_2 > 0$ depending on $R$, $n$, $\Lambda$, $\lambda$, and $C$.
\end{thm}

With the proof of the preceding two inequalities, we will move to the corresponding bridge step 
\begin{equation}
 ({1 \ \over |R|} \int_R u^{p_0} \ dxdt )^{1/p_0} \  \leq \  \gamma_3 ({1 \ \over |R|} \int_R u^{-p_0} \ dxdt)^{-1/p_0}
 \end{equation}
as in the parabolic case above. Together, these three inequalities will give us a Harnack inequality for weak solutions to (\ref{non_linear_equation}) in $R$. 
 %-------------------------------------------------------------
% You cannot use the \chapter{} or \section{} commands because they put a header that the grad school doesn't allow
\begin{center} 
\textbf{2.2 Non linear Cacciopoli Inequality}
\end{center}
%-------------------------------------------------------------
 
We will begin the proof of a Harnack inequality in the same manner as above with the Cacciopoli inequality. For convenience, we will let $R = R((x,t),r)$.
 
\begin{thm}{Local bounds for weak solutions to $Lu = 0$}\label{nonlinearlocalbounds}

With the hypotheses above, including the structural conditions $(S)$, let $u \in S^{1,2}(R)$ be a positive, weak solution to 
\begin{equation}
 Lu = 0
\end{equation}
 in $R = \{B(x, r ) \times (0, \tau) \}$. Also define $R' = \{B(x, r' ) \times (0,\tau') \}$ for $r  > r' > 0$, $\tau > \tau ' > 0$. Then there exists constants $C_1$, $C_2$ depending on $\lambda$, $n$, and $R$ such that 
\begin{equation}
\label{nonlinear_cacc_gradv} \int_{R'} | \hg u |^2 dxdt \ \leq C_1 \int_{R} u^2 dxdt
\end{equation}
and 
\begin{equation}
\label{nonlinearcaccmaxt} \underset{t \in [0 ,    \sigma ]}{\max} \int_{B(x, \rho')} u^2 dxdt \ \leq C_2 \int_{R} u^2 dxdt
\end{equation}

Where the max is the essential supremum over this interval.
\end{thm}
\emph{Proof:}

\noindent Consider a positive, weak solution $u$ to $Lu =0$: 
\begin{equation}
\label{nonlineardefinition_of_weak_solution}\int_{R } \phi \ u_t \ dxdt+  \sum ^{m}_{j=1} \int_{R } A_j (x, \hg u) X_j \phi \,  dx dt = 0
\end{equation}  
Now we let the test function $\phi = u \psi ^2$ where $\psi$ is a function compact in the $x$ variable and independent of $t$. Also, we construct the test function to give $\psi = 1$ on $R'$. Substituting gives:
$$ \int_{R } u \psi^2 \ u_t \ dxdt+  \sum ^{m}_{j=1} \int_{R } A_j (x, \hg u) X_j(u \psi^2) \,  dx dt = 0 $$
Notice that ${(u^2)_t \over 2} = u u_t$ and $X_j(u \psi^2) = X_ju \psi^2 + 2u \psi X_j \psi$. This gives us:
\begin{equation}
\int_{R } {(u^2)_t \over 2} \psi^2 \ dxdt +  \sum ^{m}_{j=1} \int_{R } A_j (x, \hg u) X_ju \psi^2  \ dxdt 
\end{equation}  
$$+ 2\sum ^{m}_{j=1} \int_{R } A_j (x, \hg u) u \psi X_j \psi \,  dx dt = 0 $$
\begin{equation}
\label{nonlinearcacc_before_schwarz}\int_{R } {(u^2)_t \over 2} \psi^2  \ dxdt  +  \sum ^{m}_{j=1} \int_{R } A_j (x, \hg u) X_ju \psi^2 \ dxdt
\end{equation}  
$$= -  2\sum ^{m}_{j=1} \int_{R } A_j (x, \hg u)u \psi X_j \psi \,  dx dt
$$
Now we apply the Schwarz and Young's inequalities to the RHS of (\ref{nonlinearcacc_before_schwarz}).
$$\int_{R } {(u^2)_t \over 2} \psi^2 \ dxdt +  \sum ^{m}_{j=1} \int_{R } A_j (x, \hg u) X_ju \psi^2 \ dxdt$$
$$  
 \leq   {1 \over 2}  \int_{R } |A_j (x, \hg u)|^2 \psi^2  \ dxdt+ 2 \int_{R } u^2 | \hg \psi|^2  \,  dx dt
$$
and recall the bounds on $A_j$ from $(S)$: 
$$\int_{R } {(u^2)_t \over 2} \psi^2   \ dxdt+   \lambda \int_{R } |\hg u|^2 \psi^2 \ dxdt$$
$$  
 \leq    {C^2 \over 2} \int_{R } |\hg u|^2 \psi^2 \ dxdt +  2\int_{R }|\hg \psi|^2 u^2  \,  dx dt$$
And adding like terms we arrive at:
\begin{equation}
\label{reference_for_nonlinear_powers}\int_{R } {(u^2)_t \over 2} \psi^2  \ dxdt +   (\lambda -{C^2 \over 2} )\int_{R } |\hg u|^2 \psi^2 \ dxdt \leq  2\int_{R }|\hg \psi|^2 u^2  \,  dx dt
\end{equation}
Notice that $\frac{1}{2}(\psi ^2 u^2)_t = \psi \psi _t u^2 + u u_t \psi ^2$. And use this fact by adding $\psi \psi_t u^2 $ to both sides: 
\begin{equation}
\label{nonlinearcacc_before_bounds_on_psi}\int_{R } \frac{1}{2}(\psi ^2 u^2)_t  \ dxdt +   (\lambda -{C^2 \over 2} )\int_{R } |\hg u|^2 \psi^2 \ dxdt \leq  2\int_{R } u^2( |\hg \psi|^2 + \psi \psi_t)\,  dx dt
\end{equation}
 
\noindent Now integrate this over $R_{\sigma} = \{B(x, r) \times 0 \leq t \leq \sigma < \tau  \}$ and drop the second integrand on the LHS. It is $\geq 0$. We'll choose $\sigma$ so that 

\begin{equation}
\label{nonlinearmaxt}\int_{B(x, p')} u^2 |_{t =   \sigma} dx \ \geq \ \frac{1}{2} \underset{t \in (0 ,   \sigma)}{\max} \int_{B(x, r ')} u^2 dx.
\end{equation}
 \noindent This gives us from (\ref{nonlinearcacc_before_bounds_on_psi}):
\begin{equation}
\label{nonlinearmaxtbounds}\int _{R_{\sigma}} \frac{1}{2} (u^2 \psi ^2)_t \, dx dt \leq 2 \int_{R_{\sigma}} u^2 (|\hg \psi|^2 \ + \ |\psi \, \psi _t|) \, dxdt.
\end{equation}

\noindent Now consider the 
\begin{equation}
\underset{t \in (0 ,  \sigma)}{\max} \int_{B(x, r ')} u^2 \, dx \ \leq \ 2 \int_{B(x, r ')} u^2 \  |_{t =  \sigma} \ dx
\end{equation} 
which we constructed above in (\ref{nonlinearmaxt}). $\psi = 1$ on $R'$ and $t$ is fixed, so 

$$\underset{t \in (0,   \sigma)}{\max} \int_{B(x, r ')} u^2 \, dx \, \leq \, 2 \int_{B(x, r ')} u^2 \,  |_{t =  \sigma} \ dx $$

$$\leq \, 2 \int_{B(x, r ')} u^2 \psi ^2 |_{t =   \sigma} \, dx$$

$$\leq \ 2 \int_{B(x, r)} u^2 \psi ^2 |_{t =   \sigma} \, dx$$

\noindent Now recall (\ref{nonlinearmaxtbounds}) integrated over $R_{\sigma}$. We will combine that with the $\frac{1}{2  }$ factor as well as the fact that the above is twice the (\ref{nonlinearmaxtbounds}) inequality. Putting these together gives:
\begin{equation}
\underset{t \in (0 ,   \sigma)}{\max} \int_{B(x, r ')} u^2 dx \ \leq \ 8 \lambda \iint_{R_{\sigma}} u^2 (|\hg \psi|^2 + |\psi \psi _t|) \, dx dt.
\end{equation}

\noindent Then use the non-negativity of the RHS and take the integral over $R$ for
\begin{equation}
\leq 8   \iint_{R} u^2 (|\hg \psi|^2 + |\psi \psi _t|) \, dx dt.
\end{equation}

\noindent Now we use as result on the bounds of the horizontal gradient of test functions in Citti, Garofalo, and Lanconelli \cite{CiL} in order to attain a constant independent of $ \psi $.

\begin{lem}{ Existence of Carnot-Carath\'eodory test functions.} There exists $r _0 > 0$ such that given a metric ball $B(x, r) \subset \subset R$, with $r \leq r _0$ and $0 < r ' < r$, there exists a function $\psi \in C^{\infty}_0 (B(x, r))$ such that 
$0 \leq \psi \leq 1$, $\psi \equiv 1$ in $B(x, r ')$ and $|\hg \psi| \leq \frac{C_{HG}}{r - r '}$. Here, $C_{HG} > 0$ is a constant independent of $r$ and $r '$.   
\end{lem}
\noindent Define $\psi (x,t) := \psi _1 (t) \psi _2 (x)$ in the following way, using the cut-off functions of Citti, Garofalo and Lanconelli \cite{CiL}:
$$
%\begin{equation*}
\psi _1 (t) = \left\{
 \begin{array}{rl}
  1 & \text{if }  \tau' < t < \tau \\
   0 & \text{if } t \leq 0 \\
 \end{array} \right.
%\end{equation*}
$$

$$
%\begin{equation*}
\psi _2 (x) = \left\{
 \begin{array}{rl}
  1 & \text{if } d_{cc}(x, 0) < r '\\
   0 & \text{if } d_{cc}(x, 0) > r \\
 \end{array} \right.
%\end{equation*}
$$

\noindent and as before we have a bound on the horizontal gradient. This allows us to say
\begin{equation}
|\hg \psi |^2 \ + \ |\psi \psi _t | \ \leq \ \bigl( \frac{C_{HG}^2}{(r - r ')^2} + \frac{1}{\tau - \tau '} \bigr)
\end{equation}

\noindent Which gives us
\begin{equation}
\underset{t \in (0 ,   \sigma)}{\max} \int_{R '} u^2 \, dx \ \leq \ 8   (\frac{C_{HG}^2}{(r - r ')^2} + \frac{1}{\tau - \tau '}) \iint_{R} v^2 \, dxdt.
\end{equation}
This proves (\ref{nonlinearcaccmaxt}). We note that the result is similar to the parabolic case since the approach was the same. The difference is that the bounds on the $a_{ij}(x,t)$ matrix are gone and we find only the bounds from the structural conditions $(S)$. This is also true the following proof of the bounds on the derivatives.

\noindent Now we return to (\ref{nonlinearcacc_before_bounds_on_psi})
$$
\int_{R } \frac{1}{2}(\psi ^2 u^2)_t   \ dxdt+   (\lambda -{C^2 \over 2} )\int_{R } |\hg u|^2 \psi^2 \ dxdt \leq  2\int_{R } u^2( |\hg \psi|^2 + \psi \psi_t)\,  dx dt
$$
and drop the first integrand. Integrated over $R'$ gives
\begin{equation}
\label{nonlinearbounds_on_hgu_before_any_calc}(\lambda -{C^2 \over 2} )\int_{R' } |\hg u|^2 \psi^2 \ dxdt \leq  2\int_{R' } u^2( |\hg \psi|^2 + \psi \psi_t)\,  dx dt
\end{equation}
Again use the fact that $\psi = 1$ on $R'$
$$
 (\lambda -{C^2 \over 2} )\int_{R'}  |\hg u|^2 \, dxdt \leq 2   \int_{R'} u^2 \, (|\hg \psi |^2 + |\psi \psi _t|) \, dxdt
 $$

\noindent Which gives us
$$
\int_{R'}  |\hg u|^2 \, dxdt \leq 2 (\lambda -{C^2 \over 2} )^{-1}  \int_{R'} u^2 \, (|\hg \psi |^2 + |\psi \psi _t|) \, dxdt.
$$
And using the previous bounds on $\psi$, $\psi_t$ and $|\hg \psi|$
\begin{equation}
\int_{R'}  |\hg u|^2 \, dxdt \leq 2 (\lambda -{C^2 \over 2} )^{-1}  (\frac{C_{HG}^2}{(r - r ')^2} + \frac{1}{\tau - \tau '}) \int_{R} u^2 \ dxdt.
\end{equation}

\noindent And the integral on the RHS is taken over $R$ since it is non-negative. This proves the result in (\ref{nonlinear_cacc_gradv}), and so Theorem \ref{nonlinearlocalbounds} is proved.
 %-------------------------------------------------------------
% You cannot use the \chapter{} or \section{} commands because they put a header that the grad school doesn't allow
\begin{center} 
\textbf{2.3 Nonlinear Cacciopoli for $u^p$ }
\end{center}
%-------------------------------------------------------------

Since we have a Cacciopoli inequality for weak solutions $u$, now we look at powers of $u$, $v:=u^p$.
\begin{cor}
The same inequality of Theorem \ref{nonlinearlocalbounds} holds, up to a constant, for powers of a weak solution $u$, $v=u^p$. 
\end{cor}
\noindent \emph{Proof:}
Let $\phi = \psi^2 u^p$ for $p > {1 \over 2}$ and $\psi$ be defined as above. Calculate $$v_t = (u^p)_t = pu^{p-1}u_t$$ and $$X_j(\phi) = X_j (\psi^2 u^p) = 2 \psi u^p X_j \psi + p \psi^2 u^{p-1} X_j u.$$ We substitute the definitions of $\phi$ and $v$ into (\ref{nonlineardefinition_of_weak_solution}):
$$
\int_{R} \psi^2 u^p (p)u^{p-1}u_t \ dxdt+  \sum ^{m}_{j=1} \int_{R} A_j (x, \hg u) (2 \psi u^p X_j \psi + p \psi^2 u^{p-1} X_j u)  \,  dx dt 
$$
which gives
$$
\int_{R} \psi^2 u^p (p) u^{p-1}u_t  \ dxdt+ \sum ^{m}_{j=1} \int_{R} A_j (x, \hg u) (p \psi^2 u^{p-1} X_j u) \ dxdt \ $$
$$ \leq \ -2 \sum ^{m}_{j=1} \int_{R} A_j (x, \hg u) \psi u^p X_j \psi \ dxdt.
$$
We use the Cauchy and Young's inequality on the RHS:
$$
\int_{R} \psi^2 u^p (p) u^{p-1}u_t  \ dxdt+ \sum ^{m}_{j=1} \int_{R} A_j (x, \hg u) (p \psi^2 u^{p-1} X_j u) \ dxdt \ $$
$$
\leq \ {1 \over 2}  \int_{R} |A_j (x, \hg u)|^2  \psi^2 \  \ dxdt + 2 \int_{R} u^{2p} |\hg \psi|^2 \  dxdt
$$
Use the bounds on $|A_j|$ from $(S)$ and substitute for $v$:
$$
\int_{R} {1 \over 2} \psi^2 (v^2)_t \ dxdt + \lambda \int_{R} |\hg v|^2  \psi^2  \ dxdt \ $$
%%%%%%%%%%%%%%%%%%%%%%%%%%%%%%%%%%%%%%%%%%%%%%%%%%%%%%%%%%
%%%%%%%%%%if you insert the calcs insert them here %%%%%%
%%%%%%%%%%And end them down there %%%%%%%%%%%%%%%%%%%%%%%
$$
\int_{R} {1 \over 2} \psi^2 (v^2)_t  \ dxdt+ \lambda \int_{R} |\hg v|^2  \psi^2  \ dxdt \ $$
$$
\leq \ {C^2 \over 2}  \int_{R} |\hg v|^2   \psi^2 \  dxdt + 2 \int_{R} v^2 |\hg \psi|^2 \  dxdt
$$
Adding like terms to get
$$
\int_{R} {1 \over 2} \psi^2 (v^2)_t  \ dxdt+ (\lambda  - {C^2 \over 2}) \int_{R} |\hg v|^2  \psi^2  \ dxdt \ $$
$$
\leq \ 2 \int_{R} v^2 |\hg \psi|^2 \  dxdt
$$
%From which point the proof follows as above in the case $v=u$, (\ref{reference_for_nonlinear_powers}) and following.
Notice that $\frac{1}{2}(\psi ^2 u^2)_t = \psi \psi _t u^2 + u u_t \psi ^2$. And use this fact by adding $\psi \psi_t u^2 $ to both sides: 
\begin{equation}
\label{nonlinearcacc_before_bounds_on_psiPOWERSOFU}\int_{R } \frac{1}{2}(\psi ^2 u^2)_t   \ dxdt+   (\lambda -{C^2 \over 2} )\int_{R } |\hg u|^2 \psi^2 \ dxdt \leq  2\int_{R } u^2( |\hg \psi|^2 + \psi \psi_t)\,  dx dt
\end{equation}
 
\noindent Now integrate this over $R_{\sigma} = \{B(x, r) \times 0 \leq t \leq \sigma < \tau  \}$ and drop the second integrand on the LHS. It is $\geq 0$. We'll choose $\sigma$ so that 

\begin{equation}
\label{nonlinearmaxtPOWERSOFU}\int_{B(x, p')} u^2 |_{t =   \sigma} dx \ \geq \ \frac{1}{2} \underset{t \in (0,   \sigma)}{\max} \int_{B(x, r ')} u^2 dx.
\end{equation}
 \noindent This gives us from (\ref{nonlinearcacc_before_bounds_on_psiPOWERSOFU}):
\begin{equation}
\label{nonlinearmaxtboundsPOWERSOFU}\int _{R_{\sigma}} \frac{1}{2} (u^2 \psi ^2)_t \, dx dt \leq 2 \int_{R_{\sigma}} u^2 (|\hg \psi|^2 \ + \ |\psi \, \psi _t|) \, dxdt.
\end{equation}

\noindent Now consider the 
\begin{equation}
\underset{t \in (0 ,   \sigma)}{\max} \int_{B(x, r ')} u^2 \, dx \ \leq \ 2 \int_{B(x, r ')} u^2 \  |_{t = \sigma} \ dx
\end{equation} 
which we constructed above in (\ref{nonlinearmaxtPOWERSOFU}). $\psi = 1$ on $R'$ and $t$ is fixed, so 

$$\underset{t \in (0 ,  \sigma)}{\max} \int_{B(x, r ')} u^2 \, dx \, \leq \, 2 \int_{B(x, r ')} u^2 \,  |_{t = \sigma} \ dx $$

$$\leq \, 2 \int_{B(x, r ')} u^2 \psi ^2 |_{t = \sigma} \, dx$$

$$\leq \ 2 \int_{B(x, r)} u^2 \psi ^2 |_{t = \sigma} \, dx$$

\noindent Now recall (\ref{nonlinearmaxtboundsPOWERSOFU}) integrated over $R_{\sigma}$. We will combine that with the $\frac{1}{2  }$ factor as well as the fact that the above is twice the (\ref{nonlinearmaxtboundsPOWERSOFU}) inequality. Putting these together gives:
\begin{equation}
\underset{t \in (0 ,  \sigma)}{\max} \int_{B(x, r ')} u^2 dx \ \leq \ 8 \lambda \iint_{R_{\sigma}} u^2 (|\hg \psi|^2 + |\psi \psi _t|) \, dx dt.
\end{equation}

\noindent Then use the non-negativity of the RHS and take the integral over $R$ for
\begin{equation}
\leq 8   \iint_{R} u^2 (|\hg \psi|^2 + |\psi \psi _t|) \, dx dt.
\end{equation}

\noindent Now we use a result on the bounds of the horizontal gradient of test functions in order to attain a constant independent of $ \psi $. Define $\psi (x,t) := \psi _1 (t) \psi _2 (x)$ in the following way, using the cut-off functions of \cite{CiL}:

$$
%\begin{equation*}
\psi _1 (t) = \left\{
 \begin{array}{rl}
  1 & \text{if }  \tau' < t < \tau\\
   0 & \text{if } t \leq 0 \\
 \end{array} \right.
%\end{equation*}
$$

$$
%\begin{equation*}
\psi _2 (x) = \left\{
 \begin{array}{rl}
  1 & \text{if } d_{cc}(x, 0) < r '\\
   0 & \text{if } d_{cc}(x, 0) > r \\
 \end{array} \right.
%\end{equation*}
$$

\noindent and as before we have a bound on the horizontal gradient. This allows us to say
\begin{equation}
|\hg \psi |^2 \ + \ |\psi \psi _t | \ \leq \ (\frac{C_{HG}^2}{(r - r ')^2} + \frac{1}{\tau - \tau '})
\end{equation}

\noindent Which gives us
\begin{equation}
\underset{t \in (0 ,  \sigma)}{\max} \int_{R '} u^2 \, dx \ \leq \ 8   (\frac{C_{HG}^2}{(r - r ')^2} + \frac{1}{\tau - \tau '}) \iint_{R} v^2 \, dxdt.
\end{equation}
This proves the result for (\ref{nonlinearcaccmaxt})
%%%%%%%%%%%%%%%%%%%%%%%%%%%%%%%%%%%%%%%%%
\noindent Now we return to (\ref{nonlinearcacc_before_bounds_on_psiPOWERSOFU})
$$
\int_{R } \frac{1}{2}(\psi ^2 u^2)_t  \ dxdt +   (\lambda -{C^2 \over 2} )\int_{R } |\hg u|^2 \psi^2 \ dxdt \leq  2\int_{R } u^2( |\hg \psi|^2 + \psi \psi_t)\,  dx dt
$$
and drop the first integrand. Integrated over $R'$ gives
\begin{equation}
\label{nonlinearbounds_on_hgu_before_any_calcPOWERSOFU}(\lambda -{C^2 \over 2} )\int_{R' } |\hg u|^2 \psi^2 \ dxdt \leq  2\int_{R' } u^2( |\hg \psi|^2 + \psi \psi_t)\,  dx dt
\end{equation}
Again use the fact that $\psi = 1$ on $R'$
$$
 (\lambda -{C^2 \over 2} )\int_{R'}  |\hg u|^2 \, dxdt \leq 2   \int_{R'} u^2 \, (|\hg \psi |^2 + |\psi \psi _t|) \, dxdt
 $$

\noindent Which gives us
$$
\int_{R'}  |\hg u|^2 \, dxdt \leq 2 (\lambda -{C^2 \over 2} )^{-1}  \int_{R'} u^2 \, (|\hg \psi |^2 + |\psi \psi _t|) \, dxdt.
$$
And using the previous bounds on $\psi$, $\psi_t$ and $|\hg \psi|$
\begin{equation}
\int_{R'}  |\hg u|^2 \, dxdt \leq 2 (\lambda -{C^2 \over 2} )^{-1}  (\frac{C_{HG}^2}{(r - r ')^2} + \frac{1}{\tau - \tau '}) \int_{R} u^2 \ dxdt.
\end{equation}

\noindent And the integral on the RHS is taken over $R$ since it is non-negative. This proves the result in (\ref{nonlinear_cacc_gradv}), and so the corollary is proved.
%%%%%%%%%%%%%%%%%%%% end insert here %%%%%%%%%%%%%%%%%%%
For the case $p < p_0 < 0$ we let $v= u^{p/2}$ be a weak supersolution to (\ref{nonlineardefinition_of_weak_solution}) in $R$. We then have for $\phi = v \psi ^2$ the following:
$$
v_t = (u^{p/2}) = \frac{p}{2} u^{p/2 -1} u_t,
$$
$$
X_j v = X_j (u^{p/2}) = \frac{p}{2} u^{p/2 -1} X_j u
$$
and 
$$
X_j \phi = X_j (v \psi ^2) = X_j (u^{p/2} \psi ^2) = \frac{p}{2} u^{p/2 -1} X_j u \psi ^2 + 2 \psi X_j \psi u^{p/2}
$$
And we substitute these into (\ref{nonlineardefinition_of_weak_solution}):
$$
\int_{R} (v \psi^2) \frac{p}{2} u^{p/2-1} u_t  \ dxdt+  \sum ^{m}_{j=1} \int_{R} A_j (x, \hg u) (\frac{p}{2} u^{p/2-1} X_j u \psi^2 + 2 \psi X_j \psi u^{p/2}) \  dx dt \ \geq 0. 
$$
which gives us
$$
\int_{R} (v \psi^2) \frac{p}{2} u^{p/2-1} u_t  \ dxdt+  \sum ^{m}_{j=1} \int_{R} A_j (x, \hg u) \frac{p}{2} u^{p/2-1} X_j u \psi^2 \ dxdt
$$
$$
\geq -2 \sum ^{m}_{j=1} \int_{R} A_j (x, \hg u)  \psi X_j \psi u^{p/2} \  dx dt.    
$$
At this point, notice that $p<0$ and the sign on the $2$ in the RHS. We multiply by $-1$ and proceed.
$$
\int_{R} (v \psi^2) \frac{|p|}{2} u^{p/2-1} u_t \ dxdt +  \sum ^{m}_{j=1} \int_{R} A_j (x, \hg u) \frac{|p|}{2} u^{p/2-1} X_j u \psi^2 \ dxdt
$$
$$
\leq 2 \sum ^{m}_{j=1} \int_{R} A_j (x, \hg u)  \psi X_j \psi u^{p/2} \  dx dt.    
$$
We calculate ${1 \over 2} (v_t)^2 = {p \over 2}u^{p/2}u^{p/2-1}u_t$ and substitute into the first integrand: 
$$
\int_{R} {1 \over 2} \psi^2 (v_t)^2  \ dxdt+  \sum ^{m}_{j=1} \int_{R} A_j (x, \hg u) \frac{|p|}{2} u^{p/2-1} X_j u \psi^2 \ dxdt
$$
$$
\leq 2 \sum ^{m}_{j=1} \int_{R} A_j (x, \hg u)  \psi X_j \psi u^{p/2} \  dx dt.    
$$
Again we use the Cauchy and Young's inequalities on the RHS and recognize that $u^p = v^2$ and $\hg v = \frac{|p|}{2} u^{p/2-1} X_j u$.
$$
\int_{R} {1 \over 2} \psi^2 (v_t)^2  \ dxdt+  \sum ^{m}_{j=1} \in_{R}t A_j (x, \hg u) \hg v \psi^2 \ dxdt
$$
$$
\leq {1 \over 2}  \int_{R} |A_j (x, \hg u)|^2  \ dxdt + 2 \int_{R} v^2 |\hg \psi|^2 \  dx dt.    
$$
Now use the $(S)$ bounds on $A_j$
$$
\int_{R} {1 \over 2} \psi^2 (v^2)_t \ dxdt + \lambda \int_{R} |\hg v|^2  \psi^2  \ dxdt \ 
\leq \ {C^2 \over 2}  \int_{R} |\hg v|^2   \psi^2 \  \ dxdt + 2 \int_{R} v^2 |\hg \psi|^2 \  dxdt
$$
And from this point, the proof proceeds exactly like above. The case for $0 < p < 1$ is modified in the following way: we must reverse $t$ and use $v:= u^p(x,-t)$. Then our test function must be $\phi = u^p \psi^2$. The proof proceeds in the same manner except there is a factor of $ \frac{\lambda (2p-1)}{p^2} $ in front of the $|\hg v|^2$ term.
 %-------------------------------------------------------------
% You cannot use the \chapter{} or \section{} commands because they put a header that the grad school doesn't allow
\begin{center} 
\textbf{2.4 Nonlinear Iteration and bridge step  }
\end{center}
%-------------------------------------------------------------
 
In order to run a Moser type iteration scheme, we needed the local bounds on
\begin{equation}
 \int_{R'} | \hg u |^2 dxdt \ \leq C_1 \int_{R} u^2 dxdt
\end{equation}
and 
\begin{equation}
\underset{t \in [0 ,  \sigma ]}{\max} \int_{B(x, r')} u^2 dxdt \ \leq C_2 \int_{R} u^2 dxdt
\end{equation}
Saloff-Coste shows in \cite{SC} that the parabolic Harnack inequality is equivalent to the doubling property $$\mu B(x, 2R) \leq A \mu B(x,R)$$ and the Poincar\'e inequality $$\int_{B(x,r)} |v -V_B|^2 \, dx \leq Cr^2 \int_{B(x,r)} |\nabla v|^2 dx$$ as we have used in the proof above. This is done in the Riemannian manifold setting. He assumes the matrix of coefficients, $a_{ij}$, are bounded and measurable, symmetric and elliptic. Saloff-Coste addresses the Moser iteration and on its dependence on both the Sobolev inequality and the Poincare inequality. Also, Grigory'an \cite{Gr} showed the same for non-compact Riemannian manifolds. Now we will use the local bounds above to run the Moser iteration. 

\begin{thm}\label{NLmaxu} 
Let $u$ be a positive, weak subsolution to $$\nlparab$$ in $R$ with $R'$ defined as above. If the structural hypotheses $(S)$ and the general assumptions from the local bounds theorem hold, then for $p>1$, there exists a constant $\gamma$ depending on $\lambda$, $C$, $n$, and $R$ such that
\begin{equation}
\label{nonlinearmaxuleqpnorm}\underset{R'}{\max} \ u \leq \gamma (\frac{1}{|R|} \iint_{R} u^p \, dxdt )^{1/p}. 
\end{equation}
\end{thm} 
%%%%%%%%%%%%%%%%%% insert calcs for max u %%%%%%%%%%%%%%%%%%%%%%%%%%%%%
%%%%%%%%%%%%%%%%%%%%%%%%%%%%%%%%%%%%%%%%%%%%%%%%%%%%%%%%%%%%%%%%%%%%%%%%
The proof follows just as in the parabolic case. The constant $\gamma$ is slightly different. This is due to the fact that in the original parabolic case, the $a_{ij}(x,t)$ are bounded and elliptic with constants $\lambda$ and ${1 \over \lambda}$ whereas the function $A_j(x, \hg u)$ in the nonlinear case is bounded by a constant $C$ and another $\lambda$ which is not an eigenvalue. Here $\gamma$ is the finite sum of terms
\begin{equation}
\gamma = [C_{B}(\frac{r}{r '})^Q \, 8  (\frac{C_{HG}^2}{(r - r ')^2} + \frac{1}{\tau - \tau '}) +1 ] \bigl(  8 (\frac{C_{HG}^2}{(r - r ')^2} + \frac{1}{\tau - \tau '}) \,  C_{B} \, (\frac{r}{r '})^Q   \bigr)^{\frac{2}{n}}.
\end{equation}
\emph{Proof:}
In order to prove the theorem, we use the lemma modeled after Moser \cite{Mo1}. Now for $v \in S^{1,2}(R)$, define as before:
\begin{equation}
H_{r,  \tau}(v):=  \normbrt \iint_R v^2 \, dxdt
\end{equation}
\begin{equation}
M_{r, \tau} (v) :=  \normbr \, \underset{t \in [0, \tau ]}{\max} \int_{B_{r}} v^2 dx
\end{equation}
\begin{equation}
D_{r, \tau}(v) := \normbrt \, \iint_R |\hg v|^2 dxdt
\end{equation}

\begin{lem}
Let $v \in S^{1,2}(R)$. Let $k=1 + \frac{2}{n}$ for $n \in \mathbb{N}$, $\alpha = \frac{n}{n-2}$, $\beta = \frac{n}{2}$ Then, we have for some constants $C_1$, $C_2 >0$: 
\begin{equation}
\label{preliminaryiteration}H_{r',  \tau}(v^k) \leq C_1 \, M_{r,  \tau}(v)^{\frac{2}{n}} \bigl[ C_2 H_{r',  \tau}(v) + D_{r',  \tau}(v) \bigr]
\end{equation}

\end{lem}
which we proved above. In terms of $H_{r, \tau}(v)$, $D_{r, \tau}(v)$ and $M_{r, \tau}(v)$ we calculate our local bounds (\ref{nonlinear_cacc_gradv}) and (\ref{nonlinearcaccmaxt}) as:
\begin{equation}
M_{r , \tau } (v) = \normbrp \underset{t \in (0, \tau )}{\max} \int_{B_{r }} v^2 dx 
\end{equation}
\begin{equation}
 \leq  8   (\frac{C_{HG}^2}{(r - r ')^2} + \frac{1}{\tau - \tau '}) \normbr \iint_{R} v^2 dx dt
\end{equation}
\begin{equation}
\leq  8   (\frac{C_{HG}^2}{(r - r ')^2} + \frac{1}{\tau - \tau '})  H_{r, \tau}(v)
\end{equation}
\begin{equation}
\leq  8  (\frac{C_{HG}^2}{(r - r ')^2} + \frac{1}{\tau - \tau '}) \,  H_{r ,\tau}(v)
\end{equation}

\noindent Then for some $C_3$ depending on $\lambda$, $R$ and $n$ we have 
\begin{equation}
M_{r , \tau } (v) \leq C_3 \, H_{r ,\tau}(v)
\end{equation}

\noindent And for 
$$D_{r , \tau '} (v) = {1 \over |B_{r} \times (0, \tau ')|}  \, \iint_{R} |\hg v|^2 \, dxdt $$
we calculate:
$$D_{r , \tau '} (v) \leq \, {1 \over |B_{r} \times (0, \tau ')|} \, \frac{| (0, \tau)|}{| (0, \tau)|} \, 2(\lambda - \frac{C2}{2})^{-1} (\frac{C_{HG}^2}{(r - r ')^2} + \frac{1}{\tau - \tau '}) \iint_{R} v^2  \ dxdt$$
$$\leq \, {1 \over |B_{r} \times (0, \tau )|} \, \frac{| (0, \tau)|}{| (0, \tau')|} \, 2(\lambda - \frac{C2}{2})^{-1} (\frac{C_{HG}^2}{(r - r ')^2} + \frac{1}{\tau - \tau '}) \iint_{R} v^2  \ dxdt$$

$$= \frac{ \tau}{\tau '} \, 2(\lambda - \frac{C2}{2})^{-1} (\frac{C_{HG}^2}{(r - r ')^2} + \frac{1}{\tau - \tau '}) \, H_{r ,\tau} (v)$$

$$\leq \, \frac{ \tau}{\tau '} \, 2(\lambda - \frac{C2}{2})^{-1} (\frac{C_{HG}^2}{(r - r ')^2} + \frac{1}{\tau - \tau '}) \, H_{r ,\tau} (v)$$
gives for some $C_4$ dependent on $\lambda$, $n$ and $R$
\begin{equation}
\label{localboundsonD}D_{r ', \tau '} (v) \leq \, C_4 \, H_{r ,\tau} (v)
\end{equation}

\noindent So substituting these into the results of the lemma above, we get

$$H_{r ', \tau '} (v^{k}) \leq  C_1 C_3 H_{r,\tau }(v)^{\frac{2}{n}} [C_4 H_{r, \tau }(v)+ C_2H_{r, \tau }(v)]  $$

$$=  C_1 C_3[C_4 + C_2 ] H_{r,\tau }(v)^{\frac{2}{n}}  H_{r, \tau }(v)$$

$$= \gamma \, H_{r , \tau }(v)^{\frac{2}{n} + 1}$$
\begin{equation}
= \gamma \, H_{r , \tau }(v)^k
\end{equation}
where 
$$
\label{gamma}\gamma = 4 C_S C_B (\frac{r}{r'})^Q [ (\lambda - \frac{C2}{2})^{-1} \frac{ \tau}{\tau '} (\frac{C_{HG}^2}{(r - r ')^2} + \frac{1}{\tau - \tau '}) +\frac{1}{2r^2} ] \bigl(  8|B_r|     (\frac{C_{HG}^2}{(r - r ')^2} + \frac{1}{\tau - \tau '}) \,     \bigr)^{\frac{2}{n}}.
$$

%$H_{r ' \tau '} (w^{2k}) \leq  \Bigl( 8 \lambda ^2 (\frac{C_{HG}^2}{(r - r ')^2} + \frac{1}{\tau - \tau '}) +1 \Bigr) H_{r  \tau } (w) \Bigl( 8 \lambda (\frac{C_{HG}^2}{(r - r ')^2} + \frac{1}{\tau - \tau '}) H_{r  \tau }(w) \Bigr) ^{\frac{2}{n}}$

%$H_{r ' \tau '} (w^{2k}) \leq  \Bigl( 8 \lambda ^2 (\frac{C_{HG}^2}{(r - r ')^2} + \frac{1}{\tau - \tau '}) +1 \Bigr)[8 \lambda (\frac{C_{HG}^2}{(r - r ')^2} + \frac{1}{\tau - \tau '})]^{\frac{2}{n}} H_{r  \tau } (w)   H_{r  \tau }(w) ^{\frac{2}{n}}$

%$H_{r ' \tau '} (w^{2k}) \leq  C(\lambda, R, n) \,  H_{r  \tau }(w) ^{\frac{2}{n} + 1} = C(\lambda, R, n) \,  H_{r  \tau }(w) ^k $

\noindent We can repeat this exponentiation by $k$ as in the parabolic case. Take
\begin{equation}
\label{vnu} v_{\nu} := v^{p_{\nu} /2}
\end{equation}
a positive subsolution in $R$ where $p_{\nu} = p_0 k^{\nu}$, $p_{\nu} \geq p_0 \geq p' >1$ and $k = 1 + 2/n$, $n \in \mathbb{N}$ is as above. 

Now we set up the domains for each of the iterations. Each iteration is defined in its own rectangle. 

\noindent Let 
\begin{equation}
\tau _{\nu} := \tau \frac{1 + \tau ' v}{1+ \tau v},
\end{equation} 
and
\begin{equation}
r _{\nu} := r \frac{1 + r ' v}{1 + r v}
\end{equation}
Let 
\begin{equation}
R_{\nu} = B_{r _{v}} \times (0, \tau _{\nu} )
\end{equation}

\noindent Set 
\begin{equation}
H_{\nu} := H_{r _{\nu} \tau _{\nu}} (v_{\nu}) = \frac{1}{|R_{\nu}|} \iint_{R_{\nu}} v^{p_{\nu}} dx dt.
\end{equation}
\noindent Now we will show the iteration. We notice that
\begin{equation}
H_{\nu + 1} = H_{r _{\nu +1} \tau _{\nu +1}} (v_{\nu +1}) = \frac{1}{|R_{\nu +1}|} \iint_{R_{\nu + 1}} v^{p_{\nu +1}} dx dt \leq \gamma H_{\nu} ^k
\end{equation}
where $\gamma$ depends on $R$, $n$, $\lambda$ only. 

Now we will use this fact to achieve a maximum. First iterate down to $H_0$. Take any $H_{\nu }$ and notice that 
$$H_{\nu} \leq \gamma _1 H^k_{\nu -1} \leq \gamma _2 ^{\nu} H^k _{\nu -1}$$ 
for some $\gamma_2 \geq \gamma^{\frac{1}{\nu}}_1.$
So $$H_{\nu } \leq \gamma _2 ^{\nu} H^k _{\nu -1} $$
Now use $(H_{\nu-1})^k$ to iterate again.
$$= \gamma _2 ^{\nu} (H_{\nu-1})^k \leq \gamma _2 ^{\nu } (\gamma_2 ^{\nu - 1} H^k_{\nu - 2})^k$$

$$= \gamma _2 ^{(\nu ) + (\nu -1)k }  H^{k^2}_{\nu - 2}$$

$$\leq \gamma _2 ^{(\nu ) + (\nu -1)k }  (\gamma _2 ^{(\nu - 2)}H^k_{\nu - 3})^{k^2}$$

$$= \gamma _2 ^{(\nu ) + (\nu -1)k + (\nu -2)k^2} H_{\nu - 3}^{k^3}$$

$$\leq \gamma _2 ^{(\nu) + (\nu -1)k + (\nu -2)k^2 + ... + (1)k^{\nu - 1}} H_{0}^{k^{\nu}}$$

$$\leq \gamma _3 ^{k^{\nu}} H_{0}^{k^{\nu}}$$ 
for some 
$$\gamma_3 \geq \gamma_2^{\frac{\nu + (\nu -1)k + ... + (1)k^{\nu -1}}{k^{\nu}}}$$
So we have 
$$H_{\nu} \leq \gamma_3^{k^{\nu}} H_{0}^{k^{\nu}}$$
which says that
$$ \frac{1}{|R_{\nu}|} \iint_{R_{\nu}} v^{p_0 k^{\nu} / 2 } dxdt \leq \gamma_3^{k^{\nu}} \bigl( \frac{1}{|R|} \iint_{R} v^{p_0 / 2 } dxdt \bigr)^{k^{\nu}} $$
Then taking $k^{\nu}$ roots: 
\begin{equation}\bigl( \frac{1}{|R_{\nu}|} \iint_{R_{\nu}} |v^{p_0 / 2 }|^{k^{\nu}} dxdt \bigr)^{1/k^{\nu}} \leq \gamma_3  \frac{1}{|R|} \iint_{R} v^{p_0 / 2 } dxdt. 
\end{equation} 
letting $\nu \rightarrow \infty$ we conclude:
$$ \underset{R'}{\max} \, u \, \leq \gamma \bigl( \frac{1}{|R|} \iint_{R} u^{p } dxdt \bigr)^{1/p}. $$
This proves Theorem \ref{NLmaxu}.

Similarly we are able to achieve the inequalities
\begin{equation}
\gamma \bigl( \frac{1}{|R|} \iint_{R} u^{-p } dxdt \bigr)^{1/-p} \, \leq \, \underset{R'}{\min} \, u 
\end{equation}
and the crucial bridge step
\noindent 
\begin{equation}
 \frac{1}{|R_{-}|} \iint_{R_{-}} v^{p_0 k^{\nu} / 2 } dxdt \leq \gamma \frac{1}{|R^{+}|} \iint_{R^{+}} v^{-p_0 k^{\nu} / 2 } dxdt 
\end{equation}
which are now seen to be independent of the nonlinear function $A_j(x, \hg u)$. So the iterations for each of these inequalities may be run in a similar way with the slight difference that the constants will be altered in the manner explained above.
 %-------------------------------------------------------------
% You cannot use the \chapter{} or \section{} commands because they put a header that the grad school doesn't allow
\begin{center} 
\textbf{2.5 Nonlinear Estimates on log(u) and H\"older continuity}
\end{center}
%-------------------------------------------------------------
 
Next we wish to obtain estimates on the function $v:= log(u)$ for a weak supersolutions $u$ to 
\begin{equation}\label{nlparab}\nlparab. \end{equation} 
If we are able to obtain estimates of the form
\begin{equation}
 {1 \over  |R^+|} \int_{R^+} \sqrt{(v-C(x,r))}  \ dxdt\leq N(v)
\end{equation} and 
\begin{equation}
{1 \over  |R^-|} \int_{R^-} \sqrt{(C(x,r)-v)}  \ dxdt\leq N(v)
\end{equation}
for some constant $N$ depending on the function $v$, some constant $C$ and rectangles $R^+$ and $R^-$ subsets of $R((x,t), r)$. Then we will be in a position to apply Aimar's \cite{Ai} result for functions satisfying a BMO property with a lag mapping used above and complete the proof of the Harnack inequality. 

To begin, we will first let $v=-log(u)$, $\phi = \psi^2 u^{-1}$, with $\psi \geq 0$    and of compact support on $B(x,2r)$. We will let 
$$\psi = \psi (x) = \prod _{\nu = 1}^n \chi _{\nu} (x_{\nu})$$ and  
$$
%\begin{equation*}
\chi _{\nu} (x_{\nu}) = \left\{
 \begin{array}{rl}
  1 & \text{if } |x_{\nu} - x_{\nu}^0| < r \\
   0 & \text{if } |x_{\nu} - x_{\nu}^0| \geq 2r \\
 \end{array} \right.
%\end{equation*}
$$

\noindent following \cite{CiL}. Here we use $(x_{1}^0, x_{1}^0, ...x_{n}^0)$ to be the center of the ball in space. Then we calculate
$$v_t = (-log(u))_t = {-u_t \over u}$$

$$X_j \phi = X_j (\psi^2 u^{-1}) = 2 \psi X_j \psi u^{-1} - \psi^2 u^{-2} X_j u$$

$$X_j (v) = X_j (-log(u)) = -X_j u (u^{-1}) $$
And with (\ref{nlparab}) we have

$$
\iint_{R} \psi^2 u^{-1} u_t  \ dxdt + 2\sum_{j=1}^{m} \iint_{R} A_j \psi X_j \psi u^{-1} \ dxdt - \sum_{j=1}^{m} \iint_{R} A_j \psi^2 u^{-2} X_j u \ dxdt \ \geq 0.
$$
Moving the second integrand to the RHS

$$
\iint_{R} \psi^2 u^{-1} u_t  \ dxdt  - \sum_{j=1}^{m} \iint_{R} A_j \psi^2 u^{-2} X_j u \ dxdt \ \geq -2\sum_{j=1}^{m} \iint_{R} A_j \psi X_j \psi u^{-1} \ dxdt
$$ 
which gives us by substituting $v_t$

$$
-\iint_{R} \psi^2 v_t  \, dxdt  - \sum_{j=1}^{m} \iint_{R} A_j \psi^2 u^{-2} X_j u \ dxdt \ \geq -2\sum_{j=1}^{m} \iint_{R} A_j \psi X_j \psi u^{-1} \ dxdt
$$
and hence

$$
\iint_{R} \psi^2 v_t   \ dxdt  + \sum_{j=1}^{m} \iint_{R} A_j \psi^2 u^{-2} X_j u \ dxdt \ \leq  \ 2\sum_{j=1}^{m} \iint_{R} A_j \psi X_j \psi u^{-1} \ dxdt
$$
Again we use the Cauchy and Young inequalities on the RHS to obtain

$$
\iint_{R} \psi^2 v_t  \ dxdt   + \sum_{j=1}^{m} \iint_{R} A_j \psi^2 u^{-2} X_j u \ dxdt \
$$

$$
 \leq  \ {1 \over 2} \sum_{j=1}^{m} \iint_{R} |A_j|^2 u^{-2} \psi^2 \ dxdt + 2 \iint_{R} |\hg \psi |^2 \ dxdt
$$
Now use the boundedness of $A_j(x, \hg u)$ from the $(S)$ conditions on the second integrand on the LHS and on the first integrand of the RHS
$$
\iint_{R} \psi^2 v_t  \ dxdt   + \lambda \iint_{R} \psi^2 |\hg u|^2 u^{-2}  \ dxdt \
$$

$$
 \leq  \ { \Lambda \over 2}  \iint_{R} |\hg u|^2 u^{-2} \psi^2  \ dxdt+ 2 \iint_{R} |\hg \psi |^2 \ dxdt
$$
At this point, notice that we have $X_j (v) = -X_j u (u^{-1}) $ from above which gives us $| \hg v| =  |\hg u| u^{-1} $.

$$
\iint_{R} \psi^2 v_t   \ dxdt  + \lambda \iint_{R} \psi^2 |\hg v|^2  \ dxdt \ \leq  \ { \Lambda \over 2}  \iint_{R} |\hg v|^2   \psi^2  \ dxdt+ 2 \iint_{R} |\hg \psi |^2 \ dxdt
$$
Gathering like terms
\begin{equation}
\iint_{R} \psi^2 v_t   \ dxdt  + (\lambda -{ \Lambda \over 2} ) \iint_{R} \psi^2 |\hg v|^2  \ dxdt \ \leq  \ 2 \iint_{R} |\hg \psi |^2 \ dxdt
\end{equation}
Integrate from $t_1$ to $t_2$ and observe that $\psi = 1$ on $B(x,r)$ is independent of $t$.

$$
\int_{B(x,r)}  v  \ dx|_{t_1}^{t_2}  + (\lambda -{ \Lambda \over 2} ) \int_{t_1}^{t_2} \int_{B(x,r)}   |\hg v|^2  \ dxdt \ \leq  \ 2  (t_2 - t_1) \int_{B(x,r)} |\hg \psi |^2 \ dx
$$
Now we recall the Poincar\'e lemma and the definition of $V_B(t)= {1 \over |B|} \int v dx$.

$$
[V_B(t_2) - V_B(t_1)]|B(x,r)|   + C_P r^2(\lambda -{ \Lambda \over 2} ) \int_{t_1}^{t_2} \int_{B(x,r)}   (v - V_B)^2  \ dxdt $$
$$ \leq  \ 2  (t_2 - t_1) \int_{B(x,r)} |\hg \psi |^2 \ dx
$$
Divide by $(t_2 - t_1)$ 

$${[V_B(t_2) - V_B(t_1)] \over (t_2 - t_1)}|B(x,r)|   + C_P r^2(\lambda -{ \Lambda \over 2} ){1 \over (t_2 - t_1)} \int_{t_1}^{t_2} \int_{B(x,r)}   (v - V_B)^2  \ dxdt \ $$
$$\leq  \ 2  \int_{B(x,r)} |\hg \psi |^2 \ dx$$
Take $t_2 \rightarrow t_1$

$$
{dV_B \over dt}|B(x,r)|   + C_P r^2(\lambda -{ \Lambda \over 2} ) \int_{B(x,r)}   (v - V_B)^2  \ dxdt \ \leq  \ 2  \int_{B(x,r)} |\hg \psi |^2 \ dx.
$$
Divide through by  $|B(x,r)|$ and use the bounds on $|\hg \psi| \leq 1/r$

$$
{dV_B \over dt}    + {C_P r^2 \over |B(x,r)|}(\lambda -{ \Lambda \over 2} ) \int_{B(x,r)}   (v - V_B)^2  \ dxdt \ \leq  \ {2 \over |B(x,r)|r^2}  |B(x,r)|.
$$ 
We note that this inequality is the same form as (\ref{dV_plus_poincare_LEQ_const}) and that the proof is identical for both positive and negative $t$. This shows that $v = - log(u)$ meets the BMO condition from Aimar \cite{Ai}. So we may apply his results to finish the bridge step. Thus we have proved the Harnack inequality for positive, weak solutions to $$\nlparab.$$ With the Harnack inequality in place, the H\"older continuity with respect to the CC metric follows by the same argument as before following the pattern of J. Moser \cite{Mo1}.

\newpage

%You have to type the bibliography here instead of calling it from a separate file because it has to be single spaced within the entries and
%double spaced between the entries.  As far as I know, there is no way to have that formatting when you call the bibliography from a .bib file.
%However, the bibliography cannot have the Bibliography heading at the top of the pages, so while you have to put the bibliography here so that 
%the references are called properly within the dissertation, you have to type them in an entirely separate document for printing.
\singlespace


\begin{thebibliography}{99}
\bibliographystyle{alpha} 
%\bibitem[A1]{andrews:cont}
%ANDREWS, B., Contraction of hypersurfaces in Euclidean space.  \emph{Calc. Var. %Partial Differential Equations 2} (1994), 151-171.\\
 

\bibitem {Ai}AIMAR, H.,\! Elliptic\! and Parabolic\! BMO and Harnack's Inequality. \emph{Transactions American Mathematical Society Vol. 306, Number 1} (March 1988), 265-276.  
  
\bibitem {AS}ARONSON, D., and SERRIN, J., Local behavior of solutions of quasilinear parabolic \ equations.  \emph{Archive \ for \ Rat.\  Mech. \ and \ Anal. \ Volume 25, Number 2} (1967), 81-122. 
   
\bibitem {Bo}BONY, J., \ Principe du maximum, in\a'{e}galit\a'{e} de Harnack \ et unicit\a'{e} du probl\a'{e}me de Cauchy pour les op\a'{e}rateurs elliptiques d\a'{e}g\a'{e}n\a'{e}r\a'{e}s, \emph{Annales de l'institut Fourier Volume 19, Number 1} (1969), 277-304.   
 
\bibitem {Br}BRAMANTI, M., and UGUZZONI, F., Heat kernels for non-divergence operators of H\"{o}rmander type \emph{C. R. Acad. Sci. Paris Ser. I Volume 343} (2006).   
\bibitem {CDG1} CAPOGNA, L., DANIELLI, D., and GAROFALO, N., Capacitary estimates and the local behavior of solutions on nonlinear subelliptic equations. \emph{Am. J. Math. 118, 6} (1996), 1153-1196. 
\bibitem {CDG2}CAPOGNA, L., DANIELLI, D., and GAROFALO, N., The geometric Sobolev embedding for vector fields and the isoperimetric inequality. \emph{Comm. in Analysis and Geometry 2} (1994), 203-215. 

\bibitem {CDG3}CAPOGNA, L., DANIELLI, D., and GAROFALO, N., An embedding the-  orem and the Harnack inequality for nonlinear subelliptic equations. \emph{Comm. Partial Diff. Eq. 18} (1993), 1765-1794.



\bibitem {CaD}CAPOGNA, L., DANIELLI, D., PAULS, S., and TYSON, J., \emph{An introduction to the Heisenberg group and the sub-Riemannian isoperimetric problem}, vol. 259 of \emph{Progress in Mathematics}. Birkh\"{a}user, 2007.

\bibitem {ChH} CHENG, J., \ HWANG, J., \ MALCHIODI, A. and YANG, P., Minimal surfaces in pseudohermitian geometry. Ann. Scuola Norm. Sup. Pisa Volume 5 Number 4 , (2005), 129-177.

\bibitem {CiL}CITTI, G., GAROFALO, N., and\, LANCONELLI, E. \,\,Harnack's\,\, inequality for sum of squares of vector\, fields \,plus\, a potential.\, \emph{Amer.\, Jour.\, of \,Math. 115}, 3 (1993), 699-734.



\bibitem {FG}FABES, E. and\, GAROFALO N., Parabolic B.M.O. and Harnack's Inequality. \emph{Proceedings of the American Mathematical Society, Vol. 95, No. 1} (September 1985), 63-69.


\bibitem {FrL}FRANCHI, B., and LANCONELLI, E., \ H\"{o}lder regularity theorem for a class  of \ nonuniformly elliptic \ operators with measurable coefficients. \emph{Ann. Scuola norm. sup. pisa. 10} (1983), 523-541.


\bibitem {GaN}GAROFALO, N. and\, NHIEU, D. M. Isoperimetric and Sobolev inequalities for\,\, Carnot-Carath\a'{e}adory\,\, metrics\,\, and \,the\, existence \,of\, minimal\, surfaces. \emph{Comm. in Pure and Appl. Math. 49} (1996), 1081-1144.


\bibitem {GiT}GILBARG, D., and TRUDINGER, N., \emph{Elliptic Partial Differential Equations of Second Order.} Springer-Verlag, Berlin, 2001.


\bibitem {Gr}GRIGORY'AN, A., The heat equation on non-compact Riemannian Manifolds. \emph{Math. USSR Sb 72} (1992), 47-77.


\bibitem {Had}HADAMARD, J., Extension a l'\a'{e}quation de la chaleur d'un th\a'{e}oreme de A. Harnack. \emph{Rendiconti del Circolo Matematico di Palermo  3} Number 3, (1954), 337-346.

\bibitem {Hal}HALLER, E. Comparison principles for fully-nonlinear parabolic equations and regularity theory for weak solutions of parabolic systems in Carnot groups. PhD Dissertation, University of Arkansas - Fayetteville (2007).

\bibitem {Ho}H\"{O}RMANDER,\, L.\, Hypoelliptic\,\, second\,\, order\,\, differential\,\, equations. \emph{Acta Math.}, 119 (1967), 147-171.


\bibitem {Jer}JERISON, D.\! The\! Poincar\a'{e} inequality\! for vector fields satisfying H\"{o}rmander's condition. \emph{Duke Math. J. 53} (1986), 503-523.


\bibitem {JeS}JERISON, D., and SANCHEZ-CALLE, A., Sub-elliptic second order differential operators. \emph{Lectures Notes in Mathematics 1277} (1985), 46-77. 


\bibitem {KrS}KRYLOV, N., \ and SAFONOV, \ M., \ A certain property \ of solutions of parabolic \ equations with measurable coefficients. \ \emph{Math. USSR \ Izvestiya Volume 16, Number 1} (1980), 151-164. 


\bibitem {KS}KUSUOKA, S., and STROOCK, D., \ \ Long time estimates for the heat kernel associated with a uniformly subelliptic symmetric second order operator.  \emph{Ann. of Math. Volume 2, Number 127} (1988), 165-189. 

\bibitem{L}LIEBERMAN, G., \emph{Second Order Parabolic Differential Equations}. World Scientific Publishing Co., River Edge, NJ, 1998.


\bibitem {Ma}MANFREDINI, \ M., \ A \ note \ on the Poincare inequality \ for Lipschitz vector ields of step two. \emph{Proc. Amer. Math. Soc. 138}  (2010), 567-575.


\bibitem {Mo1}MOSER, J., \!A Harnack inequality for parabolic differential equations. \emph{Comm. Pure Appl. Math. 24} (1964), 101-134. 


\bibitem {Mo2}MOSER, J., \!On Harnack's theorem for elliptic differential equations. \emph{Comm. Pure Appl. Math. 14} (1961), 577-591.


\bibitem {Mo3}MOSER, J., \ On \ a pointwise \ estimate \ for \ parabolic \ differential equations. \emph{Comm Pure Appl. Math. Volume 24, Number 5} (1971), 727-740.

\bibitem {Nag}NAGEL, A., and STEIN, E., Differentiable control metrics and scaled bump functions. \emph{Journal of Differential Geometry 57} (2001), 465-492.

\bibitem {NaW}NAGEL, A., STEIN, E., and WAINGER, S. Balls\, and\, metrics\, defined by vector fields I: Basic properties. \emph{Acta Math. 155} (1985), 103-147.


\bibitem {Nas}NASH, J., \!Continuity of solutions of parabolic and elliptic equations. \emph{Amer. J. Math. 80} (1958), 931-954.

\bibitem {Pa} PAULS, S., Minimal \ Surfaces in the \ Heisenberg Group. \emph{Geometriae Dedicata Volume 104, Number 1} (2004), 201-231.

\bibitem {SC}SALOFF-COSTE, L., \!Parabolic Harnack inequality for divergence form second order differential operators. \emph{Potential Anal. 4} (1995) 429-467. 


\bibitem {Se}SERRIN, J., \!A Harnack inequality for nonlinear equations. \emph{Bull. Amer. Math. Soc. Volume 69, Number 4} (1963), 481-486.
\bibitem {Sh}SHORES, E. Hypoellipticity \,for \,linear\, degenerate \,elliptic systems in Carnot groups and applications. e-print. arXiv.org: math.AP/0502569.
\bibitem {Sh1}SHORES, E. Regularity \,theory for \,weak solutions\, of systems in Carnot groups. PhD Dissertation, University of Arkansas - Fayetteville (2005).


\bibitem {Tr1}TRUDINGER, N., \!On Harnack-type inequalities and their applications to quasilinear elliptic equations. \emph{Comm. Pure Appl. Math. 20} (1967), 721-747.


\bibitem {Tr2}TRUDINGER, N., \!Some existence theorems for quasi-linear, non-uniformly elliptic equations in divergence form. \emph{J. Math. Mech. Anal. 18} (1968), 909-919. 


\bibitem {Tr3}TRUDINGER, N., Linear Elliptic operators with measurable coefficients. \emph{Ann. Scuola Norm. Sup. Pisa Volume 3, Number 27} (1973), 265-308.


\bibitem {UL1}UGUZZONI, F., BONFIGLIOLI, A., and LANCONELLI, E., Uniform Gaussian estimates of the fundamental solutions for heat operators on Carnot groups. \emph{Adv.  Differential Equations 7} (2002), 1153-1192.


\bibitem {UL2}UGUZZONI, F., BONFIGLIOLI, A., and LANCONELLI, E.,  Fundamental solutions for non-divergence form operators on stratified groups. \emph{Trans.  Amer. Math. Soc. 356} (2004), 2709-2737.


\bibitem {VSC}VAROPOULOS, N., SALOFF-COSTE, L., and COULHON, T., \emph{Ananlysis and Geometry on Groups.} Cambridge University Press, Cambridge, UK, 1992.  
 
 
 
 
 
 

\end{thebibliography}
\end{document}